\documentclass[12pt]{article}
\usepackage{amsmath,amsfonts,amsthm,amssymb}
\usepackage{typearea}
\newcommand{\stsets}[1]{\mathbb{#1}}
\newcommand{\R}{\stsets{R}}

\theoremstyle{definition}
\newtheorem{definition}{Definition}
\theoremstyle{remark}
\newtheorem{remark}[definition]{Remark}
\newtheorem{example}[definition]{Example}

\newtheoremstyle{mytheorem}{0.5cm}{0.2cm}{\slshape}{ }{\bfseries}{.}{ }{}
\theoremstyle{mytheorem}
\newtheorem{Th}[definition]{Theorem}
\newtheorem{lemma}[definition]{Lemma}
\newtheorem{Cor}[definition]{Corollary}

\renewcommand{\P}{\mathbf{P}}
\renewcommand{\Re}{\mathrm{Re\,}}

\DeclareMathOperator{\E}{{\bf E}}
\DeclareMathOperator{\var}{{\bf var}}
\DeclareMathOperator{\supp}{supp}
\DeclareMathOperator{\co}{co}
\DeclareMathOperator{\one}{{\mathbf{1}}}

\newcommand{\restr}[1][r]{|^{#1}}

\DeclareMathOperator*{\ssum}{{\textstyle \sum}} 
\newcommand{\comp}{c}

\renewcommand{\epsilon}{\varepsilon}
\renewcommand{\vec}[1]{\mathbf{#1}}
\renewcommand{\phi}{\varphi}
\newcommand{\ti}{\to\infty}
\newcommand{\ssp}{\hspace{2pt}}
\newcommand{\seg}{see, \hbox{e.\ssp g.,}\ }
\newcommand{\ie}{\hbox{i.\ssp e.}\ }
\newcommand{\eg}{\hbox{e.\ssp g.,}\ }
\newcommand{\cf}{\hbox{c.\ssp f.}\ }
\newcommand{\iid}{\hbox{i.\ssp i.\ssp d.}\ }
\newcommand{\as}{\hbox{a.\ssp s.}\ }

\newcommand{\sB}{\mathcal{B}}
\newcommand{\sK}{\mathcal{K}}
\newcommand{\sM}{\mathcal{M}}

\newcommand{\sF}{\mathcal{F}}
\newcommand{\DD}{\mathbb{D}}
\newcommand{\EE}{\mathbb{E}}
\newcommand{\GG}{\mathbb{G}}
\newcommand{\BB}{\mathbb{B}}
\newcommand{\CC}{\mathbb{C}}
\newcommand{\KK}{\mathbb{K}}
\newcommand{\FF}{\mathbb{F}}
\newcommand{\TT}{\mathbb{T}}
\renewcommand{\SS}{\mathbb{S}}
\newcommand{\sas}{S\alpha S}

\newcommand{\deq}{\overset{\scriptscriptstyle\mathcal{D}}{=}}
\newcommand{\vto}{\overset{{\mathrm{v}}\ }{\rightarrow}}
\newcommand{\wto}{\Rightarrow}

\newcommand{\eps}{\varepsilon}
\renewcommand{\EE}{\KK'}
\newcommand{\cset}{B^r}
\newcommand{\neutral}{\mathbf{e}}
\newcommand{\zero}{\mathbf{0}}
\newcommand{\KDU}{\tilde\KK}
\newcommand{\ddual}{\sharp}
\newcommand{\KDD}{\KK^\ddual}
\newcommand{\sMzero}{\mathcal{M}_\zero}
\newcommand{\setr}[1][r]{A_{#1}}

\newcommand{\cl}{\mathrm{cl}}

\newcommand{\xir}{\xi^{(r)}}

\numberwithin{equation}{section}
\numberwithin{definition}{section}

\begin{document}

\bibliographystyle{acmtrans-ims}

\thispagestyle{empty}
\vspace*{1.5cm}
\begin{center}
  
  {\bfseries\LARGE Strictly stable distributions\\[3mm] on convex
    cones\footnote{Supported by the British Council/Alliance Fran\c
      caise
      grant {\em Geometric interpretation of stable laws}.} }\\[10mm]
  
  \textsc{By Youri Davydov, Ilya Molchanov\footnote{Supported 
      by the Swiss National Science
      Foundation Grant 200021--100248.} and Sergei Zuyev}

  \textit{University  Lille 1, University of Berne and University of
    Strathclyde}
\end{center}

\begin{quotation}
  \small Using the LePage representation, a symmetric $\alpha$-stable
  random element in Banach space $\BB$ with $\alpha\in(0,2)$ can be
  represented as a sum of points of a Poisson process in $\BB$. This
  point process is union-stable, \ie the union of its two independent
  copies coincides in distribution with the rescaled original point
  process.  This shows that the classical definition of stable random
  elements is closely related to the union-stability property of point
  processes.
  
  These concepts makes sense in any convex cone, \ie in a semigroup
  equipped with multiplication by numbers, and lead to a construction
  of stable laws in general cones by means of the LePage series. We
  prove that random samples (or binomial point processes) in rather
  general cones converge in distribution in the vague topology to the
  union-stable Poisson point process. This convergence holds also in a
  stronger topology, which implies that the sums of points converge in
  distribution to the sum of points of the union-stable point process.
  Since the latter corresponds to a stable law, this yields a limit
  theorem for normalised sums of random elements with $\alpha$-stable
  limit for $\alpha\in(0,1)$.
  
  By using the technique of harmonic analysis on semigroups we
  characterise distributions of $\alpha$-stable random elements and
  show how possible values of the characteristic exponent $\alpha$
  relate to the properties of the semigroup and the corresponding
  scaling operation, in particular, their distributivity properties.
  It is shown that several conditions imply that a stable random
  element admits the LePage representation.  The approach developed in
  the paper not only makes it possible to handle stable distributions
  in rather general cones (like spaces of sets or measures), but also
  provides an alternative way to prove classical limit theorems and
  deduce the LePage representation for strictly stable random vectors
  in Banach spaces.

  \renewcommand{\thefootnote}{}


\end{quotation}

\tableofcontents

\section{Introduction}
\label{sec:introduction}

Stability of random elements is one of the basic concepts in
probability theory. A random vector $\xi$ with values in a Banach
space $\BB$ has a strictly stable distribution with characteristic
exponent $\alpha\neq0$ (notation $\sas$) if, for all $a,b>0$,
\begin{equation}
  \label{eq:a1alph-a+b1-}
  a^{1/\alpha}\xi_1+b^{1/\alpha}\xi_2\deq
  (a+b)^{1/\alpha}\xi\,, 
\end{equation}
where $\xi_1,\xi_2$ are independent copies of $\xi$, and $\deq$
denotes equality in distribution.  This stability concept (in a more
general form) was introduced by Paul L\'evy and thereafter has been
actively studied in relation to limit theorems for sums of random
variables, \seg \cite{haz:sieb01,sam:taq94,zol86} for the
finite-dimensional case and \cite{ara:gin,sam:taq94} for random
elements in Banach spaces. The following basic results are available:
\begin{itemize}
\item a complete characterisation of $\sas$ random elements in terms
  of characteristic functionals;
\item a complete description in terms of LePage (or
  Khinchin--L\'evy--LePage) expansions;
\item a complete description of the domains of attraction using tail
  behaviour. 
\end{itemize}
Recall that the LePage series representation of a $\sas$ vector
$\xi\in\R^d$ for $\alpha\in(0,2)$ says
\begin{equation}
  \label{eq:lpg}
  \xi \deq \sum_{k=1}^\infty \Gamma_k^{-1/\alpha}\eps_k\,,
\end{equation}
where $\Gamma_1,\Gamma_2,\ldots$ are the successive times of jumps of
a homogeneous Poisson process on the positive half-line, and
$\eps_1,\eps_2,\ldots$ are \iid unit random vectors independent of the
$\Gamma$'s. Note that in \cite{sam:taq94} the notation $\sas$ stands
for \emph{symmetric} $\alpha$-stable. We prefer to use the same
notation for \emph{strictly} $\alpha$-stable random elements, since
for semigroup-valued random elements the symmetry requirement may be
too restrictive and does not necessarily imply the strict stability.

Definition (\ref{eq:a1alph-a+b1-}) of stability makes sense in any
space, where addition of elements and multiplication by positive
scalars are defined, \ie in any convex cone. So far this case has been
thoroughly investigated for the cone of compact convex subsets of a
separable Banach space with the topology generated by the Hausdorff
metric and the main operation being the Minkowski (elementwise)
addition, \seg \cite{dav:pau:rac00,gin:hah85b,gin:h:z}.  In this case
the principal results (LePage representation, domains of attractions,
etc.) are completely analogous to those well-known for general Banach
spaces. The only exception is the lack of non-trivial strictly stable
distributions for $\alpha>1$, which is explained by impossibility to
centre sums of sets and the fact that the exact dual operation to the
Minkowski addition cannot be defined.

A wealth of information about stable probability measures on Euclidean
spaces and locally compact groups can be found in \cite{haz:sieb01},
see also \cite{haz86} and \cite{hey04}. Strictly stable distributions
are necessarily infinitely divisible. Infinite divisibility of
random objects in positive convex cones was studied in \cite{jon98}.
Infinite divisible elements in semigroups were studied in
\cite{rus88a,rus88b} and are comprehensively covered in the monograph
\cite{hoeg:muk95}. As we show in this paper, the studies of stability
of random elements in semigroups bring into the play further
properties of semigroups, in particular, the relationships between the
neutral element and the origin, distributivity laws and the metric
structure on semigroups.

Max-stable random variables appear if the addition in
(\ref{eq:a1alph-a+b1-}) is replaced by the maximum operation and
min-stable random variables in case of the minimum operation
\cite{res87}.  A max-stable $\sas$ laws exists for every positive
$\alpha$, while min-stable laws may have every $\alpha<0$ as a
possible value for the characteristic exponent.  On the other hand, it
is well known that possible values of $\alpha$ for a Banach
space-valued random element fill $(0,2]$.  One of the aims of this
paper is to identify basic algebraic and topological properties of the
carrier space that control the range of the characteristic exponent
$\alpha$ of strictly stable random elements.  We also give answers to
the questions that concern series representations and domains of
attraction of stable laws. The key idea is the relationship between
the stability concept formulated in (\ref{eq:a1alph-a+b1-}) and the
concept of stability of point processes (and random sets) with respect
to the union \cite{ma,mo93l}.  The obtained results cover not only the
classical cases of linear spaces and extremes of random variables.
They provide a unified framework for considering additive- and
max-stable laws for random variables as special cases of
semigroup-valued random elements. This framework, in particular,
includes random closed sets stable with respect to Minkowski addition
or union operations and random measures stable with respect to
addition or convolution operations.

The content of this paper can be outlined as follows.
Section~\ref{sec:convex-cones} introduces the main algebraic concepts:
the convex cone, the neutral element and the origin, and some metric
properties that are important in the sequel.

Section~\ref{sec:lepage-series-cone} shows that the classical
stability concept (\ref{eq:a1alph-a+b1-}) is, in a sense, secondary to
the union-stability.  The LePage representation (\ref{eq:lpg}) of
stable laws provides an expression of a $\sas$ random vector $\xi$ as
the sum of points of a Poisson point process $\Pi_\alpha$ with support
points
\begin{displaymath}
  \supp\Pi_\alpha=\{\Gamma_k^{-1/\alpha}\eps_k,\; k\geq 1\}\,,
\end{displaymath}
where $\alpha\in(0,2)$.  The distribution of the random set
$\kappa_\alpha=\supp\Pi_\alpha$ is union-stable, \ie
\begin{equation}
  \label{eq:a1alph-b1alph-a+b1}
  a^{1/\alpha}\kappa'_\alpha\,\cup\, b^{1/\alpha}\kappa''_\alpha\deq
  (a+b)^{1/\alpha}\kappa_\alpha\,,
\end{equation}
where $\kappa'_\alpha$ and $\kappa''_\alpha$ are independent copies of
$\kappa_\alpha$. Distributions of union-stable random closed sets in
$\R^d$ have been completely characterised in \cite{mo93l}, see also
\cite[Ch.~4]{mo1}. Because $\kappa'_\alpha$ and $\kappa''_\alpha$
possess common points with probability zero,
(\ref{eq:a1alph-b1alph-a+b1}) immediately translates into
(\ref{eq:a1alph-a+b1-}) by taking the sums over the support points.
It should be noted that (\ref{eq:a1alph-b1alph-a+b1}) makes sense for
all $\alpha\neq0$ (see \cite{mo93l}), while (\ref{eq:a1alph-a+b1-}) in
a Banach space holds for $\alpha\in(0,2]$ only. One of the reasons for
this is that the series in (\ref{eq:lpg}) might diverge.  In
Section~\ref{sec:lepage-series} it is shown that the LePage series
absolutely converges for all $\alpha\in(0,1)$ if the semigroup
possesses a sub-invariant norm. This representation makes it possible
to define L\'evy processes with values in a cone, see
Section~\ref{sec:gener-lepage-seri}.

It is well known that, under a regular variation type assumption, a
normalised random sample (or binomial point process) converges in
distribution to the Poisson point process $\Pi_\alpha$, \seg
\cite{res87}. This result is generalised for point processes in Polish
spaces.  Furthermore, we show that in case $\alpha\in(0,1)$ this
convergence holds in a stronger topology that ensures the convergence
of sums of points from point processes. This complements the result of
\cite{dav:egor05}, where this type of convergence was studied for
point processes in $\R^d$.  From this fact we derive that normalised
sums of random elements converge in distribution to the LePage
representation of the corresponding stable law, \ie a limit theorem
for normalised samples in cones. This also yields a new proof of the
limit theorem in Banach spaces with $\sas$ limits for
$\alpha\in(0,1)$. These results are described in
Section~\ref{sec:conv-stable-laws-1}.

Section~\ref{sec:distr-stable-rand} explores the distributions of
stable random elements on semigroups, including the ranges of the
stability parameter. We first define the Laplace transform of a random
element as a functional that acts on the family of the characters, and
confirm its uniqueness. The infinite divisibility property implies
that the Laplace transform has an exponential form. The main result
establishes the equivalence between the stability property and the
homogeneity of the corresponding Laplace exponent.  Further we
describe several essential properties of cones and semigroups that
have a particular bearing in view of the properties of stable
distributions. Among these properties the most important are the
distributivity properties of the multiplication by numbers and the
relationship between the neutral element and the origin in a
semigroup. The range of possible parameters for the stable law
provides a new characteristic of a general cone. In particular, we
describe the cases when the stability parameter $\alpha$ belongs to
$(0,2]$ as in the conventional case of linear spaces and when $\alpha$
is an arbitrary positive or arbitrary negative number.

Section~\ref{sec:integr-repr-stable} exploits techniques from
harmonic analysis on semigroups, in particular, the representations of
infinite divisible and negative definite functions in view of
characterising the Laplace exponents of $\sas$ random elements. In
particular, we show that the corresponding L\'evy measure is
homogeneous, and characterise other ingredients of the integral
representations: the linear functional and the quadratic form.

Finally, Section~\ref{sec:lepage-seri-repr} aims to show that strictly
stable random elements in a rather general cone admit the LePage
representation. We address this question by comparing the integral
representations of Laplace exponents with the formula for the
probability generating functional of a stable Poisson process.  It is
shown that under rather weak conditions every $\sas$ random element
can be realised as the LePage series that corresponds to a Poisson
process on the second dual semigroup.  Its intensity measure is the
L\'evy measure of the corresponding $\sas$ random element. The key
issue here is to show that the L\'evy measure is actually supported by
the semigroup itself, which leads to the ``conventional'' LePage
series similar to (\ref{eq:lpg}). Apart from the proof of the LePage
representation in rather general semigroups, this also yields a new
proof of the LePage representation for $\sas$ random elements with
$\alpha\in(0,1)$ in reflexive Banach spaces.

To summarise, we first show that the sum of points of a union-stable
Poisson point process follows $\sas$ law, then demonstrate that
convergence of point processes yields a limit theorem with $\sas$
limiting distribution, and finally prove that in a rather general case
any $\sas$ random element with $\alpha\in(0,1)$ can be represented as
a sum of points of a Poisson process.

Section~\ref{sec:examples} describes a variety of examples drawing
analogy or contrasting with classical Euclidean or Banach space valued
$\sas$ vectors.  Clearly, if a cone is embeddable in a Banach space,
then it is possible to use the results already available for stable
distributions, see \cite{sam:taq94}. However even in this case we come
up with new proofs that further our understanding of stable laws in
linear spaces.

\section{Convex cones}
\label{sec:convex-cones}

\subsection{Basic definitions}
\label{sec:basic-definitions}

Here we summarise several basic definitions related to convex cones
and semigroups.

\begin{definition}
  An abelian \emph{topological semigroup} is a topological space $\KK$
  equipped with a commutative and associative continuous binary
  operation $+$.  It is assumed that $\KK$ possesses the \emph{neutral
    element} $\neutral$ satisfying $x+\neutral=x$ for every $x\in\KK$.
\end{definition}

Consider a family of continuous automorphisms $D_a:\KK\to \KK$ indexed
by positive real numbers $a>0$. Assume that $D_1$ is the identical map
and that $D_aD_bx=D_{ab}x$ for all $a,b>0$ and $x\in\KK$. In
\cite{rat:sch77} such $\KK$ is called an abelian semigroup over the
operator domain $(0,\infty)$.  The result of applying $D_a$ to $x\in
\KK$ can be understood as the \emph{multiplication} of $x$ by $a$ that
yields the following equivalent reformulation of the properties of
$D$.

\begin{definition}
  \label{def:cone}
  A \emph{convex cone} is an abelian topological semigroup $\KK$ being
  a metrisable Polish (complete separable) space with a continuous
  operation $(x,a)\mapsto ax$ of multiplication by positive scalars
  for $x\in\KK$ and $a>0$ so that the following conditions are
  satisfied:
  \begin{align}
    \label{cone1}
    a(x+y)&=ax+ay\,,\quad a>0,\; x,y\in\KK,\\
    \label{cone2}
    a(bx)&=(ab)x\,,\quad a,b>0,\; x\in\KK,\\
    \label{cone3}
    1x&=x\,,\quad x\in\KK\,,\\
    \label{cone4}
    a\neutral &=\neutral\,, \quad a>0\,.
  \end{align}
  $\KK$ is called a \emph{pointed cone} if there is a unique
  element $\zero$ called the \emph{origin} such that $ax\to\zero$ as
  $a\downarrow0$ for any $x\in\KK\setminus\{\neutral\}$.
\end{definition}

It should be emphasised that we do not always require the following
distributivity condition
\begin{equation}
  \label{eq:cone4}
  (a+b)x=ax+bx\,,\quad a,b>0,\; x\in\KK\,.
\end{equation}
We often call (\ref{eq:cone4}) the \emph{second distributivity law}.
Although this condition is typically imposed in the literature on
cones (\seg \cite{keim:rot92}), this law essentially restricts the
family of examples, \eg it is not satisfied for the cone of compact
(not necessarily convex) subsets of a Banach space (with Minkowski
addition) or on $\R_+=[0,\infty)$ with the maximum operation. Note
that the condition (\ref{eq:cone4}) is not natural if the
multiplication is generated by a family of automorphisms $D_a$, $a>0$,
as described above. In view of this, $nx$ means $D_n x$ or $x$
\emph{multiplied} by $n$ and \emph{not} the sum of $n$ identical
summands being $x$.

If the addition operation on $\KK$ is a group operation, then we
simply say that $\KK$ is a \emph{group}. If also the second distributivity
law holds, then $\KK$ satisfies the conventional axioms of a linear
space.  However, even if $\KK$ is a group, the second distributivity
law does not have to hold, \eg if $\KK=\R$ with the usual addition and
the multiplication defined as $D_ax=a^\beta x$ with $\beta\neq1$, see
also Example~\ref{ex:rescaled-functins}.

\subsection{Origin and neutral element}
\label{sec:orig-neutr-elem}

Unless stated otherwise we always assume that $\KK$ is a pointed cone,
\ie it possesses the origin.  The neutral element of $\KK$ does not
necessarily coincide with the origin. For instance, if $\KK$ is the
semigroup of compact sets in $\R^d$ with the union operation and the
conventional multiplication by numbers, then $\zero=\{0\}$, while the
neutral element $\neutral$ for the union operation is the empty set.
In many other cases the neutral element does coincide with the origin,
\eg if $\KK$ is a linear space.  Note that the definition of the
origin implies that $\zero+\zero=\zero$, and $\zero$ and $\neutral$
are the only elements of $\KK$ satisfying $ax=x$ for all $a>0$.

\begin{lemma}
  \label{lemma:n=0}
  Let $\KK$ be a pointed cone.
  \begin{description}
  \item[\textsl{(i)}] If the second distributivity law holds, then
    $\neutral=\zero$.
  \item[\textsl{(ii)}] If there exists $x\neq\neutral$ which possesses an
    inverse $(-x)$, \ie $x+(-x)=\neutral$, then $\neutral=\zero$.
  \end{description}  
\end{lemma}
\begin{proof}
  \textsl{(i)} By~(\ref{eq:cone4}), any $x\neq \neutral$ can be
  decomposed as
  \begin{displaymath}
    x=\frac{n-1}{n}x+\frac{1}{n}x\,,\quad n\geq 1\,.
  \end{displaymath}
  By letting $n\ti$ and using the continuity of multiplication we
  arrive at $x=x+\zero$. Thus, $\zero=\neutral$ by the uniqueness of
  the neutral element.

  \textsl{(ii)} By (\ref{cone1}) and (\ref{cone4}),
  $n^{-1}x+n^{-1}(-x)=\neutral$. The left-hand side converges as $n\ti$
  to $\zero+\zero$, whence $\zero=\neutral$.
\end{proof}

In particular, Lemma~\ref{lemma:n=0}\textsl{(ii)} implies that
$\zero=\neutral$ if $\KK$ is a group.

\begin{definition}
  An element $z\in\KK$ is called $\alpha$-stable with $\alpha\neq0$,
  if 
  \begin{equation}
    \label{eq:selfsas}
    a^{1/\alpha}z + b^{1/\alpha} z=(a+b)^{1/\alpha}z
  \end{equation}
  for all $a,b>0$.
\end{definition}

Throughout this paper $\KK(\alpha)$ denotes the set of $\alpha$-stable
elements of $\KK$.  Clearly, $\neutral,\zero\in\KK(\alpha)$ for any
$\alpha\neq0$.
%Each $\alpha$-stable element is infinitely divisible, \ie for every
%$n\geq2$, $z$ can be obtained as a sum of $n$ identical elements from
%$\KK$. 
In particular, $\KK(\infty)$ is the set of \emph{idempotent} elements
that satisfy $z+z=z$, and $\KK(1)$ consists of all $z\in\KK$ that
satisfy (\ref{eq:cone4}).

\begin{lemma}
  If the second distributivity law (\ref{eq:cone4}) holds, then
  $\KK(\alpha)=\{\neutral\}$ for any $\alpha\neq1$.
\end{lemma}
\begin{proof}
  Putting $a=b$ in~\eqref{eq:selfsas} and using (\ref{eq:cone4})
  yields $2^{1-1/\alpha}z=z$ and $2^{1/\alpha-1}z=z$.  By iteration,
  we obtain that $\beta^nz=z$, where $\beta=2^{1-1/\alpha}$ for
  $\alpha>1$ and $\beta=2^{1/\alpha-1}$, for $\alpha<1$. Passing to
  the limit, we have that $z=\zero$ and thus $z=\neutral$ by
  Lemma~\ref{lemma:n=0}\textsl{(i)}.
\end{proof}

\subsection{Norm and metric}
\label{sec:normed-cones}

\begin{definition}
  \label{def:norm-cone}
  A pointed cone $\KK$ is said to be a \emph{normed cone} if
  $\KK$ (or $\KK\setminus\{\neutral\}$ if $\zero\neq\neutral$) is
  metrisable by a metric $d$ which is \emph{homogeneous} at the
  origin, \ie $d(ax,\zero)=ad(x,\zero)$ for every $a>0$ and $x\in\KK$.
  The value $\|x\|=d(x,\zero)$ is called the \emph{norm} of $x$.
\end{definition}

In Sections~\ref{sec:lepage-series-cone} and~\ref{sec:conv-stable-laws-1}
it is assumed that $\KK$ is a normed cone. Note that the function
$d(x,\zero)$ should be called a \emph{gauge function} rather than a
norm, since it is not assumed to be sub-linear, \ie $d(x+y,\zero)$ is
not necessarily smaller than $d(x,\zero)+d(y,\zero)$. However, we
decided to use the word \emph{norm} in this context, because we employ
the gauge function to define the balls and spheres in exactly the same
way as the conventional norm is used. Most of further results can be
reformulated for metrics that are homogeneous of a given order $r>0$,
\ie $d(ax,\zero)=a^rd(x,\zero)$.

It is obvious that $\|x\|=0$ if and only if $x=\zero$. Furthermore,
$\|ax\|=a\|x\|$ for all $a>0$ and $x\in\KK$. If $\neutral\neq\zero$,
then \eqref{cone4} implies that
$\|\neutral\|=d(\neutral,\zero)=\infty$. It is therefore essential to
allow for $d$ to take infinite values, \cf Definition~\ref{def:cone}.
For instance, if $\KK$ is the cone $\R_+=[0,\infty)$ with the minimum
operation, then the Euclidean distance from any nonempty $x\in\R_+$ to
$\infty$ (being the neutral element) is infinite.

If $\KK$ is a linear space, the metric and the norm can be routinely
constructed using a star-shaped neighbourhood of the origin. If this
neighbourhood is convex, then $\KK$ is a locally convex topological
vector space and the corresponding norm is \emph{sub-linear}, \ie
\begin{equation}
  \label{eq:norm-sublinear}
  \|x+y\|\leq\|x\|+\|y\|\,.
\end{equation}

The open ball of radius $r$ centred at $\zero$ is denoted
by
\begin{displaymath}
  B_r=\{x\in\KK:\; \|x\|< r\}\,.
\end{displaymath}
The interior of its complement is given by
\begin{displaymath}
  \cset=\{x:\; \|x\|> r\}\,.
\end{displaymath}
If $\neutral\neq\zero$, then $\neutral\in \cset$
for all $r>0$.  The set
\begin{displaymath}
  \SS=\{x:\; \|x\|=1\}
\end{displaymath}
is called the \emph{unit sphere}.  Note that $\SS$ is complete with
respect to the metric induced by the metric on $\KK$. The existence of
the origin implies that $\|x\|<\infty$ for all
$x\in\KK\setminus\{\neutral\}$, therefore $\KK$ admits a \emph{polar
  decomposition}. This decomposition is realised by the bijection
$x\leftrightarrow(\|x\|,x/\|x\|)$ between
\begin{displaymath}
  \EE=\KK\setminus\{\zero,\neutral\}
\end{displaymath}
and $(0,\infty)\times\SS$.

\subsection{Sub-invariance}
\label{sec:sub-invar-prop}

In addition to the homogeneity property of the metric $d$, we
sometimes require that
\begin{equation}
  \label{triangleq}
  d(x+h,x) \leq d(h,\zero)=\|h\|,\quad x,h\in\KK\,.
\end{equation}
Then the metric (or the norm) in $\KK$ is said to be
\emph{sub-invariant}. This technical condition guarantees the uniform
continuity of the norm, \ie closeness of $\|x\|$ and $\|x+h\|$ if
$\|h\|$ is small.  Indeed, by the triangular inequality,
$d(x+h,\zero)\leq d(x+h,x)+d(x,\zero)$, implying the sub-linearity of
the norm (\ref{eq:norm-sublinear}) and also that $\|x+h\|-\|x\|\leq
d(x+h,x)$. On the other hand, $d(x,\zero)\leq d(x,x+h)+d(x+h,\zero)$,
so that $-d(x,x+h)\leq \|x+h\|-\|x\|$.  Therefore, in view
of~(\ref{triangleq}),
\begin{displaymath}
  \big|\,\|x+h\|-\|x\|\,\big| \leq d(x+h,x)\leq \|h\|\,.
\end{displaymath}
In particular,
\begin{equation}
  \label{eq:triang3}
  \|x\|\leq \|x+h\|+\|h\|,\quad x,h\in\KK\,.
\end{equation}
The two `$+$' in the right-hand side of (\ref{eq:triang3}) should not
be confused: the first one is the addition operation in $\KK$ while
the second one is the conventional sum of positive numbers.  If $\KK$
is a group, then an \emph{invariant} (thus also sub-invariant) metric
always exists, see \cite[Ch.~6, p.~210]{kel55}, \ie (\ref{triangleq})
holds with the equality sign.  In general, (\ref{triangleq}) is too
restrictive, \eg for $\R_+$ with the maximum operation and Euclidean
metric.

\begin{lemma}
  \label{lemma:sub-invar}
  If $\KK$ has a sub-invariant norm, then $\zero=\neutral$ and, for
  any $\alpha\in(0,1)$, $\neutral$ is the only element with a finite
  norm that belongs to $\KK(\alpha)$.
\end{lemma}
\begin{proof}
  Applying (\ref{triangleq}) with $th$ instead of $h$ and letting
  $t\downarrow0$ implies that $x+\zero=x$ for all $x$, whence
  $\neutral=\zero$. The sub-invariant norm also satisfies
  (\ref{eq:norm-sublinear}), so that, for each $z\in\KK(\alpha)$
  \begin{displaymath}
    (a+b)^{1/\alpha}\|z\|=\|(a+b)^{1/\alpha}z\|
    =\|a^{1/\alpha}z+b^{1/\alpha}z\|
    \leq (a^{1/\alpha}+b^{1/\alpha})\|z\|\,.
  \end{displaymath}
  If $z\neq\zero$ and $\|z\|\neq \infty$ then $(a+b)^{1/\alpha}\leq
  a^{1/\alpha}+b^{1/\alpha}$ for all $a,b>0$ is impossible for
  $\alpha\in(0,1)$.
\end{proof}

\bigskip

Typical examples of cones that fulfil our requirements are Banach
spaces or convex cones in Banach spaces; the family of compact (or
convex compact) subsets of a Banach space with Minkowski addition; the
family of compact sets in $\R^d$ with the union operation; the family
of all finite measures with the conventional addition operation and
multiplication by numbers.  Another typical example is the set
$\R_+=[0,\infty)$ with the maximum operation $x+y=x\vee y=\max(x,y)$.
In order to distinguish this example from the conventional cone
$(\R_+,+)$, we denote it by $(\R_+, \vee)$. These and other examples
are discussed in Section~\ref{sec:examples}.

\bigskip

We say that $\KK$ can be isometrically embedded in a Banach space
$\BB$ if there exists an injection $I:\KK\to\BB$ such that $I(ax+by) = aI(x)
+ bI(y)$ for all $a,b>0$ and $x,y \in \KK$ and $d(x,y)= \|I(x) -
I(y)\|$ for all $x,y\in\KK$.  However, this embedding is possible only
under some conditions on the cone $\KK$ and the corresponding metric.

\begin{Th}
  \label{prop:4}
  A convex cone $\KK$ with a metric $d$ can be embedded into a Banach
  space if and only if the second distributivity law~(\ref{eq:cone4})
  holds and $d$ is homogeneous and invariant, \ie $d(ax,ay)=ad(x,y)$
  and $d(x+z,\,y+z) = d(x,y)$ for all $a>0$ and all $x,y,z\in\KK$.
\end{Th}
\begin{proof}
  Necessity is obvious. Sufficiency follows from H\"ormander's theorem
  \cite{hor54}, see also~\cite{keim67}.
\end{proof}

\section{LePage series on a cone}
\label{sec:lepage-series-cone}

\subsection{Point processes on a cone}
\label{sec:point-processes-cone}

Consider a normed cone $\KK$ with its Borel $\sigma$-algebra
$\sB(\KK)$.  Let $\sMzero$ (respectively $\sM$) be the family of
\emph{counting measures} $m$ on $\sB(\KK)$ such that
$m(\cset)<\infty$ (respectively $m(B_r)<\infty$) for every $r>0$.
Both families $\sMzero$ and $\sM$ also contain the null-measure.
Denote by $\delta_x$ the unit mass measure concentrated at $x\in\KK$.
Any counting measure can be represented as
\begin{equation}
  \label{eq:delta}
  m=\delta_{x_1}+\delta_{x_2}+\cdots=\ssum_{i}\delta_{x_i}\,,
\end{equation}
where $x_1,x_2,\ldots$ is an at most countable collection of points
such that only a finite number of $x_i$'s lies in $\cset$ if
$m\in\sMzero$ or in $B_r$ if $m\in\sM$ for every $r>0$.  When
considering $\sMzero$ we deviate from the typical setting in the
theory of point processes (\seg \cite{da:vj}), where the point sets
are assumed to be locally finite. Our setting allows for a
concentration point at the origin for measures from $\sMzero$.  To
cover these both cases with the same notation, let $\setr$ denote
$\cset$ ($B_{r^{-1}}$, respectively) in case we consider measures from
$\sMzero$ ($\sM$, respectively). Then we always have $m(\setr)<\infty$
whenever $m\in\sM$ or $m\in\sMzero$.

The counting measure $m$ is said to be \emph{simple} if the points
$x_1,x_2,\ldots$ from (\ref{eq:delta}) are distinct. A simple counting
measure is fully characterised by its support
\begin{displaymath}
  \supp m=\{x\in\KK :\ m(\{x\})>0\}\,.
\end{displaymath}

A \emph{point process} $\mu$ is a measurable map from some probability
space into $\sMzero$ (or $\sM$) with the $\sigma$-algebra generated by
the sets of measures $m\in\sMzero$ (or $m\in\sM$) such that $m(B)=n$
for Borel sets $B\subset\KK$ and $n\geq0$.  The distribution of $\mu$
is denoted by $\P$.  The process is \emph{simple} if almost all its
realisations are simple.  The \emph{probability generating functional}
of $\mu$ is defined as
\begin{displaymath}
  G_\mu(u)=\E\exp\Bigl\{\int_\KK \log u(x)\,
  \mu(dx)\Bigr\}=\E\Bigl[\prod_{x_i\in\supp\mu}
    u(x_i)^{\mu(\{x_i\})}\Bigr]\,,  
\end{displaymath}
where $\E$ is the expectation with respect to $\P$ and
$u:\KK\mapsto(0,1]$ is a function that identically equals 1 on the
complement of $\setr$ for some $r>0$, \cf \cite[Sec.~7.4]{da:vj}.

If $F$ is a Borel set, then $\int_F x\mu(dx)$ is the random element
obtained as the sum of all the points from $F\cap \supp\mu$ taking
into account possible multiplicities. If $\mu(F)$ is \as finite, this
integral is a well-defined finite sum. Otherwise, \eg if one considers
\begin{displaymath}
  \int x\mu(dx)=\int_{\KK} x\mu(dx)
\end{displaymath}
with the integration over $F=\KK$, the almost sure convergence of this
integral is understood as the existence of the integral for almost all
$\mu$. The \emph{absolute convergence} requires the existence of
$\int\|x\|\mu(dx)$. A weaker condition is the convergence of the
\emph{principal value} which requires that $\int_{\setr} x\mu(dx)$
converges as $r\downarrow0$. Similar definitions are applicable for
the integral $\int g(x)\mu(dx)$, where $g:\KK\mapsto\KK$ is a
measurable function.

\subsection{Stable Poisson process}
\label{sec:stable-poiss-proc}

Let $\Lambda$ be a measure on $\KK$ which is finite on all $\setr$,
$r>0$. A point process $\Pi$ is called a \emph{Poisson process} with
\emph{intensity measure} $\Lambda$ if, for any disjoint family of Borel
sets $F_1,\dots,F_n$, the random variables $\Pi(F_1),\dots,\Pi(F_n)$
are jointly independent Poisson distributed with means
$\Lambda(F_1),\dots,\Lambda(F_n)$, respectively.  The Poisson process
is simple if and only if its intensity measure is non-atomic.

Given the automorphisms $D_a: x\mapsto ax$ introduced in
Section~\ref{sec:basic-definitions}, $D_am$ denotes the image of $m$,
\ie $(D_am)(A)=m(D_a^{-1}A)=m(D_{a^{-1}}A)$ for every Borel $A$.
If $m=\sum_i \delta_{x_i}$ is a counting measure, then
$D_am=\sum_i \delta_{ax_i}$, in particular,
$\supp(D_am)=a(\supp m)$ and
\begin{equation}
  \label{eq:intmu}
  \int x D_am(dx)=a\int xm(dx)\,.
\end{equation}
An important property of a Poisson process $\Pi$ is that $D_a\Pi$ is
again a Poisson process driven by the intensity measure $D_a\Lambda$,
if $\Lambda$ is the intensity measure of $\Pi$.

\begin{example}[Stable Poisson point process]
  \label{example2}
  Recall that $\EE=\KK\setminus\{\zero,\neutral\}$ can be identified
  with $(0,\infty)\times\SS$ using the polar decomposition. Define a
  measure $\Lambda$ (also denoted by $\Lambda_{\alpha,\sigma}$) on
  $\EE$ as the product of the measure $\theta_\alpha$ on $(0,\infty)$
  such that $\alpha\neq 0$ and
  \begin{equation}
    \label{eq:defm}
    \begin{cases}
      \theta_\alpha((r,\infty))=r^{-\alpha} & \text{if}\; \alpha>0 \,,
      \\
      \theta_\alpha((0,r))=r^{-\alpha} & \text{if}\; \alpha<0 \,, 
    \end{cases}
    \quad r>0\,,
  \end{equation} 
  and a finite measure $\sigma$ on the Borel $\sigma$-algebra
  $\sB(\SS)$ induced on $\SS$. The Poisson process on $\EE$ (and thus
  also on $\KK$) with intensity measure $\Lambda_{\alpha,\sigma}$ is
  denoted by $\Pi_{\alpha,\sigma}$. We omit the subscript $\sigma$ and
  write simply $\Pi_\alpha$ when no confusion occurs. The measure
  $\sigma$ is called its \emph{spectral measure}.  If $\alpha>0$, then
  $\Pi_\alpha\in\sMzero$, \ie it has the concentration point at the
  origin. If $\alpha<0$, then $\Pi_\alpha\in\sM$, and if
  $\neutral\neq\zero$, then the support points of $\Pi_\alpha$ have
  the concentration point at $\neutral$.
\end{example}

The importance of the process $\Pi_\alpha$ in our context stems from
the fact that it is \emph{stable} with respect to the addition
operation applied to the corresponding counting measures.
  
\begin{Th}
  \label{prop:5}
  Let $\Pi'_\alpha$ and $\Pi''_\alpha$ be two independent
  copies of $\Pi_\alpha$. Then
  \begin{equation}
    \label{eq:7}
    D_{a^{1/\alpha}}\Pi'_\alpha+D_{b^{1/\alpha}}\Pi''_\alpha \deq
    D_{(a+b)^{1/\alpha}}\,\Pi_\alpha
  \end{equation}
  for all $a,b>0$.
\end{Th}
\begin{proof}
  By~(\ref{eq:defm}), the intensity measure
  $\Lambda=\Lambda_{\alpha,\sigma}$ of $\Pi_\alpha$ satisfies
  \begin{equation}
    \label{eq:lams}
    (D_{a^{1/\alpha}}\Lambda)(A)=\Lambda(D_{a^{-1/\alpha}}A)
    = a\Lambda(A)
  \end{equation}
  for any Borel $A$.  Note that the left-hand side is the result of a
  contraction of the phase space, while the right-hand side is a
  measure obtained by multiplying the values of $\Lambda$ by a number.
  Thus the processes on the both sides of (\ref{eq:7}) are Poisson
  with the same intensity measure $(a+b)\Lambda$.
\end{proof}

Property (\ref{eq:7}) can be reformulated as the stability
property of $\Pi_\alpha$ with respect to the union operation applied
to its support set.  Let $\kappa_\alpha=\supp \Pi_\alpha$ and let
$\kappa'_\alpha$ and $\kappa''_\alpha$ be two independent copies of
$\kappa_\alpha$. Since the intensity measure $\Lambda$ of $\Pi_\alpha$ is
non-atomic, $a^{1/\alpha}\kappa'_\alpha\cap
b^{1/\alpha}\kappa''_\alpha=\emptyset$ with probability~1. Then
(\ref{eq:7}) is equivalent to
\begin{equation}
  \label{eq:7bis}
  a^{1/\alpha}\kappa'_\alpha \cup
  b^{1/\alpha}\kappa''_\alpha \deq (a+b)^{1/\alpha}\kappa_\alpha
\end{equation}
for all $a,b>0$. This means that $\kappa_\alpha$ is a union-stable
random closed set, \ie a $\sas$ random element in the cone of closed
sets with the union operation, see \cite{mo93l}.

The following theorem provides a useful representation of the process
$\Pi_\alpha$, especially for simulation purposes. Its proof relies
on basic facts on transformations of a Poisson process.
% and can be found in \cite[p.~26]{sam:taq94}.

\begin{Th}
  \label{prop:2}
  Let $\{\zeta_k\}$ and $\{\eps_k\}$ be two independent sequences of
  \iid random variables, where $\zeta_k$, $k\geq1$, have exponential
  distribution with mean 1, and $\eps_k$, $k\geq1$ are distributed on
  $\SS$ according to $\hat{\sigma}(\cdot)=\sigma(\cdot)/\sigma(\SS)$,
  where $\sigma$ is a finite measure on $\SS$. Define
  $c=\sigma(\SS)^{1/\alpha}$ and $\Gamma_k=\zeta_1+\dots+\zeta_k$,
  $k\geq1$.  Then, for any $\alpha\neq0$,
  \begin{equation}
    \label{eq:prepr}
    \Pi_{\alpha,\sigma}\deq \sum_{k=1}^{\infty} \delta_{\Gamma^{-1/\alpha}_k
      \eps_k c}\,.
  \end{equation}
\end{Th}

By conditioning with respect to the number of points of
$\Pi_\alpha=\Pi_{\alpha,\sigma}$ in $\setr$, it is easy to calculate
its probability generating functional
\begin{equation}
  \label{eq:pgfp}
  G_{\Pi_\alpha}(u)=\exp\left\{-\int_\KK(1-u(x))\Lambda(dx)\right\}\,,
\end{equation}
where $\Lambda=\Lambda_{\alpha,\sigma}$ and $u:\KK\mapsto[0,1]$ is a
function that is identically equal to $1$ outside of $\setr$, see
\cite[Ex.~7.4(a)]{da:vj}. By passing to a limit, it is possible to
show that (\ref{eq:pgfp}) holds for any function $u$ such that $1-u$
is integrable with respect to $\Lambda$.

\begin{lemma}
  \label{lem:multiple}
  Let $\Pi_\alpha^{(\eta)}$ be the point process obtained by
  multiplying all points from $\Pi_{\alpha,\sigma}$ by \iid
  realisations of a positive random variable $\eta$ with
  $a=\E\eta^\alpha<\infty$. Then $\Pi_\alpha^{(\eta)}$ has spectral
  measure $a\sigma$ and coincides in distribution with the point
  process $D_{a^{1/\alpha}}\Pi_\alpha$.
\end{lemma}
\begin{proof}
  It suffices to show that the probability generating functionals of
  the both processes coincide. By conditioning with respect to the
  realisation of $\Pi_\alpha$, one easily obtains that
  \begin{displaymath}
    G_{\Pi_\alpha^{(\eta)}}(u)=\exp\left\{-\int_\KK(1-\E u(x\eta))
      \Lambda(dx)\right\}.
  \end{displaymath}
  By changing variables $x\eta=y$ and using (\ref{eq:lams}), it is
  easily seen that the exponent in the formula for the probability
  generating functional is
  \begin{displaymath}
    -\E\eta^\alpha \int_\KK (1-u(x))\Lambda(dx)
  \end{displaymath}
  that corresponds to the probability generating functional of
  $D_{a^{1/\alpha}}\Pi_\alpha$.
\end{proof}

\begin{remark}
  \label{rem:pi-alpha}
  A stable Poisson process $\Pi_\alpha$ can be defined on any convex
  cone $\KK$ without assuming that $\KK$ possesses the origin or a
  norm. It suffices to consider a Poisson point process whose
  intensity measure $\Lambda$ is homogeneous, \ie satisfies
  \eqref{eq:lams}. 
\end{remark}

\subsection{LePage series}
\label{sec:lepage-series}

A $\KK$-valued random element $\xi$ is said to have $\sas$
distribution if it satisfies (\ref{eq:a1alph-a+b1-}) with the addition
and multiplication operations defined on $\KK$.  Theorem~\ref{th:1}
below provides a rich family of $\sas$ distributions by their series
decomposition.

\begin{Th}
  \label{th:1}
  Let $\{\zeta_k,\, k\geq1\}$ be \iid exponentially distributed random
  variables with mean 1 and let $\Gamma_k=\zeta_1+\dots+\zeta_k$,
  $k=1,2,\dots$. Furthermore, let $\{\eps_k\}$ be independent of
  $\{\zeta_k\}$ \iid random elements on the unit sphere $\SS$ in $\KK$
  with a common distribution $\hat\sigma$. If the principal value of
  the integral $\int x\,\Pi_{\alpha,\hat\sigma}(dx)$ is finite with
  probability 1, then for any $z\in\KK(\alpha)$ and $c\geq0$, the
  series
  \begin{equation}
    \label{eq:9}
    \xi_\alpha=z+c\sum_{k=1}^\infty\Gamma_k^{-1/\alpha}\eps_k
  \end{equation}
  converges almost surely and $\xi_\alpha$ has a $\sas$ distribution
  on $\KK$.
  
  If the norm on $\KK$ is sub-invariant, then the infinite sum
  in~(\ref{eq:9}) converges absolutely \as for any $\alpha\in(0,1)$.
\end{Th}
\begin{proof}
  The case $c=0$ is trivial.  Without loss of generality assume
  that $z=\neutral$ in (\ref{eq:9}).  Note that
  \begin{displaymath}
    c\sum_{k=1}^\infty\Gamma_k^{-1/\alpha}\eps_k=\int x\,
    \Pi_\alpha(dx)\,,
  \end{displaymath}
  where $\Pi_\alpha=\Pi_{\alpha,\sigma}$ with
  $\sigma=c^\alpha\hat\sigma$.  Recall that $\setr$ denotes $\cset$ if
  $\alpha>0$ and $B_{r^{-1}}$ if $\alpha<0$.  By~(\ref{eq:7}) and
  almost sure finiteness of $\Pi_\alpha(\setr)$ we have that
  \begin{displaymath}
    \int_{\setr} x\,D_{a^{1/\alpha}}\Pi'_\alpha(dx)+
      \int_{\setr} x\,D_{b^{1/\alpha}}\Pi''_\alpha(dx) \deq
        \int_{\setr} x\,D_{(a+b)^{1/\alpha}}\Pi_\alpha(dx)\,.
  \end{displaymath}
  Since the principal value of $\int x\Pi_{\alpha,\hat\sigma}(dx)$ is
  finite by the condition, \eqref{eq:intmu} implies that the integrals
  above also converge. Then we can let $r\downarrow 0$ to obtain that
  \begin{equation}\label{eq:stable}
    a^{1/\alpha}\xi'_\alpha + b^{1/\alpha}\xi''_\alpha\deq
    (a+b)^{1/\alpha}\xi_\alpha\,,
  \end{equation}
  where $\xi'_\alpha$ and $\xi''_\alpha$ are two independent copies of
  $\xi_\alpha$, \ie $\xi_\alpha$ is $\sas$.
  
  Now assume that the norm is sub-invariant, \ie it satisfies
  (\ref{triangleq}).  If $S_n=\sum_{k=1}^n\Gamma_k^{-1/\alpha}\eps_k$,
  then
  \begin{displaymath}
    d(S_n,S_{n+m})\leq
    \Bigl\|\sum_{k=n+1}^{n+m}\Gamma_k^{-1/\alpha}\eps_k\Bigr\| \leq
    \sum_{k=n+1}^{n+m}\Gamma_k^{-1/\alpha}\,.
  \end{displaymath}
  Note that the first sum corresponds to the semigroup addition in
  $\KK$ while the second is the ordinary sum of positive numbers.  The
  last expression vanishes almost surely as $n\ti$, since
  $\Gamma_k^{-1/\alpha} \sim k^{-1/\alpha}$ with probability 1 by the
  strong law of large numbers. Referring to the completeness of $\KK$
  confirms the a.\,s.  convergence of the series.
\end{proof}

The distribution of $\xi_\alpha$ given by (\ref{eq:9}) is determined
by its deterministic part $z$ and the measure
$\sigma=c^\alpha\hat\sigma$ on $\SS$ called the \emph{spectral
  measure}.  The series (\ref{eq:9}), also written as
\begin{equation}
  \label{eq:lepage-int}
  \xi_\alpha=z+\int x\Pi_\alpha(dx)\,,
\end{equation}
is called the \emph{LePage series} on $\KK$. The convergence of the
principal value of the integral in (\ref{eq:lepage-int}) is equivalent
to the convergence of the series in (\ref{eq:9}). Theorem~\ref{th:1}
implies that for $\alpha\in(0,1)$ under the sub-invariance condition,
the integral in~(\ref{eq:lepage-int}) converges absolutely a.\ssp s.
The convergence property for other $\alpha$ depends on finer
properties of the cone.  For instance, in $(\R_+,+)$ the series
(\ref{eq:9}) converges if and only if $\alpha\in(0,1)$. In a Banach
space the properly centred (or symmetrised) series also converges for
$\alpha\in[1,2)$, see~\cite[Sec.~1.5]{sam:taq94}.  In some cones the
series (\ref{eq:9}) converges for all positive or negative $\alpha$ as
the following example shows.

\begin{example}
  \label{ex:rvee}
  Consider the cone $(\R_+,\vee)$. If $\eps_k=1$ for all $k\geq1$ and
  $\alpha>0$, then $\xi_\alpha=c\Gamma_1^{-1/\alpha}$ has the Fr\'echet
  distribution which is max-stable. The Weibull distribution arises in
  $([0,\infty],\min)$ cone for $\alpha<0$. Here any
  $\sas$-distribution for $\alpha>0$ is trivial concentrated at~$0$.
\end{example}

\begin{remark}
  The point $z$ in~\eqref{eq:lepage-int} does not have to be
  $\alpha$-stable itself in order to produce a $\sas$ element
  $\xi_\alpha$.  For instance, $\xi_\alpha$ has the Fr\'echet
  distribution which is stable in the cone $(\R,\vee)$ for $z$ being
  any negative number and $\eps_k=1$ for all $k\geq 1$, \cf
  Corollary~\ref{cor:decomp}.
\end{remark}

\begin{remark}
  A convergent LePage series yields a (possibly degenerate)
  $\sas$ random element in any convex cone $\KK$. It is not essential
  that the topology on $\KK$ is metrisable or $\KK$ is equipped with a norm.
  The only assumption is that $\KK$ possesses a polar decomposition.
\end{remark}

A simple example of a $\sas$ law in any Banach space is provided by
the random element $\xi=\eta x_0$, where $\eta$ is non-negative
$\sas$ random variable and $x_0$ is any non-zero vector.  The same
construction is applicable to any convex cone that satisfies the
second distributivity law (\ref{eq:cone4}). However, if
(\ref{eq:cone4}) does not hold, this construction no longer works.
For instance, if $\KK=(\R_+,\vee)$, then $\xi=\eta x_0$ with a
non-negative $\sas$ random variable $\eta$ is not $\sas$, since for
$\xi_1$ and $\xi_2$ being two independent copies of $\xi$, the random
elements $\xi_1+\xi_2=(\eta_1\vee\eta_2)x_0$ and $2^{1/\alpha}\xi$
have different distributions.  In contrast, Theorem~\ref{th:1} yields
that
\begin{displaymath}
  Y_\alpha=\max_{k\geq1}\Gamma_k^{-1/\alpha}x_0
  =\Gamma_1^{-1/\alpha}x_0  
\end{displaymath}
is $\sas$ and so provides an example of a $\sas$ distribution in this
cone.

\subsection{L\'evy process}
\label{sec:gener-lepage-seri}

In this section it is assumed that $\KK$ has a sub-invariant norm.
Let $\{\Gamma_k\}$ and $\{\eps_k\}$ be as in Theorem \ref{th:1}. If
$\alpha\in(0,1)$ and $\{\eta_k\}$ are \iid copies of a non-negative
random variable $\eta$ with $c=\E\eta^\alpha<\infty$, then
Theorem~\ref{th:1} and Lemma~\ref{lem:multiple} imply that
\begin{equation}
  \label{eq:9a}
  Y_\alpha^{(\eta)}=\sum_{k=1}^\infty\Gamma_k^{-1/\alpha}\eta_k\eps_k
\end{equation}
absolutely converges almost surely and has $\sas$ distribution on
$\KK$ with spectral measure $\sigma=c\hat\sigma$.

Representation~(\ref{eq:9a}) can be used to construct a $\KK$-valued
$\sas$ L\'evy process, \ie $\KK$-valued process with independent and
stationary $\sas$ increments. A $\KK$-valued stochastic process $X_t$,
$t\geq0$, is said to have independent increments, if, for every
$n\geq1$ and $0\leq t_0<t_1<\cdots<t_n$, there exist jointly
independent $\KK$-valued random elements
$\xi_{t_0t_1},\dots,\xi_{t_{n-1}t_n}$ (representing the increments)
such that the joint distributions of $X_{t_0},X_{t_1},\dots,X_{t_n}$
and $X_{t_0},X_{t_0}+\xi_{t_0t_1},\dots,X_{t_0}+\xi_{t_0t_1}
+\cdots+\xi_{t_{n-1}t_n}$ coincide.

\begin{Th}
  \label{th:levy}
  Let $\{\eps_k\}$ and $\{\Gamma_k\}$ be as in Theorem~\ref{th:1} and
  $\alpha\in(0,1)$.  Let $\{\tau_k, k\geq1\}$ be a sequence of \iid
  random variables uniformly distributed on $[0,1]$ and independent of
  the sequences $\{\eps_k\}$ and $\{\Gamma_k\}$.  Then the process
  \begin{equation}
    X_\alpha (t)
    = \sum_{k=1}^\infty\Gamma_k^{-1/\alpha}
    \one_{[0,t]}(\tau_k)\eps_k\,,\;\;\; t\in [0,1],
  \end{equation}
  has independent increments given by
  \begin{equation}
    \xi_{ts} = \sum_{k=1}^\infty\Gamma_k^{-1/\alpha}
    \one_{(t,s]}(\tau_k)\eps_k\,, \;\;\; t<s,
  \end{equation}
  which are $\sas$ distributed with the spectral measure
  $(s-t)\hat\sigma$.
\end{Th}
\begin{proof}
  The proof follows immediately from the construction (\ref{eq:9a})
  applied to $\eta_k=\one_{(t,s]}(\tau_k)$, $k\geq1$.
\end{proof}

\begin{example}
  Consider $(\R_+,\vee)$. Since $\SS=\{1\}$,
  \begin{displaymath}
    X_\alpha(t) = \Gamma^{-1/\alpha}_{k(t)}\,, \quad t\in[0,1]\,,
  \end{displaymath}
  where $k(t)=\min\{k:\; \tau_k\leq t\}$. Then $X_\alpha$ is a Markov
  process with \as piecewise constant non-decreasing trajectories on
  $(\eps,1]$ for every $\eps>0$. This process has been studied
  in~\cite{res:roy94s} under the name of a \emph{super-extremal
    process} for $\KK$ being the family of upper semicontinuous
  functions with pointwise maximum operation.
\end{example}

\section{Convergence to stable laws}
\label{sec:conv-stable-laws-1}

\subsection{Weak convergence of point processes}
\label{sec:point-processes}

A stable Poisson process $\Pi_\alpha$ has always infinite number of
support points since its intensity measure is infinite. If $\alpha<0$
this process has finite number of points in any bounded set (and so
has realisations from $\sM$), while for $\alpha>0$ its concentration
point is at the origin $\zero$ and only a finite number of points lies
outside any ball centred at $\zero$. The first case corresponds to
well-known $\sigma$-finite point processes.  In order to study
convergence of processes with possible concentration point at $\zero$
(\ie with realisations from $\sMzero$), we need to amend some
conventional definitions.  

Consider the family $\CC$ of continuous bounded functions $f:\KK\to\R$
such that $f(x)=0$ for all $x\notin \setr$ with some $r>0$.  A
sequence of counting measures $\{m_n,n\geq1\}$ from $\sMzero$ (or from
$\sM$) is said to converge \emph{vaguely} to $m$ (notation $m_n\vto
m$) if
\begin{equation}
  \label{eq:vaguedef}
  \int f(x)\,m_n(dx) \to \int f(x)\,m(dx)\quad \text{as}\ n\ti
\end{equation}
for every $f\in\CC$. For counting measures from $\sM$ this definition
turns into the conventional definition of vague convergence, see
\cite{m:k:m}.  Using the inversion transformation $x\mapsto x/\|x\|^2$
it is easy to confirm that conventional properties of the vague
convergence hold for point processes from $\sMzero$. In particular,
the family $\sMzero$ with the vague topology is a Polish space, see
\cite[Prop.~1.15.5]{m:k:m}.

Note that the vague convergence is often formulated for counting
measures in locally compact spaces, where the functions from $\CC$ have
\emph{compact} supports, \seg~\cite{res87}.  Counting measures
and point processes in general Polish spaces have been systematically
studied in \cite{m:k:m} and \cite{da:vj}.

The following result extends (\ref{eq:vaguedef}) for $\KK$-valued
functions.

\begin{lemma}
  \label{lem:kappa1}
  If $m_n\vto m$ with a finite measure $m$, then 
  \begin{equation}
    \label{eq:new-kappa}
    \int g(x)\,m_n(dx) \to \int g(x)\,m(dx)
  \end{equation}
  for any continuous function $g:\KK\mapsto\KK$ such that
  $g(x)=\neutral$ for all $x\in\setr$ for some $r>0$.
\end{lemma}
\begin{proof}
  Note that $m$ is a finite measure on $\setr$.  By
  \cite[Sec.~1.15]{m:k:m} %RussianProp.~1.9.10]
  it is possible to order the support points of $m_n$ so that they
  converge to the corresponding support points of $m$. This implies
  (\ref{eq:new-kappa}) taking into account the continuity of $g$.
\end{proof}

\emph{Weak convergence} $\mu_n\wto\mu$ of point processes is defined
using the weak convergence of the corresponding probability
distributions, \ie $\E h(\mu_n)\to \E h(\mu)$ as $n\ti$ for any
bounded and continuous in the vague topology function $h$ that maps
$\sM$ or $\sMzero$ into $\R$.

Let $\mu\restr$ denote the restriction of $\mu$ onto $\setr$.
It is shown in \cite[Th.~3.1.13]{m:k:m} %RusTh4.2.9 
that $\mu_n\wto\mu$ for a simple point process $\mu$ if and only if
$\mu_n(B)$ weakly converges to $\mu(B)$ for all $B$ from a certain
subring of $\{B\in\sB(\KK):\; \mu(\partial B)=0\; \text{a.s.}\}$ such
that for each closed $F$ and each open neighbourhood $U\supset F$, one
has $F\subset B\subset U$ for some $B$ from this subring.  Since any
such set $B$ may be chosen to be separated from the origin by a
positive distance $r$ (in case of $\sMzero$) or contained in a ball of
radius $r^{-1}$ (in case of $\sM$), and $\mu(B)=\mu\restr(B)$, the
following result holds.

\begin{lemma}
  \label{prop:nd1}
  A sequence $\{\mu_n,n\geq1\}$ of point processes weakly converges
  to a simple point process $\mu$ if and only if
  $\mu_n\restr\wto\mu\restr$ as $n\ti$ for all $r>0$ such that
  $\mu(\partial A_r)=0$ almost surely.
\end{lemma}

\subsection{Convergence of binomial processes}
\label{sec:conv-binom-proc}

Let $\{\xi_k,\, k\geq1\}$ be a sequence of \iid $\KK$-valued random
elements. For every $n\ge1$, $\sum_{i=1}^n\delta_{\xi_i}$ is called
the \textit{binomial point process}.  It is simple if and only if the
distribution of $\xi_k$ is non-atomic.  The scaled versions of the
binomial process are defined by
\begin{equation}
  \label{eq:beta_n=s-delt-}
  \beta_n=\sum_{k=1}^n \delta_{\xi_k/b_n}\,, \quad n\geq1\,,
\end{equation}
for a sequence $\{b_n,\ n\geq1\}$ of normalising constants. We shall
typically have 
\begin{equation}
  \label{eq:b_n=n1-quad-ngeq1}
  b_n=n^{1/\alpha}L(n)\,,\quad n\geq1\,,
\end{equation}
with $\alpha\neq0$ and a slowly varying at infinity function $L$.

The following theorem shows that the Poisson point process
$\Pi_\alpha$ with $\alpha>0$ arises as a weak limit for binomial
processes $\beta_n$ defined by (\ref{eq:beta_n=s-delt-}) if the
$\xi_k$'s have regularly varying tails.

\begin{Th}
  \label{prop:1}
  Let $\xi,\xi_1,\xi_2,\ldots$ be \iid $\EE$-valued random elements.
  Then $\beta_n\wto \Pi_\alpha$ as $n\ti$ for $\alpha>0$ if and only
  if there exists a finite measure $\sigma$ on $\sB(\SS)$ such that
  \begin{equation}
    \lim_{n\ti} n\P\Bigl\{\frac{\xi}{\|\xi\|}\in G,\
    \|\xi\|>rb_n\Bigr\}=\sigma(G)r^{-\alpha}
    \label{eq:resnick}
  \end{equation}
  for all $r>0$ and $G\in\sB(\SS)$ with $\sigma(\partial G)=0$, where
  the $b_n$'s are given by (\ref{eq:b_n=n1-quad-ngeq1}).
\end{Th}

This result is similar to Proposition~3.21 in \cite[p.~154]{res87} for
$\KK$ being a locally compact space.  Note that the condition of type
(\ref{eq:resnick}) is typical in limit theorems for sums of random
elements, see \cite[p.~167]{ara:gin}. The condition of $\xi$ being
$\EE$-valued ensures that $0<\|\xi\|<\infty$ \as

\begin{proof}[Proof of Theorem~\ref{prop:1}]
  \textsl{Sufficiency.} By Lemma~\ref{prop:nd1}, it suffices to
  show that $\beta_n\restr\wto\Pi_\alpha\restr$ for all $r>0$.  Since
  $\SS$ is a Polish space, for every $\delta>0$, there exists a
  compact set $\SS_\delta\subset\SS$ such that
  $\sigma(\SS\setminus\SS_\delta)<\delta$. Denote by
  $\beta_n\restr[r,\delta]$ (respectively
  $\Pi_\alpha\restr[r,\delta]$) the restriction of $\beta_n$
  (respectively $\Pi_\alpha$) onto the set
  $[r,\infty)\times \SS_\delta$.
  
  Consider a continuous in the vague topology function
  $F:\sMzero\mapsto\R$ with the absolute value bounded by $a$. Since
  $[r,\infty)\times \SS_\delta$ is a locally compact space,
  \cite[Prop.~3.21]{res87} yields that $\E
  F(\beta_n\restr[r,\delta])\to\E F(\Pi_\alpha\restr[r,\delta])$ as
  $n\ti$.  Note that $r^{-\alpha}$ in (\ref{eq:resnick}) coincides
  with $\theta_\alpha((r,\infty))$.
    
  Furthermore,
  \begin{displaymath}
    |\E F(\beta_n\restr)-\E F(\beta_n\restr[r,\delta])| 
    \leq 2a\P\{\beta_n\restr[r,\delta]\neq\beta_n\restr\}\,.
  \end{displaymath}
  The latter probability is bounded by $n\P\{\|\xi\|\geq
  rb_n,\xi/\|\xi\|\notin \SS_\delta\}$. By (\ref{eq:resnick}), its
  upper limit as $n\ti$ does not exceed $\delta r^{-\alpha}$, which
  can be made arbitrarily small for every fixed $r>0$ by the choice of
  $\delta$.
  
  Similarly, the absolute difference between $\E F(\Pi_\alpha\restr)$
  and $\E F(\Pi_\alpha\restr[r,\delta])$ is bounded from above by
  \begin{displaymath}
    2a\P\{\Pi_\alpha\restr[r,\delta]\neq\Pi_\alpha\restr\}\leq
    2a\Big(1-\exp\{-\theta_\alpha(r,\infty)\,
    \sigma(\SS\setminus\SS_\delta)\}\Big)
    \leq 2ar^{-\alpha}\delta\,.
  \end{displaymath}
  Therefore $\E F(\beta_n\restr)\to\E F(\Pi_\alpha\restr)$ as $n\ti$.
  
  \textsl{Necessity.} Let $B=(r,\infty)\times G$, where $G$ is a Borel
  subset of $\SS$ with $\sigma(\partial G)=0$. Then
  $\P\{\Pi_\alpha(\partial B)>0\}=0$, whence
  by~\cite[Th.~3.1.13]{m:k:m} %RusTh4.2.9 
  $\beta_n(B)$ converges weakly to $\Pi_\alpha(B)$. The former has the
  binomial distribution with mean
  \begin{displaymath}
    \E \beta_n(B)=n\P\{\|\xi\|>rb_n,\; \xi/\|\xi\|\in G\}\,,
  \end{displaymath}
  while the latter has the Poisson distribution with mean
  $r^{-\alpha}\sigma(G)$. Therefore,
  \begin{displaymath}
    n\P\{\|\xi\|>rb_n,\; \xi/\|\xi\|\in G\}\to r^{-\alpha}\sigma(G)\,.
  \end{displaymath}
\end{proof} 

By applying the inversion transformation $x\mapsto x\|x\|^{-2}$ it is
possible to convert all measures from $\sMzero$ to measures from
$\sM$. Therefore, an analogue of Theorem~\ref{prop:1} holds for
$\alpha<0$.

\begin{Cor}
  \label{prop:1-minus}
  Let $\xi,\xi_1,\xi_2,\ldots$ be \iid $\EE$-valued random elements.
  Then $\beta_n\wto \Pi_\alpha$ as $n\ti$ for $\alpha<0$ if and only
  if there exists a finite measure $\sigma$ on $\sB(\SS)$ such that
  \begin{equation}
    \lim_{n\ti} n\P\Bigl\{\frac{\xi}{\|\xi\|}\in G,\
    \|\xi\|< rb_n\Bigr\}=\sigma(G)r^{-\alpha}
    \label{eq:resnick-minus}
  \end{equation}
  for all $r>0$ and $G\in\sB(\SS)$ with $\sigma(\partial G)=0$, where
  the $b_n$'s are given by (\ref{eq:b_n=n1-quad-ngeq1}).
\end{Cor} 

\subsection{Convergence of sums}
\label{sec:convergence-sums}

Normalised sums of random elements in $\KK$ can be represented as sums
of points of the corresponding binomial processes $\beta_n$. In order
to derive the limit theorem for the normalised sums from the
convergence $\beta_n\wto\Pi_\alpha$ shown in Theorem~\ref{prop:1} we
need to prove that the convergence holds in a stronger topology than
the vague topology described in Section~\ref{sec:point-processes}.
Indeed, the sum of points of a point process $\mu$ can be written as
\begin{displaymath}
  \int g d\mu=\int g(x)\,\mu(dx)=\sum_{x\in\supp\mu} g(x)\,,
\end{displaymath}
where $g(x)=x$ is a continuous function $\KK\mapsto\KK$ whose support
is neither bounded nor separated from the origin. Therefore, the weak
convergence $\mu_n\wto\mu$ does not imply the weak convergence of the
integrals $\int gd\mu_n$ to $\int gd\mu$.  

All point processes in this section are assumed to belong to
$\sMzero$, \ie they have almost surely at most a finite number of
points outside $B_r$ for every $r>0$.  Note that the sum and integrals
of $\KK$-valued functions are understood with respect to the addition
operation in $\KK$ and their absolute convergence with respect to the
norm on $\KK$.

\begin{lemma}
  \label{lem:gsup}
  Assume that $\KK$ possesses a sub-invariant norm.  Let
  $\mu_n\wto\mu$ for a point process $\mu$ such that
  $\mu(\cset)<\infty$ \as for all $r>0$. Let
  $g:\KK\mapsto\KK$ be a continuous function such that $\int gd\mu$
  converges absolutely and
  \begin{equation}
    \label{eq:unb}
    \limsup_{n}\P\left\{
      \left \|\int_{B_r} gd\mu_n\right\|\geq\eps\right\} \to 0\quad
    \mathrm{as}\ r\downarrow 0 
  \end{equation}
  for each $\eps>0$.  Then $\int gd\mu_n$ weakly converges to $\int
  gd\mu$.
\end{lemma}
\begin{proof}
  Since the space $\sMzero$ is Polish, by the Skorohod theorem it is
  possible to define $\mu_n$ and $\mu$ on the same probability space
  so that $\mu_n\vto\mu$ almost surely. In particular,
  $\mu_n|^r\vto\mu|^r$ \as for any $r>0$ such that
  $\P\{\mu(r\SS)>0\}=0$.  By the triangle inequality,
  \begin{multline*}
    d\Bigl(\int gd\mu_n,\int gd\mu\Bigr)\\ \leq d\Bigl(\int
    gd\mu_n,\int gd\mu_n|^r\Bigr) + d\Bigl(\int gd\mu_n|^r,\int
    gd\mu|^r\Bigr) + d\Bigl(\int gd\mu|^r,\int gd\mu\Bigr)\,.
  \end{multline*}  
  By~(\ref{triangleq}), the first summand is at most $\|\int_{B_r}
  gd\mu_n\|$.  The second summand $\zeta_n(r)$ converges to zero \as
  as $n\ti$ for any $r>0$ by Lemma~\ref{lem:kappa1} and continuity of
  $g$.  The third summand $\gamma(r)$ converges to zero \as as
  $r\downarrow0$, since $\int gd\mu$ exists.  Thus,
  \begin{align*}
    \P\left\{d\Bigl(\int gd\mu_n,\int gd\mu\Bigr)\geq\eps\right\}
    &\leq \P\Bigl\{\Bigl\|\int_{B_r}
    gd\mu_n\Bigr\|\geq\eps-\zeta_n(r)-\gamma(r)\Bigr\}\\
    &\leq \P\Bigl\{\Bigl\|\int_{B_r}
    gd\mu_n\Bigr\|\geq\eps/2\Bigr\} +\P\{\zeta_n(r)+\gamma(r)\geq
    \eps/2\}\,.
  \end{align*}
  By (\ref{eq:unb}) the probability on the left can be made
  arbitrarily small by the choice of $n$ and~$r$.
\end{proof}

Now we prove that the result of Theorem~\ref{prop:1} can be
strengthened to show a stronger type of convergence if
$\alpha\in(0,1)$. The corresponding topology in the special case of
point processes in the Euclidean space was called the
$\delta$-topology in \cite{dav:egor05}.

\begin{Th}
  \label{prop:3}
  Assume that the norm on $\KK$ is sub-invariant.  If
  (\ref{eq:resnick}) holds with $\alpha\in(0,1)$, then $\int
  x\,\beta_n(dx)$ weakly converges to $\int x\,\Pi_\alpha(dx)$.
\end{Th}
\begin{proof}
  Let $\vec{\rho}=(\rho_{1}, \rho_{2},\dots)$ be the decreasing 
  sequence of the norms of the support points of $\Pi_\alpha$ and
  $\vec{\rho}_n=(\rho_{1,n}, \rho_{2,n}, \dots)$ be the non-increasing
  infinite sequence of the norms of the support points of the process
  $\beta_n$ filled with $0$'s starting from the index $n+1$.
  
  Notice that~(\ref{eq:resnick}) implies that
  \begin{displaymath}
    \lim_{n\ti} n\P\{\|\xi\|>rb_n\}=cr^{-\alpha}
  \end{displaymath}
  with a constant $c$, which implies that $\|\xi\|$ belongs to the
  domain of attraction of $\alpha$-stable one-sided law on $\R_+$. It
  is well known (\seg \cite[Lemma~1]{lep:wood:zin81}) %\cite{jan89,kham:nag02})
  that in this case the finite-dimensional distributions of the
  sequences $\vec{\rho}_n$ converge to those of the sequence
  $\vec{\rho}$.  Using the Skorohod theorem, one can define
  $\vec{\rho}_n$ and $\vec{\rho}$ on the same probability space so
  that $\rho_{k,n}\to\rho_k$ for all $k=1,2,\dots$ almost surely.
  
  Recall a well-known convergence criterion in the space $L^1(\mu)$ of
  positive functions integrable with respect to a $\sigma$-finite
  measure $\mu$: if $f_n\geq0$, $f_n\to f$ $\mu$-almost everywhere and
  $\int f_n\,d\mu\to\int f\,d\mu$ then $f_n\to f$ in $L^1(\mu)$.
  Taking $\mu=\sum_k \delta_{\{k\}}$ to be a counting measure on
  $\{1,2,\dots\}$, this translates into
  \begin{equation}\label{eq:l1}
    \sum_{k=1}^\infty |\rho_{k,n}-\rho_k|\to0\quad\text{as $n\ti$
      almost surely.}
  \end{equation}
  
  Since the intensity measure of $\Pi_\alpha$ is non-atomic, with
  probability~1 none of the $\rho_k$'s coincides with~$r$.
  Conditionally on this event, in view of~(\ref{eq:l1}) we have that
  $\Pi_\alpha(\cset)=\beta_n(\cset)=k_0$ for all sufficiently large
  $n$, whence
  \begin{multline*}
    J_n(r)=\Bigl\|\int_{B_r} x\,\beta_n(dx)\Bigr\|\leq
    \int_{B_r}\|x\|\,\beta_n(dx)=\sum_{k=k_0+1}^\infty \rho_{k,n}
    \to \sum_{k=k_0+1}^\infty\rho_k\\
    =\int_{B_r}\|x\|\,\Pi_\alpha(dx)=I(r)
  \end{multline*}
  as $n\ti$ with probability~1. Hence, $ \limsup_n J_n(r)\leq I(r)$
  \as Given arbitrary $\delta>0$, fix a small $r>0$ so that
  \begin{displaymath}
%    \P \Bigl\{\int_{B_r} \|x\|\,\Pi_\alpha(dx)> \eps\Bigr\}< \delta\,.
    \P\{I(r)> \eps\}< \delta\,.
  \end{displaymath}
  This is possible since the integral $\int\|x\|\,d\Pi_\alpha$ with
  $\alpha\in(0,1)$ converges a.\,s.  and thus in probability. Then by
  Fatou's lemma,
  \begin{multline*}
    \limsup_n \P\{J_n(r)>\eps\}\leq \P\{\limsup_n \one_{J_n(r)>\eps}\}\\
    =\P\{\limsup_n J_n(r)> \eps\}\leq \P\{I(r)>\eps\}<\delta\,.
  \end{multline*}
  Thus (\ref{eq:unb}) holds and an application of Lemma~\ref{lem:gsup}
  completes the proof.
\end{proof}

\subsection{Domains of attraction for $\sas$ laws}
\label{sec:doma-attr-sas}

A random $\KK$-valued element $\zeta$ is said to belong to the domain
of attraction of a $\sas$ random element $\xi_\alpha$ if, for a
sequence $\{\zeta_n,n\geq1\}$ of \iid copies of $\zeta$,
\begin{displaymath}
  b_n^{-1}(\zeta_1+\dots+\zeta_n) \wto \xi_\alpha\,,
\end{displaymath}
where $\{b_n,n\geq1\}$ is a sequence of positive normalising constants
and $\wto$ denotes the weak convergence of $\KK$-valued random
elements. The following result confirms that conventional sufficient
conditions on the domain of attraction of $\sas$ laws are also
applicable for distributions on convex cones. It is proved by an
immediate application of Theorem~\ref{prop:3}.

\begin{Th}
  \label{th:2}
  Assume that $\KK$ has a sub-invariant norm.  Let $\{\zeta_n\}$ be a
  sequence of independent copies of a random element $\zeta\in\EE$.
  If $\zeta$ satisfies~(\ref{eq:resnick}) with $\alpha\in(0,1)$, then
  $\zeta$ belongs to the domain of attraction of $\sas$ random element
  $\xi_\alpha$ with the spectral measure $\sigma$ and representation
  (\ref{eq:9}) with $z=\neutral$.
\end{Th}

As a simple example, one sees that $\eta\eps$ belongs to the domain of
attraction of $\sas$ law if $\eta$ is a $\sas$ positive random
variable with $\alpha\in(0,1)$ and $\eps$ is any random element with
values in $\SS$ and independent of $\eta$.

The condition (\ref{eq:resnick}) for the domain of attraction appears
also in \cite[Th.~7.11]{ara:gin} for $\KK$ being a type $p$-Rademacher
space. It also characterises the domains of attraction of multivariate
max-stable distributions, see \cite{kot:nad00}.

The following result shows that it is possible to deduce the LePage
representation of a stable random element from the existence of the
corresponding spectral measure.

\begin{Th}
  \label{thr:esm}
  Assume that $\KK$ has a sub-invariant norm. Let $\xi$ be a $\sas$
  random element in $\EE$, such that (\ref{eq:resnick}) holds with
  $\alpha\in(0,1)$, so that $\xi$ possesses the spectral measure
  $\sigma$.  Then $\xi$ admits the LePage representation given by
  (\ref{eq:9}) with $z=\neutral$.
\end{Th}
\begin{proof}
  By the stability property, $\xi$ coincides in distribution with
  $n^{-1/\alpha}(\xi_1+\cdots+\xi_n)$ for each $n\geq1$. The latter is
  the sum of support points for the binomial process $\beta_n$ (with
  all the points distributed as $\xi$), so by Theorem~\ref{prop:3} it
  converges to the sum of points of $\Pi_\alpha$ being the LePage
  series~(\ref{eq:9}). 
\end{proof}

The existence of the spectral measure for $\sas$ random elements in
separable Banach spaces is a well-known fact, see
\cite[p.~152]{ara:gin}. Together with Theorem~\ref{thr:esm}, this
provides an alternative way to deduce the LePage representation for
$\sas$ laws with $\alpha\in(0,1)$ in separable Banach spaces,
\cf \cite{rosi90}. In order to derive the existence of the LePage
representation for more general semigroups, we use the technique of
harmonic analysis on semigroups that is explained in the following
sections.

\section{Distributions of stable random elements}
\label{sec:distr-stable-rand}

\subsection{Characters on semigroups}
\label{sec:char-semigr}

Assume that the semigroup $\KK$ is equipped with \emph{involution},
\ie a continuous map $\star:\KK\mapsto\KK$ satisfying
$(x+y)^\star=y^\star+x^\star$ and $(x^\star)^\star=x$ for all
$x,y\in\KK$. Assume also that $(ax)^\star=ax^\star$ for every
$x\in\KK$ and $a>0$. Note that the involution can be the identical
map. 

If $\KK$ is a normed cone, we also assume that $\|x^\star\|=\|x\|$,
\ie the sphere $\SS$ is invariant with respect to the involution
operation.  However, in
Sections~\ref{sec:distr-stable-rand}--\ref{sec:lepage-seri-repr} we do
not require that $\KK$ is a normed cone unless stated otherwise.  It
is not even assumed that $\KK$ is a pointed cone, \ie $\KK$ possesses
the origin.

\begin{definition}
  \label{def:char}
  A function $\chi$ that maps $\KK$ into the unit disk $\DD$ on the
  complex plane is called a \emph{character} if $\chi(\neutral)=1$,
  $\chi(x+y)=\chi(x)\chi(y)$ and $\chi(x^\star)=\overline{\chi(x)}$
  (the complex conjugate of $\chi(x)$) for all $x,y\in\KK$.
\end{definition}

We often encounter the following three cases.
\begin{itemize}
\item If the involution is the identity, then the characters take real
  values from $[-1,1]$. 
\item If $\KK$ is an idempotent semigroup, \ie $x+x=x$ for each
  $x\in\KK$, then all characters take values $0$ or $1$.
\item If the involution $x^\star$ is the inverse element to $x$, \ie
  $x+x^\star=\neutral$ for all $x\in\KK$, then the characters take
  values in the unit complex circle $\TT=\{z:\; |z|=1\}$.
\end{itemize}
It is also possible that the characters take values in the whole unit
disk $\DD$, for instance, if $\KK$ is the semigroup of probability
measures with the convolution operation, where the characters are
given by characteristic functions, see Example~\ref{ex:rprobm}. Note
that one can also study possibly unbounded characters, that are often
called semicharacters, see \cite[Def.~4.2.1]{ber:c:r}.

The set $\hat\KK$ of all characters (with the pointwise multiplication
operation) is called the \emph{restricted dual semigroup} to $\KK$.
The character $\one$ (identically equal to one) is the neutral element
in $\hat\KK$. The involution on $\hat\KK$ is the complex conjugate,
\ie $\chi^\star=\overline\chi$.  The family of characters is endowed
with the topology of pointwise convergence.  Then the projection
$\rho_x:\ \chi\mapsto\chi(x)$ becomes a continuous (and therefore
measurable) function $\hat\KK\mapsto\DD$ with respect to the Borel
$\sigma$-algebra on $\hat\KK$ for each $x\in\KK$.  Note that
Definition~\ref{def:char} imposes no continuity (nor even
measurability) condition on the characters.

The multiplication by $a$ in $\KK$ induces the multiplication
operation $\chi\mapsto a\circ\chi$ on $\hat\KK$ given by
$(a\circ\chi)(x)=\chi(ax)$ for all $x\in\KK$.  Note that $2\circ\chi$
is not necessarily equal to $\chi^2$, since $\chi(2x)$ is not
necessarily equal to $\chi(x+x)=\chi(x)^2$ unless the second
distributivity law (\ref{eq:cone4}) holds.

\bigskip

A family $\KDU$ of characters is said to be a \emph{sub-semigroup} of
$\hat\KK$ if $\KDU$ contains the identity character $\one$ and is
closed with respect to pointwise multiplication, \ie
$\chi_1\chi_2\in\KDU$ for all $\chi_1,\chi_2\in\KDU$. A sub-semigroup
$\KDU$ is called a \emph{cone} (of characters) if it is also closed
with respect to multiplication by numbers $\circ$, \ie
$a\circ\chi\in\KDU$ for all $\chi\in\KDU$ and $a>0$.

\begin{definition}
  \label{def:separ}
  A sub-semigroup $\KDU\subset\hat\KK$ is called \emph{separating},
  if, for any two distinct elements $x,y\in\KK$, there exists
  $\chi\in\KDU$ such that $\chi(x)\neq\chi(y)$; $\KDU$ is called
  \emph{strictly separating} if, for any two distinct elements
  $x,y\in\KK$, there exists $\chi\in\KDU$ such that
  $\chi(x')\neq\chi(y')$ for all $x'$ and $y'$ from open
  neighbourhoods of $x$ and $y$, respectively.
\end{definition}

If the characters from $\KDU$ are continuous on $\KK$, then the strict
separation follows from the simple separation condition.  It is known
\cite[Th.~V.22.17]{hew:ros63} that every locally compact abelian group
possesses a separating family of continuous characters.  However, not
all semigroups possess a separating family of characters.  For
instance, if $x+x=y+y$ and $x+x+x=y+y+y$ for some $x\neq y$, then $x$
and $y$ cannot be separated by any character, since every character
$\chi$ necessarily satisfies $\chi(x)^2=\chi(y)^2$ and
$\chi(x)^3=\chi(y)^3$, see Example~\ref{ex:comp-non-conv}.

\subsection{Laplace transform}
\label{sec:laplace-transform}

A sub-semigroup $\KDU$ of characters generates a $\KDU$-weak topology
on $\KK$ by declaring $x_n\stackrel{w}\to x$ if and only if
$\chi(x_n)\to\chi(x)$ for all $\chi\in\KDU$. The $\KDU$-weak topology
is the weakest topology that makes all characters from $\KDU$
continuous. Let $\sF(\KK;\KDU)$ be the smallest $\sigma$-algebra on
$\KK$ that makes all $\chi\in\KDU$ measurable. This $\sigma$-algebra
is generated by the cylindrical sets $\{x\in\KK:\; \chi_i(x)\in F_i,\;
i=1,\dots,n\}$, $n\geq1$, where $F_1,\dots,F_n$ are Borel subsets of
$\DD$ and $\chi_1,\dots,\chi_n\in\KDU$.

Denote by
\begin{displaymath}
  \KDD=\hat{\KDU}
\end{displaymath}
the restricted dual semigroup to $\KDU$, \ie the family of all
characters on $\KDU$. We equip $\KDD$ with the topology of pointwise
convergence, which generates the corresponding Borel $\sigma$-algebra
$\sB(\KDD)$. The \emph{evaluation map} $\imath:\KK\mapsto\KDD$ is
defined by associating every $x\in\KK$ with $\rho=\rho_x\in\KDD$ such
that $\rho_x(\chi)=\chi(x)$ for all $\chi\in\KDU$,
\cf \cite[Sec.~20]{kel:nam63}.  The evaluation map $\imath$ is
injective if and only if $\KDU$ is separating.

The \emph{Laplace transform} of a $\KK$-valued random element $\xi$ is
a complex-valued function $\chi\mapsto\E\chi(\xi)$, where $\chi$ is a
Borel measurable character from $\hat\KK$.  The following result is
well known for random elements in locally compact spaces with
continuous characters, see \cite[\S~IX.5.7]{bouVI}.  However, it also
holds in a more general framework.

\begin{Th}
  \label{thr:laplace}
  If $\hat\KK$ has a separating sub-semigroup $\KDU$ such that
  $\sF(\KK;\KDU)=\sB(\KK)$, then every probability measure $\P$
  on $\KK$ (or the distribution of a random element $\xi$) is uniquely
  determined by its Laplace transform
  \begin{equation}
    \label{eq:lap-transform}
    \E\chi(\xi)=\int_{\KK} \chi(x)\P(dx)\,, \quad \chi\in\KDU\,.
  \end{equation}
\end{Th}
\begin{proof}
  The function $\E\chi(\xi)$ is a bounded positive definite function
  on $\KDU$. By \cite[Th.~4.2.8]{ber:c:r}, there exists one and
  only one Radon measure $\mu$ on the restricted dual semigroup to
  $\KDU$ (\ie the second dual semigroup $\KDD$), such that
  \begin{displaymath}
    \E\chi(\xi)=\int_{\KDD} \rho(\chi)\mu(d\rho)\,.
  \end{displaymath}
  The separation condition implies that the evaluation map $\imath$ is
  injective. By the condition, $\imath$ is $\sB(\KK)$-measurable,
  since, for Borel sets $F_1,\dots,F_n\subset\DD$,
  \begin{multline*}
    \imath^{-1}(\{\rho\in\KDD:\;\rho(\chi_i)\in F_i,\; i=1,\dots,n\})\\
    =\{x\in\KK:\; \chi_i(x)\in F_i,\; i=1,\dots,n\}
    \in \sF(\KK;\KDU)=\sB(\KK)\,.
  \end{multline*}  
  This makes it possible to define the image measure $\mu'$ of $\P$
  under the natural map $\imath$ by $\mu'(F)=\P(\imath^{-1}(F))$ for
  all Borel $F\subset\KDD$. By definition, $\mu'$ is supported by
  $\imath(\KK)$.  After substitution $\rho=\imath(x)$
  \eqref{eq:lap-transform} can be re-written as
  \begin{displaymath}
    \int_{\KDD} \chi(\imath^{-1}(\rho))\,\mu'(d\rho)=\int_{\KDD}
    \rho(\chi)\,\mu'(d\rho) 
  \end{displaymath}
  and the uniqueness property implies that $\mu=\mu'$.
\end{proof}

\begin{remark}
  \label{rem:fffbbb}
  Note that $\sF(\KK;\KDU)=\sB(\KK)$ if $\KDU$ is a separating
  \emph{countable} family of $\sB(\KK)$-measurable functions, see
  \cite[Prop.~I.1.4]{vak:tar:cho87}. If $\KDU$ consists of
  \emph{continuous} characters, then the separation condition already
  implies that $\sF(\KK;\KDU)=\sB(\KK)$, see
  \cite[Th.~I.1.2]{vak:tar:cho87}. If $\KK$ is a linear space,
  Theorem~\ref{thr:laplace} turns into
  \cite[Th.~IV.2.2]{vak:tar:cho87}.
\end{remark}

\begin{remark}
  \label{rem:no-suff}
  It should be noted that Theorem~\ref{thr:laplace} is the pure
  uniqueness result. It does not assert that a positive definite
  function on $\KDU$ is the Laplace transform of a certain random
  element, \ie an analogue of the Bochner theorem may not hold, see
  \cite{bis00}.
\end{remark}

\begin{example}
  \label{ex:vee}
  For the cone $(\R_+,\vee)$, the collection of indicator functions
  $\one_{[0,a]}$ with $a\geq0$ may be taken as a separating family
  $\KDU$.  Then $\chi(x_n)\to\chi(x)$ if and only if $x_n\uparrow x$
  as $n\ti$, \ie the corresponding weakly-open sets in $\R_+$ are
  $(a,b]$ with $a<b$.  They generate the same $\sigma$-algebra as the
  metric Borel $\sigma$-algebra $\sB(\R_+)$, so that
  Theorem~\ref{thr:laplace} applies. The Laplace transform of $\xi$ is
  its cumulative distribution function $\P\{\xi\leq a\}$, $a\geq0$.
\end{example}

A $\sas$ random element $\xi$ in $\KK$ is necessarily \emph{infinitely
  divisible}. The Laplace transform of $\xi$ is a \emph{positive definite}
infinitely divisible function of $\chi$, \ie $(\E\chi(\xi))^{1/n}$ is
a positive definite function of $\chi\in\KDU$ for every
$n\geq1$.  The results on infinitely divisible functions in semigroups
\cite[Th.~3.2.2, Prop.~4.3.1]{ber:c:r} imply that
\begin{equation}
  \label{eq:echixi=exp-phichi-}
  \E\chi(\xi)=\exp\{-\phi(\chi)\}\,, \quad \chi\in\KDU\,,
\end{equation}
where $\phi$ is a \emph{negative definite} complex-valued function on
$\KK$ with $\Re\phi\in[0,\infty]$ and $\phi(\one)=0$. We call $\phi$
the \emph{Laplace exponent} of $\xi$. The dominated convergence
theorem implies that $\phi(\chi)$ is continuous with respect to
pointwise convergence in $\KDU$.

\begin{definition}
  \label{def:idemp}
  A random element $\xi$ is \emph{idempotent} if $\xi$ coincides in
  distribution with the sum of its two \iid copies, \ie
  $\xi\deq\xi_1+\xi_2$, see \cite[p.~41]{gren63}.
\end{definition}

It is easy to see from the definition that $\xi$ is idempotent if and
only if its Laplace transform assumes only values $0$ and $1$.  The
deterministic $\xi$ being equal to an idempotent element of $\KK$
clearly has an idempotent distribution.  It is also the case if $\xi$
is distributed according to a finite Haar measure on any subgroup of
$\KK$.  In a sense, an idempotent random element is $\sas$ with
$\alpha=\infty$.  The following two useful results show that in some
cases non-trivial $\sas$ random elements with a finite $\alpha$ cannot
be idempotent.

\begin{lemma}
  \label{lem:idem}
  Assume that $\KK$ is a normed cone.  If $\xi$ is an idempotent
  $\sas$ random element with a finite $\alpha$ and \as finite norm,
  then $\xi=\zero$ \as
\end{lemma}
\begin{proof}
  If $\xi$ is idempotent and stable, then (\ref{eq:a1alph-a+b1-})
  implies that $\xi\deq 2^{1/\alpha}\xi$, so that $\|\xi\|=0$ \as
  Note that $\xi=\neutral$ \as with $\neutral\neq\zero$ is
  impossible, since $\|\neutral\|=\infty$ if $\neutral\neq\zero$.
\end{proof}

\begin{lemma}
  \label{lem:c2}
  Assume that $\KK$ possesses a strictly separating family $\KDU$ of
  characters with values in $[0,1]$. Then each idempotent $\sas$
  random element $\xi$ with $\alpha\neq\infty$ is necessarily a
  deterministic idempotent element in $\KK$.
\end{lemma}
\begin{proof}
  Assume that the support of $\xi$ contains at least two distinct
  points and consider a character $\chi$ that strictly separates their
  neighbourhoods. Then $\E\chi(\xi)$ is distinct from $0$ and $1$,
  which is impossible for an idempotent $\xi$. Hence $\xi$ is
  deterministic.
\end{proof} 

\subsection{Characterisation of Laplace transforms for $\sas$
  elements}
\label{sec:char-lapl-transf}

It should be noted that many semigroups do not possess any non-trivial
continuous character.  For example, the only continuous character in
$(\R_+,\vee)$ is the one identically equal to $1$, while non-trivial
characters $\one_{[0,t)}$ and $\one_{[0,t]}$ are only semicontinuous.
The following assumption imposes a weaker form of the continuity
property on the characters. For this, note that any character
$\chi\in\hat\KK$ can be decomposed as the product
\begin{equation}
  \label{eq:chi-prod}
  \chi=\chi'\chi''
\end{equation}
of the $[0,1]$-valued character $\chi'$ (the modulus of $\chi$) and
the $\TT$-valued character $\chi''$ (corresponding to the exponent of
the argument of $\chi$).

\begin{description}
\item[(C)] There exists a cone of characters $\KDU$ such that every
  $\chi\in\KDU$ has semicontinuous modulus and continuous argument and
  $\KDU$ possesses a strictly separating countable sub-family. 
\end{description}

By Remark~\ref{rem:fffbbb}, if \textbf{(C)} is satisfied, then
Theorem~\ref{thr:laplace} holds, \ie the Laplace transform on $\KDU$
uniquely determines the distribution of a $\KK$-valued random element.
Two particular important cases are
\begin{itemize}
\item $\KDU$ is a separating cone of continuous $\DD$-valued
  characters. Then it is automatically strictly separating and also
  has a countable separating sub-family, see
  \cite[Th.~I.1.2]{vak:tar:cho87}.
\item $\KDU$ is a countable strictly separating cone of $[0,1]$-valued
  semicontinuous characters.
\end{itemize}
Probability distributions in a topological linear space with a
separating family of continuous characters have been studied in
\cite{buld80}.

\medskip

A non-degenerate random element $\xi$ in a cone may have a
\emph{self-similar} distribution, \ie $s\xi$ may coincide in
distribution with $\xi$ for some (or even all) $s>0$. For instance,
this is the case if $\KK$ is the cone of closed sets in $\R_+$ and
$\xi$ is the set of zeros of the Wiener process.  However, this is
clearly impossible if the norm of $\xi$ is \as finite, since then the
non-negative finite random variable $\|\xi\|$ would coincide in
distribution with $s\|\xi\|$ for $s\neq1$. The following result is
applicable without assuming that the norm of $\xi$ is finite and even
without assuming that $\KK$ is a normed cone. Its proof follows the
scheme used in \cite[Th.~4.1.12]{mo1}.

\begin{Th}
  \label{th:s-sim}
  Assume that Condition \textbf{(C)} holds. If $\xi$ is a
  non-idempotent $\sas$ random element in $\KK$ such that, for some
  $s>0$,
  \begin{equation}
    \label{eq:scal}
    \E\chi'(s\xi)=\E\chi'(\xi)\,,\qquad \E\chi''(s\xi)=\E\chi''(\xi)
  \end{equation}
  for all $\chi=\chi'\chi''\in\KDU$, then $s=1$.
\end{Th}
\begin{proof}
  The non-idempotency condition implies that $\alpha$ is finite.
  Assume that (\ref{eq:scal}) holds for $s>1$. The case $s<1$ is the
  same, since it is just a reformulation of \eqref{eq:scal} for
  $s\circ\chi$. The $\sas$ condition (\ref{eq:a1alph-a+b1-})
  immediately implies that $\xi_1+\cdots+\xi_n\deq a_n\xi$ for some
  $a_n>0$ and all $n\geq1$.  Writing $a_n=\delta_n s^{k(n)}$ in powers
  of $s$ with $\delta_n\in[1,s)$ we obtain that
  \begin{displaymath}
    \chi(\xi_1)\cdots\chi(\xi_n) \deq \chi(\delta_n s^{k_n}\xi)\,.
  \end{displaymath}
  By taking the absolute values we deduce that
  $\chi'(\xi_1)\cdots\chi'(\xi_n)$ coincides in distribution with
  $\chi'(\delta_n s^{k_n}\xi)$. Now (\ref{eq:scal}) applied to
  $\delta_n\circ\chi'$ implies that
  \begin{displaymath}
    (\E\chi'(\xi))^n=\E\chi'(\delta_n\xi)\,, \quad n\geq1\,,
  \end{displaymath}
  and the similar identity holds for $\chi''$. 
  
  Without loss of generality assume that $\delta_n\to\delta$ as
  $n\ti$, otherwise, consider a convergent subsequence.  Assume that
  $\chi'$ is lower semicontinuous. If $\E\chi'(\xi)<1$, then
  \begin{equation}\label{eq:semi}
    \E\chi'(\delta\xi)\leq\liminf_{n\ti}
    \E\chi'(\delta_n\xi)=\liminf_{n\ti} \bigl(\E\chi'(\xi)\bigr)^n=0\,. 
  \end{equation}
  If $\E\chi'(\delta\xi)<1$, then the strict separation condition
  implies that $\E\chi'(\delta_n\xi)<1$ for sufficiently large $n$.
  Hence $\E\chi'(\xi)<1$, and the above argument yields that
  $\E\chi'(\delta\xi)=0$. Thus, $\E\chi'(\delta\xi)$ assumes only the
  values $0$ or $1$.
  
  If $\chi'$ is upper semicontinuous, a similar argument applies to
  $\E\chi'(\delta_n^{-1}\xi)=(\E\chi'(\xi))^{1/n}$. Then
  \begin{displaymath}
    \E\chi'(\delta^{-1}\xi)\geq \limsup_{n\ti} \E\chi'(\delta_n^{-1}\xi)
    =\limsup_{n\ti}(\E\chi'(\xi))^{1/n}\,.
  \end{displaymath}
  If $\E\chi'(\xi)>0$, then $\E\chi'(\delta^{-1}\xi)=1$. Furthermore, if
  $\E\chi'(\delta^{-1}\xi)>0$, then $\E\chi'(\delta_n^{-1}\xi)>0$ for
  sufficiently large $n$, whence $\E\chi'(\xi)>0$ and consequently
  $\E\chi'(\delta^{-1}\xi)=1$.
  
  Now consider the character $\chi''$. By continuity,
  $\E\chi''(\delta_n\xi)\to \E\chi''(\delta\xi)$. Thus taking limit in
  $(\E\chi''(\delta_n\xi))^2=\E\chi''(\delta_{2n}\xi)$, we obtain that
  $\E\chi''(\delta\xi)$ is either $0$ or $1$. 
  
  If $\E\chi'(\xi)=0$, then $\chi'(\xi)=0$ a.s., so that $\chi(\xi)=0$
  a.s. and $\E\chi(\xi)=0$. If $\E\chi'(\xi)=1$, then $\chi'(\xi)=1$
  a.s., so that $\E\chi(\xi)=\E\chi''(\xi)$, which is either zero or
  one also in this case.  Thus, $\E\chi(\xi)$ assumes only values $0$
  or $1$ contradicting the non-idempotency assumption.
\end{proof}

\begin{Cor}
  \label{cor:ss}
  Assume that \textbf{(C)} holds. If a non-idempotent random element
  $\xi$ satisfies $s\xi\deq\xi$ for some $s>0$, then $s=1$. 
\end{Cor}

\begin{definition}
  A $\KK$-valued random element is said to be \emph{non-trivial} if
  $\P\{\xi=\neutral\}<1$. A non-trivial $\KK$-valued random element
  $\xi$ is said to be \emph{proper} if its Laplace transform does not
  vanish for any $\chi$ from the separating cone $\KDU$.
\end{definition}

If the characters from $\KDU$ take values from $(0,1]$, then all
non-trivial random elements are proper.  If $\xi$ has an idempotent
factor, \ie $\xi$ can be represented as a sum of an idempotent (not
identically equal to $\neutral$) random element and another (possibly
deterministic) random element, then $\xi$ is not proper. The inverse
implication holds for $\KK$ being a locally compact group, where
proper random elements can be characterised as those which do not
possess idempotent factors, see \cite[Th.~IV.4.2]{par67}.
In the studies of random sets \cite[Sec.~4.1.2]{mo1} a vanishing
Laplace transform corresponds to the random set with fixed points, see
Example~\ref{ex:usrs}. Below we summarise a few immediate properties
of proper elements.

\begin{lemma}
  \label{lemma:misc}
  Assume that conditions of Theorem~\ref{thr:laplace} are satisfied.
  \begin{itemize}
  \item [\textsl{(i)}] A proper random element $\xi$ shares the same
    distribution with $\xi+c$ for a deterministic $c$ if and only if
    $c=\neutral$.
  \item [\textsl{(ii)}] A proper $\sas$ random element $\xi$ satisfies
    $\xi\neq\zero$ \as in case $\neutral\neq\zero$.
  \item [\textsl{(iii)}] If $\xi$ and $\eta$ are two independent
    elements such that $\xi+\eta$ is proper, then both $\xi$ and
    $\eta$ are proper. If, additionally, $\xi$ and $\xi+\eta$ are
    $\sas$, then $\eta$ is $\sas$.
  \end{itemize}
\end{lemma}
\begin{proof}
  \textsl{(i)} It suffices to note that
  $\E\chi(\xi)=(\E\chi(\xi))\chi(c)$, so that $\chi(c)=1$ for all
  $\chi$. 
  
  \textsl{(ii)} By the separation property, $\neutral$ and $\zero$ are
  separated by a certain character $\chi\in\KDU$. Since both $\zero$
  and $\neutral$ are idempotent and $\chi(\neutral)=1$, we necessarily
  have $\chi(\zero)=0$. Then $\delta=\P\{\xi=\zero\}\leq
  \P\{\chi(a\xi)=0\}$ for all $a>0$.  The stability property implies
  that
  \begin{displaymath}
    \chi(a_n\xi_1)\cdots \chi(a_n\xi_n)\deq\chi(\xi)
  \end{displaymath}
  for certain $a_n>0$ and each $n\geq1$.  Then $\chi(\xi)=0$ when at
  least one of $\chi(a_n\xi_1),\dots,\chi(a_n\xi_n)$ vanishes. Since
  this happens with probability at least $\delta>0$ for each of these
  independent factors and every $n\geq1$, $\chi(\xi)=0$ with
  probability $1$, so that the Laplace transform of $\xi$ vanishes,
  \ie $\xi$ is not proper.
  
  \textsl{(iii)} Since $\E\chi(\xi+\eta)=\E\chi(\xi)\E\chi(\eta)$,
  both $\xi$ and $\eta$ are proper. Furthermore,\\
  $\E(c\circ\chi)(\xi)=\E\chi(c\xi)$ implies that any proper element
  multiplied by a constant is also proper. The definition of stability
  of $\xi+\eta$ implies that
  \begin{displaymath}
    \E\chi(a^{1/\alpha}\xi_1+b^{1/\alpha}\xi_2)\,
    \E\chi(a^{1/\alpha}\eta_1+b^{1/\alpha}\eta_2)=  
    \E\chi((a+b)^{1/\alpha}\xi)\,\E\chi((a+b)^{1/\alpha}\eta)\,.
  \end{displaymath}
  After cancelling the non-zero equal terms related to $\sas$ element
  $\xi$, we get that
  \begin{displaymath}
    \E\chi(a^{1/\alpha}\eta_1+b^{1/\alpha}\eta_2)= 
    \E\chi((a+b)^{1/\alpha}\eta)
  \end{displaymath}
  for all $\chi$ which, in view of Theorem~\ref{thr:laplace}, is
  equivalent to the fact that $\eta$ is $\sas$.
\end{proof}

\begin{Cor}
  \label{cor:decomp}
  Let $\xi$ be a proper $\sas$ with representation~(\ref{eq:9}) in a
  normed cone $\KK$.  Assume that the conditions of
  Theorem~\ref{thr:laplace} are satisfied.
  \begin{itemize}
  \item [\textsl{(i)}] Then $z$ is not an idempotent element of $\KK$
    unless $z=\neutral$. In particular, $z\neq\zero$ in case
    $\zero\neq\neutral$.
  \item [\textsl{(ii)}] If, in addition, the norm is sub-invariant and
    $0<\alpha<1$, then $z=\neutral$.
  \end{itemize}
\end{Cor}
\begin{proof}
  \textsl{(i)}. Since $\chi(z)=0$ or 1 for any idempotent element,
  then $z+\xi$ is proper if and only if $\chi(z)=1$ for all $\chi$.
  Then $z=\neutral$ by the separation condition.

  \textsl{(ii)}. By Lemma~\ref{lemma:misc}\textsl{(iii)},
  $z\in\KK(\alpha)$, but $\KK(\alpha)=\{\neutral\}$ according to
  Lemma~\ref{lemma:sub-invar}. 
\end{proof}

The following result is the key characterisation theorem for Laplace
transforms of proper $\sas$ random elements. It also establishes the
equivalence of the $\sas$ property (\ref{eq:a1alph-a+b1-}) and its
``discrete'' variant (\ref{eq:plus}).

\begin{Th}
  \label{th:rep}
  Let $\xi$ be a proper random element in $\KK$. Assume that
  \textbf{(C)} holds.  Then the following conditions are equivalent.
  \begin{itemize}
  \item[\textsl{(i)}] For every $n\geq2$ there exists a positive
    constant $a_n\neq 1$ such that
    \begin{equation}
      \label{eq:plus}
      \xi_1+\cdots+\xi_n\deq a_n\xi\,,
    \end{equation}
    where $\xi_1,\dots,\xi_n$ are \iid copies of $\xi$.
  \item[\textsl{(ii)}] $\xi$ is $\sas$ with finite $\alpha$.
  \item[\textsl{(iii)}] The Laplace transform of $\xi$ is given by
    (\ref{eq:echixi=exp-phichi-}), where $\phi$ satisfies
    \begin{equation}
      \label{eq:hom}
      \phi(s\circ\chi)=s^\alpha\phi(\chi)
    \end{equation}
    for all $s>0$.
  \end{itemize}
\end{Th}
\begin{proof}
  Note that \textsl{(ii)} immediately implies \textsl{(i)} and
  \textsl{(iii)} implies \textsl{(ii)} by Theorem~\ref{thr:laplace}
  and (\ref{eq:echixi=exp-phichi-}).
  
  It remains to prove that \textsl{(i)} implies \textsl{(iii)}. First,
  \textsl{(i)} implies that $\xi$ is not idempotent. Furthermore,
  (\ref{eq:plus}) yields that $n^{-1}\phi(a_n\circ\chi)=\phi(\chi)$.
  Define $a(s)=a_n/a_m$ for $s=n/m$. It is easy to see that $a(s)$
  does not depend on the representation of its rational argument~$s$,
  \begin{equation}
    \label{eq:aphi}
    s^{-1}\phi(a(s)\circ\chi)=\phi(\chi)
  \end{equation}
  and $\phi(a(s)a(s_1)\circ\chi)=\phi(a(ss_1)\circ\chi)$.  By
  Theorem~\ref{th:s-sim}, $a(ss_1)=a(s)a(s_1)$ for all rational
  $s,s_1>0$. To prove that $a(s)=s^{1/\alpha}$ for some $\alpha>0$, it
  is now sufficient to show that $a(s)$ is continuous on the set of
  positive rational numbers. 
  
  Let $\chi=\chi'\chi''$ for $[0,1]$-valued character $\chi'$ and
  $\TT$-valued character $\chi''$. 
  Since $\chi'$ and $\chi''$ are characters themselves,
  (\ref{eq:aphi}) holds for $\chi'$ and $\chi''$ separately.  Thus, if
  $s_n$ is a sequence of rational numbers that converges to 1, then
  $\E\chi'(a(s_n)\xi)\to\E\chi'(\xi)$ and
  $\E\chi''(a(s_n)\xi)\to\E\chi''(\xi)$.
  
  Assume that $a(s_n)$ has a finite positive limit $a$. Since $\chi''$
  is continuous, $\E\chi''(a\xi)=\lim
  \E\chi''(a(s_n)\xi)=\E\chi''(\xi)$ by the dominated convergence. If
  $\chi'$ is lower semicontinuous, then by Fatou's lemma,
  \begin{displaymath}
    \E\chi'(a\xi)\leq \liminf\E\chi'(a(s_n)\xi)=\E\chi'(\xi)\,.
  \end{displaymath}
  Furthermore, (\ref{eq:aphi}) written for $a(s_n)^{-1}$ yields that
  \begin{displaymath}
    \E\chi'(a^{-1}\xi)\leq \liminf \E\chi'(a(s_n)^{-1}\xi)=\E\chi'(\xi)\,.
  \end{displaymath}
  This also holds for $a\circ\chi'$, leading to $\E\chi'(\xi)\leq
  \E\chi'(a\xi)$. Hence $\E\chi'(a\xi)=\E\chi'(\xi)$, \ie
  (\ref{eq:scal}) holds, so that $a=1$ by Theorem~\ref{th:s-sim}.
  Similar arguments apply if $\chi'$ is an upper semicontinuous
  characters.
  
  To finish the proof it suffices to consider either the case of
  $a(s_n)\to0$ or $a(s_n)\ti$.  Let $a_m>1$ in (\ref{eq:plus}) for
  some $m>1$. Assume that $a(s_n)\to0$ as $n\ti$. Then, for any
  $n\geq1$, $a(s_n)=(a_m)^{-k(n)}\delta_n$, where $\delta_n\in[1,a_m)$
  and $k(n)\ti$. Note that $(a_m)^{k(n)}=a_{m^{k(n)}}$.  Hence
  \begin{displaymath}
    \delta_n\xi=(a_m)^{k(n)}a(s_n)\xi
    \deq a(s_n)\xi_1+\cdots+a(s_n)\xi_{m^{k(n)}}\,.
  \end{displaymath}
  Therefore, 
  \begin{equation}
    \label{eq:delkm}
    \E\chi(\delta_n\xi)=\left(\E\chi(a(s_n)\xi)\right)^{m^{k(n)}}
    =\left(\E\chi(\xi)\right)^{s_nm^{k(n)}}\,,
  \end{equation}
  where $\chi$ stands for $\chi'$ or $\chi''$.  Taking if necessary a
  convergent subsequence, assume that $\delta_n\to\delta$ as $n\ti$.
  Assume that $\chi'$ is lower semicontinuous. If
  $\E\chi'(\delta\xi)<1$, then the strict separation condition implies
  that $\E\chi'(\delta_n\xi)<1$ for sufficiently large $n$. By
  (\ref{eq:delkm}), $\E\chi'(\xi)<1$, which, in turn, implies that
  $\E\chi'(\delta_n\xi)\to0$. By Fatou's lemma,
  \begin{displaymath}
    \E \chi'(\delta\xi)\leq \liminf \E\chi'(\delta_n\xi)=0\,.
  \end{displaymath}
  Thus, $\E\chi'(\delta\xi)=0$ is either zero or one. For $\chi''$,
  (\ref{eq:delkm}) implies that $\E\chi''(\delta\xi)$ is either zero
  or one. Similarly to the argument used in the proof of
  Theorem~\ref{th:s-sim}, we conclude that $\delta\xi$ is idempotent,
  whence $\xi$ is idempotent too and so cannot be proper.
  
  If $a_m<1$ for all $m\geq1$, then the assumption $a(s_n)\ti$ leads
  to a contradiction in the same manner.  Therefore, $a(s_n)\to1$ as
  $s_n\to1$, whence $a(s)$ is continuous on the rational numbers. Thus
  $a(s)=s^{1/\alpha}$ for all $s>0$ and (\ref{eq:hom})
  follows.
\end{proof}

\begin{remark}
  In general, (\ref{eq:plus}) either holds with $a_n=1$ for all $n$
  (meaning that $\xi$ is idempotent) or $a_n\neq1$ for all $n\geq2$.
  Indeed, if $a_n=1$ for some $n\geq2$, then $a_{2n-1}=1$, and further
  $a_2=1$ implying that $a_n=1$ for all $n$.
\end{remark}

Theorem~\ref{th:rep} can also be extended to any non-idempotent $\xi$,
in which case the both sides of (\ref{eq:hom}) are allowed to be
infinite. From now on, we mainly consider proper random elements.

\subsection{Possible values of $\alpha$ for $\sas$ distributions}
\label{sec:parameter-ranges-sas}

In the sequel we make use of the following condition:
\begin{description}
\item[(E)] $\KK$ is a pointed cone such that $\chi(sx)\to1$ for all
  $\chi\in\KDU$ and $x\in\EE$, where $s\downarrow0$ if
  $\neutral=\zero$ or $s\ti$ if $\neutral\neq\zero$.
\end{description}
Condition \textbf{(E)} means that $sx$ $\KDU$-weakly converges to
$\neutral$ for all $x\in\EE$. It clearly holds if
$\zero=\neutral$ and the characters are continuous. Conversely, if
\textbf{(E)} holds and $\KK$ is a group, then $y_n\to y$ implies
$y_n-y\to\neutral$, so that
$\chi(y_n)\chi(-y)=\chi(y_n)\chi^{-1}(y)\to1$. Thus all characters are
continuous, so that \textbf{(C)} automatically holds given that the
characters form a separating family.  Condition \textbf{(E)} has a
further implication on the characterisation of proper random elements
in $\KK$.

\begin{lemma}
  \label{lemma:proper}
  If Condition~\textbf{(E)} holds, then every $\sas$ random element in
  $\EE$ necessarily has a proper distribution.
\end{lemma}
\begin{proof}
  Depending on the sign of $\alpha$ and whether $\|\neutral\|=0$ or
  $\|\neutral\|=\infty$, define $\eta_n=n^{1/\alpha}\xi$ or
  $\eta_n=n^{-1/\alpha}\xi$ so that $\|\eta_n\|\to\|\neutral\|$ a.\,s.
  Let $\E\chi(\xi)=0$. The stability property implies that
  $\E\chi(\eta_n)=0$. Now the dominated convergence theorem leads to a
  contradiction with the fact that $\E\chi(\eta_n)\to
  \E\chi(\neutral)=1$.
\end{proof}

The following important result shows that the relationship between the
origin and the neutral element has a crucial influence on the range of
possible values of the stability parameter $\alpha$.

\begin{Th}
  \label{th:alpha}
  Assume that \textbf{(C)} and \textbf{(E)} hold. Then, for every
  proper $\sas$ random element $\xi$, its characteristic exponent
  $\alpha$ is positive if and only if $\neutral=\zero$.
\end{Th}
\begin{proof}
  Assume that $\neutral=\zero$. Then $sx$ converges to $\neutral$ as
  $s\downarrow0$. Since $(s\circ\chi)(x)=\chi(sx)$, \textbf{(E)}
  ensures that $s\circ\chi\to\one$ pointwise.  The continuity of the
  Laplace exponent $\phi$ implies that the left-hand side of
  (\ref{eq:hom}) converges to $\phi(\one)=0$, so that $s^\alpha\to0$
  as $s\to0$. Thus $\alpha>0$.
  
  Assume that $\neutral\neq\zero$. By \textbf{(E)}, $\chi(sx)\to1$
  as $s\ti$ for all $x\neq\zero$. Lemma~\ref{lemma:misc}\textsl{(ii)}
  implies that $\xi\neq\zero$ \as Thus, $\E\chi(s\xi)\to1$, whence
  the right-hand side of (\ref{eq:hom}) converges to zero as $s\ti$,
  so that $\alpha<0$.
\end{proof}

Theorem~\ref{th:alpha} does not hold without assuming that $\xi$ is
proper, see Example~\ref{ex:usrs}. Neither it holds without assuming
that $\KK$ is a pointed cone, see Example~\ref{ex:isrm}. By
Lemma~\ref{lemma:proper}, it is possible to omit the requirement that
$\xi$ is proper from Theorem~\ref{th:alpha} if $\xi$ belongs to $\EE$
almost surely.

\medskip

Since the function $\phi$ from (\ref{eq:echixi=exp-phichi-}) is
negative definite, the properties of negative definite functions
together with Theorem~\ref{th:rep} yield the following result.

\begin{Th}
  \label{thr:bound}
  Assume that \textbf{(C)} and the second distributivity
  law (\ref{eq:cone4}) hold. Then
  \begin{itemize}
  \item[\textsl{(i)}] every proper $\sas$ random element $\xi$ has
    parameter $\alpha\leq 2$;
  \item[\textsl{(ii)}] if the involution is identical, then
    $\alpha\leq 1$.
  \end{itemize}
  If, in addition, \textbf{(E)} holds, then $\alpha>0$ in \textsl{(i)} and \textsl{(ii)}.
\end{Th}
\begin{proof}
  \textsl{(i)} It follows from \cite[Prop.~4.3.3]{ber:c:r} that
  \begin{displaymath}
    \sqrt{|\phi(\chi_1\chi_2)|}\leq
    \sqrt{|\phi(\chi_1)|}+\sqrt{|\phi(\chi_2)|}\,.
  \end{displaymath}
  By letting $\chi_1=\chi_2=\chi$, we see that
  \begin{displaymath}
    |\phi(\chi^2)|\leq 4|\phi(\chi)|\,.
  \end{displaymath}
  The left-hand side is equal to $|\phi(2\circ\chi)|$ by
  (\ref{eq:cone4}), so that (\ref{eq:hom}) implies that $2^\alpha\leq
  4$.
  
  \textsl{(ii)} Since $\phi$ is negative definite,
  $\phi(\chi_1)+\phi(\chi_2)-\phi(\chi_1\chi_2)$ is a positive
  definite kernel, see \cite[Prop.~4.1.9]{ber:c:r}. In particular, its
  value is non-negative if $\chi_1=\chi_2=\chi$, whence
  $\phi(2\circ\chi)=\phi(\chi^2)\leq 2\phi(\chi)$, \ie $\alpha\leq1$.

  Finally, note that Lemma~\ref{lemma:n=0} implies that
  $\neutral=\zero$, thus $\alpha>0$ by Theorem~\ref{th:alpha}
  given that \textbf{(E)} holds.
\end{proof}

If the second distributivity (\ref{eq:cone4}) does not hold, then
$\alpha$ may have various ranges of possible values.  For instance, in
$(\R_+,\vee)$ any $\alpha>0$ and for $([0,\infty],\min)$ any
$\alpha<0$ are possible, see Example~\ref{ex:rvee}. The following
example shows that it is possible to define a cone where $\sas$ laws
exist with any $\alpha$ from a given interval $(0,\beta)$ or
$(-\beta,0)$.

\begin{example}
  \label{ex:range}
  Consider the cone $\R_+$ with the addition operation given by
  $(x^\beta+y^\beta)^{1/\beta}$ with $\beta>0$ and the conventional
  multiplication by positive numbers. Then $x+x=2^{1/\beta} x$ for all
  $x\in\KK$, so that a similar argument to Theorem~\ref{thr:bound}
  implies that the stability paparemetr is at most $2/\beta$. If $\xi$
  is $\sas$ in $\R_+$ with the conventional addition and
  $\alpha\in(0,1)$, then $\eta=\xi^{1/\beta}$ is stable with parameter
  $\alpha\beta$ in the newly defined cone.  The case $\beta=\infty$
  corresponds to the maximum operation in $\R_+$, where the
  characteristic exponent $\alpha$ takes any value from $(0,\infty)$.
  The same construction with $\beta<0$ gives a cone with negative
  range of $\alpha$, \cf Example~\ref{ex:harmon}.
\end{example} 

In general, the range of possible parameters of stable laws may be
used as a characteristic of a cone that, in a sense, assesses the
extent by which the second distributivity law (\ref{eq:cone4}) is
violated.

\bigskip

The following results about possible values of $\alpha$ rely entirely
on the properties of the norm if $\KK$ is a normed cone. It does not
refer to the characters on $\KK$.

\begin{lemma}
  \label{lem:alpha-norm}
  Let $\KK$ be a normed cone such that 
  \begin{equation}
    \label{eq:sum-less}
    \|x+y\|\geq \|x\|
  \end{equation}
  for all $x,y\in\KK$. Then every $\sas$ law in $\EE$ has $\alpha>0$.
\end{lemma}
\begin{proof}
  It follows from (\ref{eq:a1alph-a+b1-}) and (\ref{eq:sum-less}) that
  $2^{1/\alpha}\|\xi\|$ is stochastically greater than $\|\xi\|$, \ie
  the cumulative distribution function of the first is not greater
  than of the second one. If $\alpha<0$, this is only possible if
  $\|\xi\|$ vanishes or is infinite, while both these cases are
  excluded by requiring that $\xi$ is $\EE$-valued.
\end{proof}

Note that (\ref{eq:sum-less}) holds in such cones like the positive
half-line with the sum or maximum operation, where the addition always
increases the norm. It is interesting to note that a cone with
$\neutral\neq\zero$ that satisfies (\ref{eq:sum-less}), \textbf{(C)}
and \textbf{(E)} does not possess any non-trivial proper $\sas$ random
element. Indeed Theorem~\ref{th:alpha} implies that $\alpha<0$. For
instance, this is the case for the cone of compact sets with the union
operation, see Example~\ref{ex:usrs}.

\begin{lemma}
  \label{lem:a01}
  If $\KK$ has a sub-invariant norm, then the range of possible values
  of the characteristic exponent $\alpha$ includes $(0,1)$.
\end{lemma}
\begin{proof}
  It suffices to refer to Theorem~\ref{th:1} which provides an
  explicit construction of $\sas$ random elements with
  $\alpha\in(0,1)$ by the LePage series.
\end{proof}

\section{Integral representations of stable laws}
\label{sec:integr-repr-stable}

\subsection{Integral representations of negative definite functions}
\label{sec:levy-function-levy}

The theory of integral representations of negative definite functions
\cite[Ch.~4]{ber:c:r} makes it possible to gain further insight into
the structure of the function $\phi$ from
(\ref{eq:echixi=exp-phichi-}), \ie the Laplace exponent of a $\sas$
random element $\xi$. The random element is always assumed to be
proper in this section, so that the corresponding Laplace exponent is
finite.  Let us first introduce several important ingredients of these
integral representations specified for the cone $\KDU$ of characters
with its dual $\KDD$ introduced in Section~\ref{sec:laplace-transform}.

A \emph{L\'evy measure} is a Radon measure on $\KDD\setminus\{\one\}$
such that
\begin{equation}
  \label{eq:lev-mes}
  \int_{\KDD\setminus\{\one\}}
  (1-\Re\rho(\chi))\,\lambda(d\rho)<\infty
\end{equation}
for all $\chi\in\KDU$. Here $\one$ is the neutral element in $\KDD$,
\ie the character identically equal to 1. A function
$\ell:\KDU\mapsto\R$ is said to be $\star$-additive if
\begin{equation}
  \label{eq:stadd}
  \ell(\chi_1\chi_2)=\ell(\chi_1)+\ell(\chi_2)
  \quad \text{and}\quad
  \ell(\overline{\chi})=-\ell(\chi)
\end{equation}
for all $\chi_1,\chi_2,\chi\in\KDU$ (recall that the involution of a
character is its complex conjugate).  The function $\chi\mapsto
e^{i\ell(\chi)}$ is thus a character on $\KDU$, \ie an element of
$\KDD$.  A function $q:\KDU\mapsto\R$ is called a \emph{quadratic
  form} if
\begin{equation}
  \label{eq:qwf}
  2q(\chi_1)+2q(\chi_2)=q(\chi_1\chi_2)+q(\chi_1\overline{\chi_2})
\end{equation}
for all $\chi_1,\chi_2\in\KDU$. A real-valued function
$L(\chi,\rho)$ defined on $\KDU\times\KDD$ is called the \emph{L\'evy
function} if $L$ is $\star$-additive with respect to $\chi$ for each
$\rho$, Borel measurable with respect to $\rho$ for each $\chi$, and
\begin{equation}
  \label{eq:int_kdds-1-rhoch}
  \int_{\KDD\setminus\{\one\}}
  |1-\rho(\chi)+iL(\chi,\rho)|\,\lambda(d\rho)<\infty
\end{equation}
for each L\'evy measure $\lambda$.  It is shown in \cite{buc86} that
a L\'evy function exists, can be chosen to be continuous with
respect to its second argument and to satisfy
$L(\chi,\overline\rho)=-L(\chi,\rho)$, see also
\cite[Th.~3.1]{berg90}. We fix a certain L\'evy function constructed
according to \cite{buc86} for the semigroup $\KDU$.

If the Laplace exponent $\phi$ from~(\ref{eq:echixi=exp-phichi-}) is
finite (\ie $\xi$ is proper), it can be represented as
\begin{equation}
  \label{eq:ldcomplex}
  \phi(\chi)=i\ell(\chi)+q(\chi)+\int_{\KDD\setminus\{\one\}}
  (1-\rho(\chi)+iL(\chi,\rho))\,\lambda(d\rho)\,,
  \quad \chi\in\KDU\,,
\end{equation}
for a unique triple $(\ell,q,\lambda)$ of a $\star$-additive function,
non-negative quadratic form and the L\'evy measure, see
\cite[Th.~4.3.19]{ber:c:r}.  Note that $\phi(\one)=0$. If $\KK$ is a
group, the quadratic form $q$ corresponds to the Gaussian component of
$\xi$, see \cite[Sec.~IV.6]{par67}. Following this terminology, we say
that $\xi$ does not have a \emph{Gaussian component} if $q$ in
(\ref{eq:ldcomplex}) vanishes. Note that the elements of the integral
representation may depend on the choice of $\KDU$. The family $\KDU$
is supposed to be fixed in the sequel.

If the involution is identical, then $\rho$ takes real values, the
L\'evy function and $\ell$ vanish, so that (\ref{eq:ldcomplex}) turns
into
\begin{equation}
  \label{eq:ldc}
  \phi(\chi)=q(\chi)
  +\int_{\KDD\setminus\{\one\}}(1-\rho(\chi))\,\lambda(d\rho)\,, 
\end{equation}
where (\ref{eq:qwf}) for $q:\KDU\mapsto\R_+$ turns into
\begin{equation}
  \label{eq:qadd}
  q(\chi_1\chi_2)=q(\chi_1)+q(\chi_2)
\end{equation}
for all $\chi_1,\chi_2\in\KDU$, see \cite[Th.~4.3.20]{ber:c:r}. In
this case $q$ is also called a quadratic form. Furthermore,
$e^{-q(\chi)}$ is a character from $\KDD$. If it is possible to
associate this character with a certain $z\in\KK$ using the evaluation
map, \ie if $e^{-q(\chi)}=\chi(z)$ for all $\chi\in\KDU$, then the
Gaussian component corresponds to the deterministic point $z$.

\subsection{L\'evy measures of $\sas$ laws}
\label{sec:levy-measures-sas}

Representations (\ref{eq:ldcomplex}) and (\ref{eq:ldc}) hold for every
infinitely divisible random element $\xi$ with finite Laplace
exponent.  The $\sas$ property of $\xi$ can be used to characterise
the elements of the triplet $(\ell,q,\lambda)$.  For this, uplift the
multiplication by numbers to $\KDD$ from $\KDU$ by setting
\begin{displaymath}
  (a\circ\rho)(\chi)=\rho(a\circ\chi)\,,\quad a>0\,.
\end{displaymath}

\begin{Th}
  \label{th:l-homog}
  Assume that Condition~\textbf{(C)} holds.  If $\phi$ is the Laplace
  exponent of a proper $\sas$ random element $\xi$, then the
  corresponding L\'evy measure is homogeneous on $\KDD$, \ie
  \begin{equation}
    \label{eq:lambda-hom}
    \lambda(s\circ B)=s^{-\alpha}\lambda(B)\,, \quad s>0\,,
  \end{equation}
  for each Borel set $B$ in $\KDD$ and $\lambda$ has infinite total
  mass.
\end{Th}
\begin{proof}
  Let $\mu_t$ be the distribution of $t^{1/\alpha}\xi$. It follows
  from Theorem~\ref{th:rep}\textsl{(iii)} that the Laplace transform
  of $\mu_t$ is given by $e^{-t\phi(\chi)}$, so that $\{\mu_t,\ t>0\}$
  is the convolution semigroup associated with $\phi$. By
  \cite[Lemma~4.3.12]{ber:c:r}, the L\'evy measure $\lambda$ is the
  vague limit as $t\downarrow0$ of the images of $t^{-1}\mu_t$ under
  the evaluation map $\imath$. Since $t^{-1}\mu_t(sA)=s^{-\alpha}
  r^{-1}\mu_r(A)$ with $t=s^\alpha r$, the corresponding vague limit
  $\lambda$ satisfies (\ref{eq:lambda-hom}).
  
  By taking $B$ with $\lambda(B)>0$ and letting $s\ti$ in
  (\ref{eq:lambda-hom}) in case $\alpha<0$ and $s\downarrow0$ in case
  $\alpha>0$, it follows that $\lambda$ has infinite total mass.
\end{proof}

Since the left-hand side of (\ref{eq:lambda-hom}) can be written as
$(s^{-1}\circ\lambda)(B)$, we say that $\lambda$ has the
\emph{homogeneity order} $\alpha$.  The following result provides an
upper bound for the possible homogeneity order of the L\'evy measure
and thereupon can complement Theorem~\ref{thr:bound} even without
using the second distributivity law. It relies instead on the local
behaviour of characters near the neutral element.

\begin{Th}
  \label{th:max-alpha}
  Assume that Condition~\textbf{(C)} holds. For some $\beta>0$, all
  $\chi\in\KDU$ and all $\rho\in\KDD\setminus\{\one\}$ define
  \begin{equation}
    \label{eq:rho-local}
    g(\chi,\rho)=\liminf_{t\downarrow0}
    \frac{1-\Re (t\circ\rho)(\chi)}{t^\beta}\,.
  \end{equation}
  Assume that for all $\rho\neq\one$ there exists $\chi\in\KDU$ such
  that $g(\chi,\rho)>0$.  Then the L\'evy measure $\lambda$ of a
  proper $\sas$ random element $\xi$ has the
  homogeneity order which is strictly less than $\beta$.
  
  If (\ref{eq:rho-local}) holds for $t\ti$ and $\beta<0$, then the
  same condition on $g$ implies that the order of homogeneity is
  strictly greater than $\beta$.
\end{Th}
\begin{proof}
  If the homogeneity order is $\alpha$, then for all sufficiently
  small $t>0$,
  \begin{align*}
    \int_{\KDD\setminus\{\one\}} (1-\Re\rho(\chi))\lambda(d\rho)
    &=\int_{\KDD\setminus\{\one\}}
    (1-\Re(t\circ\rho)(\chi))\lambda(t\circ d\rho)\\
    &=\int_{\KDD\setminus\{\one\}}
    \frac{1-\Re (t\circ\rho)(\chi)}{t^\beta}
    t^{\beta-\alpha}\lambda(d\rho) \\
    &\geq \int_{\KDD\setminus\{\one\}} g(\chi,\rho)
    t^{\beta-\alpha}\lambda(d\rho)\,.
  \end{align*}
  Since, by Theorem~\ref{th:l-homog}, $\lambda$ has infinite total
  mass, the obtained expression converges to infinity as
  $t\downarrow0$ if $\alpha\geq\beta$ contrary to (\ref{eq:lev-mes})
  unless $\lambda$ is supported by $\{\rho:\;g(\chi,\rho)=0\}$ for all
  $\chi$. The latter is, however, impossible in view of the condition on
  $g$.  The second statement is proved similarly.
\end{proof}

Define a measure $\Lambda$ on $\sB(\KK)$ as the inverse image under
the natural map of $\lambda$ restricted onto $\imath(\KK)$. The
measure $\Lambda$ is homogeneous on $\KK$ with the same order as
$\lambda$.  Although it is generally impossible to examine
(\ref{eq:rho-local}) for \emph{all} characters $\rho$ from the second
dual semigroup, it is possible to assess the homogeneity order of
$\Lambda$ by applying (\ref{eq:rho-local}) to $\rho=\rho_x$ for
$x\in\KK$. If
\begin{displaymath}
  g(\chi,x)=\liminf_{t\downarrow0}\frac{1-\Re \chi(tx)}{t^\beta}>0\,,
\end{displaymath}
for all $x\neq\neutral$ and some $\chi\in\KDU$, then the homogeneity
order of $\Lambda$ is less than $\beta$, \ie the maximal homogeneity
order corresponds to the order of decrease of $1-\Re\chi(tx)$ as
$t\downarrow0$.

\begin{Cor}
  \label{cor:max-bound}
  If \textbf{(C)} and the second distributivity law~\eqref{eq:cone4}
  hold, then the order of homogeneity of the L\'evy measure of a
  proper $\sas$ random element is strictly less than 2.  If also the
  involution is identical, then the order of homogeneity is strictly
  less than~1.
\end{Cor}
\begin{proof}
  If the second distributivity law holds, then 
  \begin{displaymath} 
    (n^{-1}\circ\rho)(\chi)=\rho(n^{-1}\circ\chi)=\rho(\chi)^{1/n}\,.
  \end{displaymath} 
  The inequality
  \begin{equation}
    \label{eq:ineq-bcr}
    \frac{1}{n^2}(1-\Re (\rho(\chi)^n))
    \leq \frac{\pi^2}{4}(1-\Re\rho(\chi))
  \end{equation}
  from \cite[p.~109]{ber:c:r} implies that 
  \begin{displaymath}
    n^2(1-\Re (n^{-1}\circ\rho)(\chi))
    \geq \frac{4}{\pi^2}(1-\Re \rho(\chi))\,.
  \end{displaymath}
  Then Theorem~\ref{th:max-alpha} is applicable with $\beta=2$ and
  $g(\chi,\rho)\geq 1-\Re \rho(\chi)$. 

  Another inequality from \cite[p.~109]{ber:c:r} reads that
  \begin{displaymath}
    \frac{1}{n}(1-|\rho(\chi)|^{2n}) \leq 1-\Re
    \rho(\chi\overline\chi)\,. 
  \end{displaymath}
  If the involution is identical, then
  \begin{displaymath}
    n(1-(n^{-1}\circ\rho)(\chi)) \geq 1-|\rho(\chi)|\,,
  \end{displaymath}
  whence the homogeneity order for the L\'evy measure is strictly less
  than~1. 
\end{proof}

By combining Corollary~\ref{cor:max-bound} with
Theorem~\ref{thr:bound}\textsl{(i)} and the continuity property of
$\phi$, we arrive at the following result.

\begin{Cor}
  \label{cor:order2}
  If \textbf{(C)} and the second distributivity law~\eqref{eq:cone4}
  hold, then every proper $\sas$ random element with $\alpha=2$ has
  the Laplace exponent given by
  \begin{displaymath}
    \phi(\chi)=i\ell(\chi)+q(\chi)\,, \quad \chi\in\KDU\,,
  \end{displaymath}
  where $\ell$ is a $\star$-additive continuous function and $q$ is a
  continuous non-negative quadratic form on $\KDU$.
  
  If the involution is identical, then every proper $\sas$ random
  element with $\alpha=1$ has the Laplace exponent
  $\phi(\chi)=q(\chi)$ being continuous non-negative linear functional
  of $\chi$.
\end{Cor}

Without the second distributivity law, the value $\alpha=2$ does not
necessarily correspond to the Gaussian component, \eg in
$(\R_+,\vee)$.

\subsection{Quadratic form and L\'evy function}
\label{sec:quadratic-form-star}

It remains to explore implications of the stability of $\xi$ on the further
ingredients of (\ref{eq:ldcomplex}), namely the quadratic form $q$,
the L\'evy function and the function $\ell$. The following result
describes some cases when $\xi$ does not possess a Gaussian component.

\begin{Th}
  \label{thr:l-and-q}
  Assume that Condition~\textbf{(C)} and the second distributivity
  law~\eqref{eq:cone4} hold. Let $q$ be the quadratic form in the
  integral representation~(\ref{eq:ldcomplex}) of a proper $\sas$
  random element.
  \begin{itemize}
  \item[\textsl{(i)}] $q$ vanishes unless $\alpha=1$ or $\alpha=2$. If
    $\KK$ is a group, then $q$ vanishes if $\alpha\neq 2$. If
    $\alpha=1$, then $q(\overline\chi)=q(\chi)$ for all $\chi$ and $q$
    satisfies (\ref{eq:qadd}).
  \item[\textsl{(ii)}] If the involution is identical, then $q$ from
    representation~(\ref{eq:ldc}) vanishes unless $\alpha=1$.
  \end{itemize}
\end{Th}
\begin{proof}
  \textsl{(i)} The $\star$-additivity of $L(\cdot,\rho)$ implies that
  $L(\chi^2,\rho)=2L(\chi,\rho)$. By the second distributivity law,
  $2\circ\chi=\chi^2$, so that
  $\rho(\chi^2)=\rho^2(\chi)=(2\circ\rho)(\chi)$.
  
  By Theorem~\ref{th:rep},
  $\phi(\chi^2)=\phi(2\circ\chi)=2^\alpha\phi(\chi)$, hence it is also
  finite. Furthermore,
  \begin{align}
    \label{eq:rescal}
    \phi(\chi^2)&=i\ell(\chi^2)+q(\chi^2)+  \int_{\KDD\setminus\{\one\}} 
    (1-\rho(\chi^2)+iL(\chi^2,\rho))\,\lambda(d\rho) \notag\\
    &= 2i\ell(\chi)+q(\chi^2)+\int_{\KDD\setminus\{\one\}}
    (1-(2\circ\rho)(\chi)+2L(\chi,\rho))\lambda(d\rho)\,,
  \end{align}
  while, by  Theorem~\ref{th:l-homog},
  \begin{align}
    \label{eq:rescal2}
    2^\alpha\phi(\chi)&=2^\alpha i\ell(\chi)+2^\alpha q(\chi)+
    \int_{\KDD\setminus\{\one\}}
    (1-\rho'(\chi)+iL(\chi,\rho'))\,\lambda(2^{-1}\circ d\rho')\notag \\
    &= 2^\alpha i\ell(\chi)+2^\alpha
    q(\chi)+\int_{\KDD\setminus\{\one\}}
    (1-(2\circ\rho)(\chi)+iL(\chi,2\circ \rho))\lambda(d\rho).
  \end{align}
  Equating the right-hand sides of \eqref{eq:rescal} and
  \eqref{eq:rescal2}, we arrive at
  \begin{equation}
    \label{eq:ReIm}
    2^\alpha i\ell(\chi)+2^\alpha q(\chi)
    =2i(\ell(\chi)-\tilde\ell(\chi))+q(\chi^2)\,,
  \end{equation}
  where the integral
  \begin{displaymath}
    \tilde\ell(\chi)
    =\frac{1}{2}\int_{\KDD\setminus\{\one\}}
    (L(\chi,2\circ\rho)-2L(\chi,\rho))\lambda(d\rho)
  \end{displaymath}
  is finite and may be regarded as a real-valued $\star$-additive
  function on $\KDU$.
  
  By comparing the coefficients of the real parts in~\eqref{eq:ReIm},
  we see that
  \begin{equation}
    \label{eq:aq}
    q(\chi^2)=2^\alpha q(\chi)\,.
  \end{equation}
  Next, (\ref{eq:qwf}) applied to a real-valued character $\chi$,
  yields that $2q(\chi)=q(\chi^2)$.  Since $\alpha\neq1$, we have
  $q(\chi)=0$ for every real-valued character $\chi$. Since
  $\chi\overline\chi$ is a real-valued character,
  $q(\chi\overline\chi)=0$. Now (\ref{eq:qwf}) applied to
  $\chi_1=\chi_2=\chi$ for a (not necessarily real-valued character)
  $\chi$ yields that $4q(\chi)=q(\chi^2)=2^\alpha q(\chi)$.  Since
  $\alpha\neq2$, $q$ vanishes identically.
  
  If $\KK$ is a group, then putting $\chi_1=\chi_2=\chi$ in
  (\ref{eq:qwf}) and noticing that $\overline{\chi}=\chi^{-1}$, we
  obtain that $q$ vanishes unless $\alpha=2$.

  Assume that $\alpha=1$. By applying (\ref{eq:qwf}) to
  $\chi_1=\chi_2=\chi$ and using~\eqref{eq:aq} we obtain that
  $q(\chi\overline\chi)=2q(\chi)$. Then 
  \begin{align*}
    2q(\chi_1)+2q(\chi_2)
    &=\frac{1}{2}\bigl(2q(\chi_1\overline{\chi_1})
    +2q(\chi_2\overline{\chi_2})\bigr)\\ 
    &=\frac{1}{2}\bigl(q(\chi_1\overline{\chi_1}\chi_2\overline{\chi_2})
    +q(\chi_1\overline{\chi_1}\overline{\chi_2\overline{\chi_2}})\bigr)\\
    &=q(\chi_1\overline{\chi_1}\chi_2\overline{\chi_2})\\
    &=2q(\chi_1\chi_2)\,.
  \end{align*}
  Combining this with (\ref{eq:qwf}) yields that
  $q(\overline\chi)=q(\chi)$.
  
  \textsl{(ii)} By Theorem~\ref{thr:bound}\textsl{(ii)}, $\alpha\leq
  1$. If $\alpha<1$, then by \cite[Th.~4.3.20]{ber:c:r},
  \begin{displaymath}
    q(\chi)=\lim_{n\ti} \frac{\phi(\chi^n)}{n}
    =\lim_{n\ti} \frac{\phi(n\circ\chi)}{n}
    =\lim_{n\ti} \frac{n^\alpha\phi(\chi)}{n}=0\,,
  \end{displaymath}
  noticing that $\phi(\chi)$ is finite.
\end{proof}

The following result concerns idempotent semigroups, where the second
distributivity law never holds.

\begin{lemma}
  \label{lem:q-idempo}
  If $\KK$ is an idempotent semigroup, then $q$ in (\ref{eq:ldc})
  vanishes identically.
\end{lemma}
\begin{proof}
  The idempotency of $\KK$ implies that all characters take values $0$
  or $1$. Thus, $\KDU$ is also an idempotent semigroup. By
  (\ref{eq:qadd}), $q(\chi)=q(\chi^2)=2q(\chi)$, whence $q$
  vanishes.
\end{proof}

In a number of cases the L\'evy function vanishes or may be set to
zero. Particular important instances of this are mentioned in the
following theorem.

\begin{Th}
  \label{thr:no-ell}
  Assume that at least one of the following conditions holds.
  \begin{itemize}
  \item[\textsl{(i)}] The involution is identical.
  \item[\textsl{(ii)}] For all $\chi\in\KDU$ 
    \begin{equation}
      \label{eq:fininteg}
      \int_{\KDD\setminus\{\one\}}(1-\rho(\chi))\,\lambda(d\rho)<\infty\,.
    \end{equation}
  \end{itemize}
  Then the Laplace exponent of any proper $\sas$ random element $\xi$ is
  given by
  \begin{equation}
    \label{eq:intrep-nolevy}
    \phi(\chi)=i\ell(\chi)+q(\chi)+\int_{\KDD\setminus\{\one\}}
    (1-\rho(\chi))\,\lambda(d\rho)\,,
    \quad \chi\in\KDU\,,
  \end{equation}
  for a $\star$-additive function $\ell$ and a quadratic form $q$. If
  \textbf{(C)} and the second distributivity law~\eqref{eq:cone4}
  hold, then $\ell$ vanishes unless $\alpha=1$.
\end{Th}
\begin{proof}
  If the involution is identical, then the L\'evy function vanishes.
  If (\ref{eq:fininteg}) holds, then
  \begin{displaymath}
    \int_{\KDD\setminus\{\one\}}L(\chi,\rho)\,\lambda(d\rho)
  \end{displaymath}
  is a finite $\star$-additive functional of $\chi$, so that it can be
  combined with $\ell(\chi)$ from (\ref{eq:ldcomplex}).
  
  By Theorem~\ref{th:l-homog}, the integral in
  (\ref{eq:intrep-nolevy}) is homogeneous of order $\alpha$. By the
  homogeneity and~\eqref{eq:cone4},
  $\phi(\chi^2)=\phi(2\circ\chi)=2^\alpha\phi(\chi)$, so that
  $2\ell(\chi)=2^\alpha\ell(\chi)$. Thus, $\ell$ vanishes unless
  $\alpha=1$.
\end{proof}

\subsection{Symmetric random elements}
\label{sec:symm-rand-elem}

\begin{definition}
  \label{def:symmetric}
  A $\KK$-valued random element is called \emph{symmetric} if $\xi$
  coincides in distribution with its involution $\xi^\star$.
\end{definition}

In case of the identical involution all random elements can be
regarded as being symmetric. More generally, a symmetric element $\xi$
can be obtained as $\xi_1+\xi^\star_2$ for \iid random elements
$\xi_1$ and $\xi_2$.  This construction is a generalisation of the
\emph{symmetrisation} procedure for random elements in Banach spaces.
By applying the involution to the both sides of
(\ref{eq:a1alph-a+b1-}) it is clear that if $\xi$ is $\sas$ then
$\xi^\star$ also is.

The \emph{principal value} of an integral over $\KDD\setminus\{\one\}$
is defined as the limit of the integrals over a sequence $\{F_n\}$ of
symmetric sets as $F_n\uparrow\KDD\setminus\{\one\}$.  The symmetry of
$F\subset\KDD$ is understood with respect to the involution on $\KDD$,
\ie with respect to the complex conjugation.

\begin{Cor}
  \label{cor:symm}
  If \textbf{(C)} holds, then the Laplace exponent of every symmetric
  proper $\sas$ random element is given by
  \begin{equation}
    \label{eq:l-exp-sym}
    \phi(\chi)=q(\chi)+\int_{\KDD\setminus\{\one\}}
    (1-\rho(\chi))\,\lambda(d\rho)\,,
    \quad \chi\in\KDU\,,
  \end{equation}
  where $\lambda$ is a symmetric homogeneous measure on $\KDD$, $q$ is
  a non-negative quadratic form, and the principal value of the
  integral converges.
\end{Cor}
\begin{proof}
  By repeating the argument from the proof of
  Theorem~\ref{th:l-homog} it is easy to see that $\lambda$ is
  symmetric with respect to the complex conjugate operation being the
  involution on~$\KDD$.
  
  The integral in (\ref{eq:ldcomplex}) converges and so does its
  principal value as the limit of the integrals over $F_n$ as
  $F_n\uparrow\KDD\setminus\{\one\}$.  The symmetry property of $F_n$
  and $\lambda$ implies that
  \begin{displaymath}
    \int_{F_n} L(\chi,\rho)\lambda(d\rho)=
    \int_{F_n} L(\chi,\overline\rho)\lambda(d\overline\rho)=
    - \int_{F_n} L(\chi,\rho)\lambda(d\rho)\,,
  \end{displaymath}
  whence
  \begin{displaymath}
    \int_{\KDD\setminus\{\one\}} L(\chi,\rho)\lambda(d\rho)=0
  \end{displaymath}
  as the principal value. Note that
  $\E\chi(\xi^\star)=\E\overline\chi(\xi)$. Since the integral of the
  L\'evy function vanishes and $\phi(\chi)=\phi(\overline\chi)$, we
  obtain that $\ell(\chi)=0$ in (\ref{eq:intrep-nolevy}).
\end{proof}

\bigskip

Consider now a Gaussian random element $\xi$ in $\KK$, \ie assume that
its Laplace transform is given by $\E\chi(\xi)=e^{-q(\chi)}$ for
$\chi\in\KDU$. Note that in the absence of the second distributivity
law a Gaussian random element is not necessarily stable.

\begin{Th}
  \label{thr:gauss-1}
  Assume that \textbf{(C)} holds. If $\xi$ is a Gaussian element, then
  $\xi$ is symmetric, $\chi(\xi)$ is deterministic for every
  real-valued character $\chi$ and $\xi+\xi^\star$ is deterministic.
  In particular, $\xi$ is deterministic if the involution is
  identical.
  
  Every character $\chi$ with values in the unit complex circle
  $\TT$ is representable as $\chi(x)=e^{i u(x)}$ with an additive
  real-valued function $u$ such that $u(\xi)$ has a normal
  distribution.
\end{Th}
\begin{proof}
  By applying (\ref{eq:qwf}) to $\chi_1=\one$, we see that
  $q(\overline{\chi_2})=q(\chi_2)$, \ie $q$ is involution symmetric.
  This immediately leads to the conclusion that
  $\E\chi(\xi)=\E\overline{\chi}(\xi)=\E\chi(\xi^\star)$, so that
  Theorem~\ref{thr:laplace} yields that $\xi\deq\xi^\star$.
  
  It is easy to see that if at least one of $\chi_1,\chi_2$ is
  real-valued in~(\ref{eq:qwf}), then
  $q(\chi_1\chi_2)=q(\chi_1)+q(\chi_2)$. Therefore, if $\chi$ is a
  real-valued character, then
  \begin{displaymath}
    \E\chi^2(\xi) = e^{-q(\chi^2)}=e^{-2q(\chi)}
    =(\E\chi(\xi))^2\,.
  \end{displaymath}
  Thus $\var \chi(\xi)=0$, hence $\chi(\xi)$ is deterministic for
  every real-valued $\chi$.  If the involution is identical, the
  strict separation condition implies that $\xi$ itself is
  deterministic.
  
  Decompose now a general $\chi\in\KDU$ as the product $\chi'\chi''$ of a
  real-valued character and a $\TT$-valued character. Since
  $\chi''\overline{\chi''}=\one$,
  \begin{displaymath}
    \E \chi(\xi+\xi^\star)=\E(\chi'(\xi))^2\,.
  \end{displaymath}
  But $(\chi')^2$ is a real-valued character, hence, by the above, the last
  expression equals
  \begin{displaymath}
    (\chi'(\xi))^2=\chi'(\xi)\chi'(\xi^\star)\chi''(\xi)\chi''(\xi^\star)
    =\chi(\xi+\xi^\star)\,.
  \end{displaymath}
  Thus, $\chi(\xi+\xi^\star)$ is deterministic, so that the conclusion
  of the theorem follows from the separation condition. 

  Let $\chi$ be a $\TT$-valued character. By iterating (\ref{eq:qwf}),
  it is seen that $q(\chi^n)=n^2 q(\chi)$ for each $n$. Interpreting
  this identity in terms of the Laplace transform of $\xi$ yields that
  \begin{displaymath}
    \E\chi^s(\xi)=(\E\chi(\xi))^{s^2}
  \end{displaymath}
  for every positive rational number $s$. Using the representation of
  $\chi$, we have 
  \begin{displaymath}
    \E e^{isu(\xi)}=(\E\chi(\xi))^{s^2}\,.
  \end{displaymath}
  The symmetry property of $\xi$ implies that $\E\chi(\xi)$ is a real
  number which does not exceed 1, so that $\E\chi(\xi)=e^{-a^2/2}$ for
  some $a$. Thus, $\E e^{isu(\xi)}=e^{-a^2s^2/2}$ meaning that $u(\xi)$
  is normally distributed with mean zero and variance $a^2$.
\end{proof}

If the Gaussian element $\xi$ from Theorem~\ref{thr:gauss-1} is
$\sas$, then deterministic $\xi+\xi^\star$ is also $\sas$, hence an
element of $\KK(\alpha)$. Also, $\xi\in\KK(\alpha)$ if the involution
is identical. Note that a Gaussian random element may be $\sas$ with
arbitrary $\alpha$, see Example~\ref{ex:new-group}.

\section{LePage series representation of $\sas$ random elements}
\label{sec:lepage-seri-repr}

\subsection{LePage series on the second dual semigroup}
\label{sec:lepage-series-second}

In Section~\ref{sec:lepage-series} we have shown that the LePage
series~(\ref{eq:9}) (or the integral $\int x\Pi_\alpha(dx)$) defines a
$\sas$ random element. It is natural to ask if any $\sas$ random
element admits such a representation as it is the case for stable
distributions in Banach spaces, see \cite[Cor.~4.10]{rosi90}. We
address this question by using the integral representations
(\ref{eq:ldcomplex}) and (\ref{eq:ldc}) for the Laplace exponent
$\phi$ of a $\sas$ random element. This idea is supported by the
formula (\ref{eq:pgfp}) for the probability generating functional of a
Poisson process, which is quite similar to (\ref{eq:ldc}). The
intensity measure of this Poisson process is the L\'evy measure
$\lambda$, so that this process lives on $\KDD$ and therefore is
denoted by $\Pi^\ddual$.

We first characterise the weak convergence of $\KDD$-valued random
elements.

\begin{lemma}
  \label{lem:kdd-weak}
  A sequence $\xi_n^\ddual$ of random elements in $\KDD$ with its
  Borel $\sigma$-algebra $\sB(\KDD)$ weakly converges to a random
  element $\xi^\ddual$ (notation $\xi_n^\ddual\wto \xi^\ddual$) if and
  only if $\E\xi_n^\ddual(\chi)\to\E\xi^\ddual(\chi)$ for all $\chi$
  from a separating semigroup $\KDU$.
\end{lemma}
\begin{proof}
  Note that $\KDU$ is a separating family of characters on $\KDD$
  acting as $\chi(\rho)=\rho(\chi)$. Furthermore,
  $\sF(\KDD;\KDU)=\sB(\KDD)$, so that Theorem~\ref{thr:laplace}
  is applicable. Each character $\chi:\KDD\mapsto\DD$ is continuous
  on $\KDD$, since $\chi(\rho_n)=\rho_n(\chi)\to\rho(\chi)=\chi(\rho)$
  if $\rho_n$  converges to $\rho$ pointwise.
  
  The result now follows from the convergence of Laplace transforms $\E
  \xi_n^\ddual(\chi)$ together with the tightness condition that is
  clearly fulfilled because of the compactness of $\KDD$, see
  \cite[Th.~IV.3.1]{vak:tar:cho87}.
\end{proof}

We also write $\xi_n^\ddual\wto\xi$ if $\xi_n^\ddual$ weakly converges
to the $\KDD$-valued random element $\imath(\xi)$ being the evaluation
image of $\xi\in\KK$. For this, the condition $\sF(\KK;\KDU)=\sB(\KK)$
of Theorem~\ref{thr:laplace} should hold in order to be able to treat
$\imath(\xi)$ as a $\KDD$-valued random element.

Denote by $\FF_\lambda(\KDD)$ the family of Borel sets
$F\subset\KDD\setminus\{\one\}$ such that
\begin{equation}
  \label{eq:fcomp}
  \int_{F}(1-\rho(\chi))\lambda(d\rho)<\infty
\end{equation}
for all $\chi\in\KDU$. With every $F$, the family
$\FF_\lambda(\KDD)$ contains all its measurable subsets. If the
involution is identical, then
$\KDD\setminus\{\one\}\in\FF_\lambda(\KDD)$. 

If $F\in\FF_\lambda(\KDD)$, then the product
\begin{displaymath}
  \xi_{F}^\ddual = \prod_{\rho\in F\cap\supp\Pi^\ddual} \rho
\end{displaymath}
exists in $\KDD$, \ie it weakly converges in the topology of pointwise
convergence on $\KDD$ meaning that
\begin{displaymath}
  \xi_{F\cap K_n}^\ddual(\chi)
  =\prod_{\rho\in F\cap K_n\cap\supp \Pi^\ddual} \rho(\chi)
\end{displaymath}
weakly converges to $\xi_{F}^\ddual(\chi)$ for every $\chi\in\KDU$ as
$K_n\uparrow \KDD\setminus\{\one\}$ with $\lambda(K_n)<\infty$. This
is easily seen by observing that the probability generating functional
of $\Pi^\ddual$ restricted onto $F\cap K_n$ coincides with the
expected value of $\xi_{F\cap K_n}^\ddual$ and using
Lemma~\ref{lem:kdd-weak}. Note that $\xi_{F\cap K_n}^\ddual$ is
well defined because $\Pi^\ddual$ has only a finite number of points
in $K_n$.

For any $F\in\FF_\lambda(\KDD)$ define
\begin{equation}
  \label{eq:gammar}
  \gamma_F(\chi)=\exp\left\{-i\left(\ell(\chi)
      +\int_{F}L(\chi,\rho)\lambda(d\rho)\right)\right\}\,,\quad 
  \chi\in\KDU\,.
\end{equation}
Note that $\gamma_F\in\KDD$.

\begin{Th}
  \label{thr:sec-no-pont}
  Assume that Condition~\textbf{(C)} holds.  Let $\xi$ be a proper
  $\KK$-valued $\sas$ random element without Gaussian component.
  Then there exists a unique measure $\lambda$ on $\KDD$ (which is
  then a homogeneous of order $\alpha$ L\'evy measure) such that, for
  the Poisson process $\Pi^\ddual$ with intensity measure
  $\lambda$, one has
  \begin{displaymath}
    \gamma_{F_n}\xi_{F_n}^\ddual\wto\xi
  \end{displaymath}
  for any sequence $F_n\uparrow\KDD\setminus\{\one\}$ such that
  $F_n\in\FF_\lambda(\KDD)$ for all $n$.
\end{Th}
\begin{proof}
  Let $\lambda$ be the L\'evy measure of $\xi$ that stems from
  (\ref{eq:ldcomplex}).  Formula (\ref{eq:pgfp}) for the probability
  generating functional applied to $\Pi^\ddual$ and the definition of
  $\gamma_F$ imply that
  \begin{align*}
    \gamma_{F_n}(\chi)\E\chi(\xi_{F_n}^\ddual)
    &=\gamma_{F_n}(\chi)\E\left[\prod_{\rho\in F_n\cap\supp\Pi^\ddual}
      \rho(\chi)\right]\\
    &=\exp\left\{-\left(i\ell(\chi)
        +\int_{F_n}(1-\rho(\chi)+iL(\chi,\rho))
        \,\lambda(d\rho)\right)\right\}\,.
  \end{align*}
  If $F_n\uparrow\KDD\setminus\{\one\}$, the right-hand side converges
  to the Laplace transform of a stable random element $\xi$ without
  Gaussian component, see (\ref{eq:ldcomplex}).  Finally,
  Lemma~\ref{lem:kdd-weak} shows that $\gamma_{F_n}\xi_{F_n}^\ddual$
  weakly converges to $\xi^\ddual=\imath(\xi)$. 
\end{proof}

In some cases (\seg Theorem~\ref{thr:no-ell}) the L\'evy function
vanishes or its contribution to the Laplace exponent may be subsumed
in $\ell$, so that $\gamma_{F_n}(\chi)=e^{-i\ell(\chi)}$.  If $\ell$
vanishes, then the normalisation is not needed and
\begin{equation}
  \label{eq:convchi}
  \xi_{F_n}^\ddual\wto\xi\,.
\end{equation}
An important instance of this concerns symmetric random elements.

\begin{Cor}
  \label{cor:no-pont-sym}
  If \textbf{(C)} holds and $\xi$ is a proper symmetric $\sas$ random
  element without Gaussian component in a normed cone $\KK$, then
  there exists a unique measure $\lambda$ on $\KDD$ (which is then
  necessarily symmetric and homogeneous of order $\alpha$) such that
  \begin{equation}
    \label{eq:gamm-co-no-norm}
    \E\left(\prod_{\rho\in \supp\Pi^\ddual} \rho\right)(\chi)
    =\E\chi(\xi)\,,\quad \chi\in\KDU\,,
  \end{equation}
  where $\Pi^\ddual$ is the Poisson process on $\KDD$ with intensity
  measure $\lambda$. The left-hand side of (\ref{eq:gamm-co-no-norm})
  is defined as the principal value, \ie as the limit of
  $\E\xi_{F_n}^\ddual(\chi)$ as $F_n\uparrow\KDD\setminus\{\one\}$,
  where $\{F_n,\, n\geq1\}$ are involution symmetric sets.
\end{Cor}

\subsection{Support of the L\'evy measure}
\label{sec:support-levy-measure}

Section~\ref{sec:lepage-series-second} establishes that, under rather
weak condition \textbf{(C)}, each proper $\sas$ random element admits
the LePage representation on the second dual semigroup $\KDD$. The
crucial further issue is to identify the elements of the second dual
semigroup $\KDD$ with elements of the original semigroup $\KK$.
However, this is not feasible in general, since (\ref{eq:ldcomplex})
involves integration over \emph{all} characters on $\KDU$, while it is
seldom possible to describe all characters, even on the real line with
the conventional addition.

If $\KK$ is a locally compact group, then the integral in
(\ref{eq:ldcomplex}) can be taken over all continuous characters on
$\KK$, see \cite{par67}. The celebrated Pontryagin reflexivity
property of locally compact groups, see \cite[Th.~V.24.8]{hew:ros63}
establishes that if $\KDD$ is the family of all continuous characters
on the family of continuous characters $\KDU$ on a locally compact
group $\KK$, then $\KK$ and $\KDD$ are isomorphic.  This means that
the L\'evy measures of random elements with values in locally compact
groups maybe thought of as being supported by $\KK$. However, the
duality theory for semigroups is much poorer, and the results are
available only in some special cases.

As noted in the proof of Theorem~\ref{th:l-homog}, the L\'evy measure
$\lambda$ is the vague limit of the images of $\nu_t=t\mu_{t^{-1}}$
under the evaluation map $\imath:\KK\mapsto\KDD$ as $t\ti$.  Although
$\mu_t$ (the distribution of $t^{1/\alpha}\xi$) is supported by $\KK$,
the vague limit of the $\imath$-images of $\nu_t$ may be supported by
the whole $\KDD$.  The following result establishes a condition under
which the L\'evy measure $\lambda$ is supported by $\imath(\KK)$. In
this case we say shortly that the L\'evy measure is supported by $\KK$
and write $\Lambda$ for the measure on $\KK$ being the inverse image
of $\lambda$ under the evaluation map.

We begin with a result on finiteness of the L\'evy measure and
related distributional properties of a $\sas$ random elements.

\begin{lemma}
  \label{lem:clf}
  Assume that Condition~\textbf{(C)} is satisfied. Fix any set
  $F\subset\KK$ such that 
  \begin{equation}
    \label{eq:one-f}
    \one\notin\cl(\imath(F))\,,
  \end{equation}
  where the closure is taken in the topology of pointwise convergence
  in $\KDD$.  If $\xi$ is a proper $\EE$-valued $\sas$ random element,
  then its L\'evy measure is finite on $\cl(\imath(F))$ and there
  exists a constant $a>0$ such that
  \begin{equation}
    \label{eq:f-hom}
    \limsup_{t\ti} t\P\{\xi\in rt^{1/\alpha}F\} =ar^{-\alpha}\,,
    \qquad r>0\,.
  \end{equation}
\end{lemma}
\begin{proof}
  The expression under the limit in \eqref{eq:f-hom} is $\nu_t(rF)$,
  where $\nu_t$ is defined above. By simple change of variable
  argument, the limit itself is a homogeneous in $r$ function, say,
  $f(r)$, so that $f(r)=r^{-\alpha}f(1)$. The statement will be
  proved, if we show that $f(1)=a$ is finite.
  
  The set $\cl(\imath(F))$ is closed and does not contain $\one$, so
  that it is compact in $\KDD\setminus\{\one\}$. Since the L\'evy
  measure $\lambda$ is a Radon measure \cite[Lemma~4.3.12]{ber:c:r},
  it is finite on compact sets, whence
  $\lambda(\cl(\imath(F)))<\infty$. According to the same reference,
  $\lambda$ is the vague limit of the images of $\nu_t$ under the
  natural map $\imath$, thus
  \begin{displaymath}
    f(1)=\limsup_{t\ti}\, \nu_t(F)
    \leq \lambda(\cl(\imath(F)))<\infty\,. \tag*{\qed}
  \end{displaymath}
  \renewcommand{\qed}{}
\end{proof}

By Lemma~\ref{lemma:proper}, properness of $\xi$ required in
Lemma~\ref{lem:clf} can be guaranteed by imposing
Condition~\textbf{(E)}. Condition (\ref{eq:one-f}) holds, in
particular, if $\imath(F)$ is closed and $\neutral\notin F$.  The
following lemma describes an important case, when $\imath(F)$ is
closed in $\KDD$.

\begin{lemma}
  \label{lem:fi-closed}
  Condition (\ref{eq:one-f}) is equivalent to the fact that $\neutral$
  does not belong to the $\KDU$-weak closure of $F$. In particular, if
  $F$ is $\KDU$-weak (sequentially) compact, then $\imath(F)$ is
  closed in $\KDD$ and (\ref{eq:one-f}) holds provided $\neutral\notin
  F$.
\end{lemma}
\begin{proof}
  The first statement is evident, since $\rho_{x_n}(\chi)\to 1$ for
  all $\chi\in\KDU$ is equivalent to $x_n\stackrel{w}\to\neutral$.
  Assume that $\rho_{x_n}(\chi)\to\rho(\chi)$ for all $\chi\in\KDU$,
  where $\{x_n\}\subset F$. By the compactness condition,
  $x_{n_k}\stackrel{w}\to x$ for a certain $x\in F$. Thus,
  $\chi(x_{n_k})\to\chi(x)$, so that $\rho(\chi)=\chi(x)$ for all
  $\chi$.
\end{proof}

In the sequel we make use of the following condition:
  \begin{description}
  \item[(S)] The neutral element $\neutral$ does not belong to the
    $\KDU$-weak closure of $\setr[1]$.
  \end{description}
In other words, Condition~\textbf{(S)} means that for any sequence
$\{x_n\}$ from $\setr[1]$ there exists a character $\chi$ such that
$\chi(x_n)$ does \emph{not} converge to $1=\chi(\neutral)$.

\begin{Th}
  \label{thr:lc-k}
  Assume that $\KK$ is a normed cone with compact unit sphere $\SS$
  such that Conditions \textbf{(C)}, \textbf{(E)} and \textbf{(S)}
  hold.  Then the L\'evy measure of any proper $\sas$ random element
  $\xi$ is supported by $\KK$.
\end{Th}
\begin{proof}
  Fix $\eps>0$ and $r>0$. Let $\nu_t\restr$ denote the restriction of
  $\nu_t$ onto $\setr$. Note that $t\P\{\xi\in
  rt^{1/\alpha}\setr[1]\}=\nu_t(r\setr[1])$, where $r\setr[1]=\setr$
  if $\alpha>0$ and $r\setr[1]=\setr[r^{-1}]$ if $\alpha<0$.  In view
  of Condition~\textbf{(S)} and Lemma~\ref{lem:fi-closed},
  Lemma~\ref{lem:clf} implies that the masses of all $\nu_t\restr$,
  $t\geq 1$, are totally bounded and there exists $R>0$ such that
  $\nu_t(R\setr[1])\leq \eps$ for all $t$ large enough.  Thus,
  $\{\nu_t\restr,\ t\geq 1\}$ is a tight family of measures.  Using
  the diagonal procedure as in \cite[Prop.~5.3.9]{linde} we come to
  the conclusion that there exists a measure $\Lambda$ on $\EE$ such
  that $\nu_{t_n}\restr$ weakly converges to $\Lambda\restr$ as
  $t_n\ti$ for any $r>0$. Thus, the L\'evy measure is supported by
  $\KK$, \ie $\lambda=\Lambda\circ\imath^{-1}$.
\end{proof}

In some other cases where the unit sphere is not compact in $\KK$, the
following result is useful. 

\begin{Th}
  \label{thr:restrsup}
  Assume that \textbf{(C)} holds and $\imath(\KK)$ is closed in the
  topology of pointwise convergence in $\KDD$. Then the L\'evy measure
  $\lambda$ of a proper $\sas$ random element in $\KK$ is supported by
  $\imath(\KK)$ and is the $\imath$-image of a homogeneous Borel
  measure $\Lambda$ on $\KK$.
\end{Th}
\begin{proof}
  By the Tikhonov theorem, $\KDD$ equipped with the pointwise
  convergence is a compact space.  By the condition, $\imath(\KK)$ is
  closed hence compact and the vague convergence definition
  immediately implies that the limiting measure $\lambda$ in
  Theorem~\ref{th:l-homog} is supported by $\imath(\KK)$.
\end{proof}

\begin{Cor}
  \label{cor:weak-comp}
  Assume that $\KK$ is a normed cone such that $\setr[1]^\comp$ is
  $\KDU$-weak relatively compact and the following condition holds:
  \begin{description}
  \item[(W)] If a sequence $\{x_n,\, n\geq1\}$ in $\KK$ is such that
    $\chi(x_n)$ converges for all $\chi\in\KDU$, then
    $\sup\|x_n\|<\infty$ in case $\neutral=\zero$ or
    $\liminf\|x_n\|>0$ in case $\neutral\neq\zero$.
  \end{description}
  Then $\imath(\KK)$ is closed in $\KDD$ and under Condition
  \textbf{(C)} the L\'evy measure of any proper $\sas$ random element
  is supported by $\KK$.
\end{Cor}
\begin{proof}
  By rescaling, $\setr^\comp$ is $\KDU$-weak relatively compact for
  all $r>0$.  If $\chi(x_n)\to \rho(\chi)$ for all $\chi\in\KDU$,
  Condition~\textbf{(W)} yields that $x_n\in \setr^\comp$ for some
  $r>0$ and all $n\geq1$. The weak relative compactness condition
  implies that $\rho(\chi)=\chi(x)$ for some $x\in\KK$ as in the proof
  of Lemma~\ref{lem:fi-closed}, so that $\imath(\KK)$ is closed.
\end{proof}

The above result is applicable if $\KK$ is a reflexive Banach space,
see Example~\ref{ex:b-space}. Although a general Banach space is also
Pontryagin reflexive if the family of continuous functionals is
equipped with the compact-open topology \cite{smit52}, we were not
able to make use of this fact to show that the L\'evy measure is always
supported by $\KK$.

\begin{remark}
  \label{rem:both-top}
  One particular simple instance, when both Theorem~\ref{thr:lc-k} and
  Corollary~\ref{cor:weak-comp} apply is when the $\KDU$-weak convergence
  is equivalent to the metric convergence in $\KK$ and the unit sphere
  is compact. Compactness of the sphere is an important requirement here
  as Example~\ref{ex:non-complete} shows.
\end{remark}

\subsection{LePage series constructed from the L\'evy measure}
\label{sec:lepage-seri-constr}

If the L\'evy measure $\Lambda$ of a proper $\sas$ random element
$\xi$ is supported by $\KK$, then the corresponding Laplace
exponent~\eqref{eq:ldcomplex} takes the form
\begin{equation}
  \label{eq:phich-1-chix+}
  \phi(\chi)=i\ell(\chi)+q(\chi)+\int_{\KK\setminus\{\neutral\}}
  (1-\chi(x)+iL(\chi,x))\,\Lambda(d x)\,, 
  \quad \chi\in\KDU\,.
\end{equation}
We slightly abuse the notation here by writing $L(\chi,x)$ instead of
$L(\chi,\rho_x)$ for $\rho_x=\imath(x)$.  Note that the L\'evy measure
$\Lambda$ satisfies
\begin{equation}
  \label{eq:Lambda-x}
  \int_{\KK\setminus\{\neutral\}}(1-\Re\chi(x))\,\Lambda(d x)<\infty
\end{equation}
for all $\chi\in\KDU$. 

In order to be able to consider the Poisson process with intensity
measure $\Lambda$, we have to ensure that $\Lambda(\setr)<\infty$, \ie
$\Lambda$ is finite on the sets in $\EE$ which are separated from
$\neutral$ by a positive distance. Condition \textbf{(C)} on its own
makes it possible to derive from (\ref{eq:Lambda-x}) that $\Lambda$ is
locally finite on $\KK\setminus\{\neutral\}$. Indeed, if
$x\neq\neutral$, then the neighbourhoods of $x$ and $\neutral$ can be
strictly separated by some character $\chi$.  Therefore,
$\Re\chi(y)\leq 1-\eps$ for some $\eps>0$ and all $y$ from a
neighbourhood of $x$. It follows from (\ref{eq:Lambda-x}) that
$\Lambda$-measure of this neighbourhood is finite.

However, the local finiteness of $\Lambda$ alone does not imply the
finiteness of $\Lambda(\setr)$ unless $\KK$ is locally compact.  In
general, the finiteness of $\Lambda(\setr)$ follows immediately from
Lemma~\ref{lem:clf} if $\neutral$ does not belong to the $\KDU$-weak
closure of $\setr$, since $\Lambda(\setr)\leq
\lambda(\cl(\imath(\setr)))$. In this case $\Lambda$ is the product of
the measure $\theta_\alpha$ given by (\ref{eq:defm}) with $\alpha\neq
0$ and a finite spectral measure $\sigma$ on $\SS$.

If $\Lambda(\setr)$ is finite, the Poisson process $\Pi$ on $\KK$ with
intensity measure $\Lambda$ \as has only a finite number of
points in $\setr$ for any $r>0$.  Define
\begin{equation}
  \label{eq:xir}
  \xir=\int_{\setr} x\,\Pi(dx)\,, \quad r> 0\,.
\end{equation}
Similarly to Theorem~\ref{thr:sec-no-pont}, a suitable normalisation
needed to ensure the convergence of $\xir$ as $r \downarrow 0$ is
provided by integrals of the L\'evy function.

\begin{Th}
  \label{thr:sumnorm}
  Let $\xi$ be a proper $\sas$ random element in a normed cone $\KK$
  without Gaussian component such that its L\'evy measure $\Lambda$ is
  supported by $\KK$. Assume that \textbf{(C)} holds and
  $\Lambda(\setr)<\infty$ for all $r>0$.  
  \begin{itemize}
  \item[\textsl{(i)}] If $\Pi$ is a Poisson process on $\KK$ with
    intensity measure $\Lambda$ and
    \begin{equation}
      \label{eq:gammar-1}
      \gamma_r(\chi)=\exp\left\{-i\left(\ell(\chi)
          -\int_{\setr}L(\chi,x)\Lambda(d x)\right)\right\}\,,\quad 
      \chi\in\KDU\,, \; r>0\,,
    \end{equation}
    then
    \begin{equation}
      \label{eq:gamma-conv}
      \gamma_r(\chi)\E\chi(\xir)\to \E\chi(\xi)\quad \text{as}
      \quad r\downarrow 0\,.
    \end{equation}
  \item[\textsl{(ii)}] If $\KDU$ consists of continuous characters
    and, for each $r>0$, there exists $x_r\in\KK$ such that
    $\gamma_r(\chi)=\chi(x_r)$ for all $\chi\in\KDU$, and if
    $x_r+\xir$ converges almost surely as $r\downarrow 0$, then the
    limit coincides in distribution with~$\xi$.
  \end{itemize}
\end{Th}
\begin{proof}
  The first statement is proved similarly to
  Theorem~\ref{thr:sec-no-pont}.  By (\ref{eq:gamma-conv}),
  $\E\chi(x_r+\xir)\to\E\chi(\xi)$. If $x_r+\xir$ \as converges to
  $\zeta$, then its limit shares the Laplace transform with~$\xi$.
\end{proof}

If $\gamma_r(\chi)=1$ for all $\chi$ (\eg if $\xi$ is symmetric), then
(\ref{eq:gamma-conv}) implies that
\begin{displaymath}
  \E\chi(\xir)\to \E\chi(\xi)\quad \text{as}
    \; r\downarrow 0\,.
\end{displaymath}
If the characters from $\KDU$ are not necessarily continuous, the
convergence of Laplace transforms implies that $\xir$ weakly
converges to $\xi$ with respect to $\KDU$-weak topology, see
\cite[Prop.~IV.3.3]{vak:tar:cho87}.  However, in general this does not
suffice to show that $\xir$ weakly converges to $\xi$.  

\subsection{LePage series representation on semigroups}
\label{sec:lepage-seri-repr-1}

The almost sure convergence of $\xir$ from (\ref{eq:xir}) to a random
element $\tilde\xi=\int x\Pi(dx)$ does not necessarily mean that
$\tilde\xi$ and $\xi$, from which the $\xir$'s are derived, share the
same distribution.  Indeed, we only have that
$\E\chi(\xir)\to\E\chi(\xi)$, while $\E\chi(\xir)$ does not
necessarily converge to $\E\chi(\tilde\xi)$ if $\chi$ is not
continuous.  Indeed, discontinuous characters cannot be simply
interchanged with infinite sums of elements from $\KK$.  The following
definition singles out characters that are interchangeable with sums
of series in $\KK$.

\begin{definition}
  \label{def:ser-cont}
  A character $\chi$ is said to be \emph{series continuous} if
  \begin{equation}
    \label{eq:sc}
    \chi\Bigl(\sum_{k=1}^\infty x_k\Bigr)=\prod_{k=1}^\infty \chi(x_k)
  \end{equation}
  for every convergent series $\sum x_k$ of elements from $\KK$.
\end{definition}

If the characters from $\KDU$ are series continuous, then
Theorem~\ref{thr:sumnorm}\textsl{(ii)} holds without assuming the
continuity of the characters. The following result shows the
uniqueness of the ingredients of the LePage series.

\begin{Th}
  \label{thr:unique}
  If \textbf{(C)} holds with $\KDU$ that consists of series continuous
  characters, then any two \as convergent (as principal values)
  LePage series (\ref{eq:9}) with sums having proper distributions are
  identically distributed if and only if the corresponding parameters
  $\alpha$ and the spectral measures coincide.
\end{Th}
\begin{proof}
  The LePage series (\ref{eq:9}) with a proper sum can be written as
  $\int x\Pi(dx)$ for the Poisson process $\Pi$ on $\KK$ with the
  intensity measure $\Lambda$. If $\xi$ admits two different LePage
  representations, this means that the Laplace exponent of $\xi$ has
  two representations with different L\'evy measures. The uniqueness
  of (\ref{eq:ldcomplex}) however implies that this is impossible.
\end{proof}

Consider relatively simple semigroups that are embeddable in a certain
group $\GG$ such that $\GG$-valued proper $\sas$ random elements admit
the LePage series representation. For instance, this is the case for
semigroups embeddable as cones in a Banach space $\BB$, see
Theorem~\ref{prop:4}. Then the following result holds if $\BB$ has a
separable dual space.

\begin{Th}
  \label{thr:cone-embed}
  Let $\GG$ be a group such that each proper $\sas$ $\GG$-valued
  random element with $\alpha\in(0,1)$ admits the LePage
  representation. Assume that there exists a countable family of
  continuous homomorphisms $f_n:\GG\to\R$, $n\geq1$, such that
  $\KK=\{x\in\GG:\; f_n(x)\geq0,\, n\geq1\}$.  Then each proper $\sas$
  random element in $\KK$ with $\alpha\in(0,1)$ admits the LePage
  series representation (\ref{eq:9}).
\end{Th}
\begin{proof}
  A $\sas$ random element $\xi$ in $\KK$ is also $\sas$ in $\GG$.  By
  the condition, $\xi$ can be represented by (\ref{eq:9}) where
  $\eps_k$ are distributed on the unit sphere in $\GG$. Now consider
  any homomorphism $f:\GG\to\R$, \ie a linear continuous functional on
  $\GG$, such that $f(x)\geq0$ for each $x\in\KK$. By applying $f$ to
  both sides of (\ref{eq:9}), using the fact that $f(\xi)$ is $\sas$
  random variable and the uniqueness of the LePage representation (see
  Theorem~\ref{thr:unique}), we obtain that $f(\eps_k)\geq 0$ \as
  Repeating this for a countable family of linear continuous
  functionals, yields that $\eps_k\in\KK$ \as
\end{proof}

In general, the main prerequisites for the existence of the LePage
representation are the fact that L\'evy measures are supported by
$\KK$ and finite on $\setr$, the continuity (or series continuity) of
the characters and the convergence of the LePage series. The latter is
particularly simple to derive if the norm is sub-invariant and
$\alpha\in(0,1)$, see Theorem~\ref{th:1}.

\begin{Th}
  \label{thr:alpha01}
  Assume that $\KK$ has a sub-invariant norm and one of the following
  cases is applicable.
  \begin{itemize}
  \item[\textsl{(i)}] Conditions~\textbf{(C)}, \textbf{(E)} and
    \textbf{(S)} hold with $\KDU$ that consists of series continuous
    characters and the unit sphere $\SS$ is compact in $\KK$.
  \item[\textsl{(ii)}] Conditions~\textbf{(C)} and \textbf{(W)} hold
    with $\KDU$ that consists of continuous characters and the unit
    sphere $\SS$ is $\KDU$-weak compact.
  \end{itemize}
  Then every $\sas$ random element in $\EE$ with $\alpha\in(0,1)$ and
  without Gaussian component can be represented as the LePage series
  (\ref{eq:9}) with $z=\neutral$.
\end{Th}
\begin{proof}
  First of all, note that $\zero=\neutral$ by
  Lemma~\ref{lemma:sub-invar} and the weak compactness of the unit
  sphere implies that the closed unit ball (being $\setr[1]^c$) is
  $\KDU$-weak compact.
  
  Next, continuity of characters in \textsl{(ii)} implies that
  \textbf{(E)} holds in this case too.  By Lemma~\ref{lemma:proper},
  $\xi$ has a proper distribution. By Theorem~\ref{thr:lc-k} or
  Corollary~\ref{cor:weak-comp}, the L\'evy measure $\Lambda$ is
  supported by $\KK$.  Let us show that also \textsl{(ii)} implies
  \textbf{(S)}. Since $\zero=\neutral$, it suffices to consider a
  sequence $\{x_n\}$ with norm greater than 1. If this sequence
  $\KDU$-weakly converges, \textbf{(W)} implies that $\sup\|x_n\|=c$
  is finite. Thus, the $x_n$'s belong to the set $F=\{x\in\KK:1\leq
  \|x\|\leq c\}$. The $\KDU$-weak compactness condition on $\SS$ and
  the continuity of the characters implies that $F$ is $\KDU$-weak
  compact and so $\KDU$-weak closed. The separation property implies
  the uniqueness of the weak limit, so that the weak limit of $x_n$
  also belongs to $F$, \ie this weak limit is not $\neutral$. Thus
  \textbf{(S)} holds in the both cases \textsl{(i)} and \textsl{(ii)},
  and Lemma~\ref{lem:clf} implies that $\Lambda(\setr)<\infty$ for all
  $r>0$.
  
  Let $\Pi_\alpha$ be the Poisson point processes with intensity
  measure $\Lambda$. By Theorem~\ref{th:1}, the sum of its points
  converges absolutely, which implies that the principal value of the
  integral $G_{\Pi_\alpha}(\chi)=\int_{\EE}(1-\chi(x))\Lambda(dx)$
  converges. The last integral being~(\ref{eq:fininteg}) allows us to
  set the L\'evy function to zero by Theorem~\ref{thr:no-ell}.
  
  Let $\xi_\alpha$ be the sum of the LePage series (\ref{eq:9}) with
  $z=\neutral$.  By (\ref{eq:pgfp}),
  \begin{displaymath}
    \E \chi(\xi_\alpha)=\exp\left\{-\int_\KK
      (1-\chi(x))\Lambda(dx)\right\}\,. 
  \end{displaymath}
  Writing a part of (\ref{eq:ldcomplex}) as above, we arrive at 
  \begin{displaymath}
    \E\chi(\xi)= e^{-(i\ell(\chi)+q(\chi))}\E\chi(\xi_\alpha)\,.
  \end{displaymath}
  It suffices to show that $\xi$ and $\xi_\alpha$ share the same
  Laplace transform.  The quadratic form $q$ vanishes, since $\xi$
  does not have a Gaussian component by the imposed condition.  The
  function $\ell$ clearly vanishes in case of either the identical
  involution, or symmetric $\xi$ or if the second distributivity law
  holds. Now we show that this also holds in general under the
  sub-invariance assumption for $\alpha\in(0,1)$.
  
  Note that $\ell$ is the only imaginary part of the Laplace exponent
  of $\xi$. By applying $\chi$ to the both sides of
  (\ref{eq:a1alph-a+b1-}) we obtain that
  \begin{displaymath}
    \E(a^{1/\alpha}\circ\chi)(\xi) \E(b^{1/\alpha}\circ\chi)(\xi) 
    = \E((a+b)^{1/\alpha}\circ\chi)(\xi)\,,
  \end{displaymath}
  whence,  using the additivity property of
  $\ell$ we have  for $a=b=1/2$ that
  \begin{displaymath}
    \ell((2^{-1/\alpha}\circ\chi)(2^{-1/\alpha}\circ\chi)
    =\ell(\chi)\,.
  \end{displaymath}
  This equality can be written shorter as $\ell(f(\chi))=\ell(\chi)$,
  where $f(\chi)$ is a character that acts on $x\in\KK$ as
  $(f(\chi))(x)=\chi(2^{-1/\alpha}(x+x))$. By iterating we obtain that
  \begin{equation}
    \label{eq:lell}
    \ell(f^k(\chi))=\ell(\chi)\,,
  \end{equation}
  where $(f^k(\chi))(x)=2^{-k/\alpha}S_{2^k})$ and $S_n$ is the sum of
  $n$ identical summands being $x$. By the sub-invariance property,
  $n^{-1}S_n$ has a norm bounded by $\|x\|$ which is finite for all
  $x$ in a normed cone. Since $\alpha<1$, we have
  $2^{-k/\alpha}S_{2^k}\to\zero=\neutral$. Being the only imaginary
  part of the continuous Laplace exponent, the function $\ell$ is also
  continuous with respect to pointwise convergence of its argument.
  Since $(f^k(\chi))(x)\to\chi(\neutral)=1$ for all $x$ and
  \textbf{(E)} holds, $f^k(\chi)$ pointwisely converges to $\one$.  By
  passing to the limit in (\ref{eq:lell}) and noticing that
  $\ell(\one)=0$, we obtain that $\ell(\chi)=0$ for all $\chi\in\KDU$.
\end{proof}                                

If $\alpha$ does not belong to $(0,1)$ or if $\KK$ does not have a
sub-invariant norm, an analogue of Theorem~\ref{thr:sumnorm} holds
assuming that the LePage series converges almost surely. In a special
case of semigroups with identical involution, we arrive at the
following result.

\begin{Th}
  \label{thr:id-invol}
  Assume that $\KK$ is a normed cone with identical involution such
  that conditions \textsl{(i)} or \textsl{(ii)} of
  Theorem~\ref{thr:alpha01} hold. Then every $\sas$ random element in
  $\EE$ without Gaussian component admits the LePage representation
  \eqref{eq:lepage-int} with $z=\neutral$, provided the principal
  value of the integral exists almost surely.
\end{Th}

If $\KK$ is idempotent, then we do not have to require that $\xi$ does
not have a Gaussian component, \cf Lemma~\ref{lem:q-idempo}. Since the
addition operation is not invertible in general, the presence of
Gaussian component in $\xi$ does not mean that $\xi$ can be decomposed
as the sum of a Gaussian element with an independent remainder.

\section{Examples}
\label{sec:examples}

\subsection{Cones with the second distributivity law}
\label{sec:examples-with-second}

These examples are the closest to ``conventional'' stable
distributions on the line or in Banach spaces. The cone $\KK$ is
typically a group or is embeddable in a group.  In particular, the
second distributivity law implies that $\zero=\neutral$ and
$\KK(\alpha)=\{\neutral\}$ for $\alpha\neq1$. If \textbf{(C)} holds,
then $\alpha\in(0,2]$ and the order of homogeneity of the L\'evy
measure is strictly less than $2$ (respectively, $1$ if the involution
is identical), see Corollary~\ref{cor:max-bound}.  By
Theorem~\ref{thr:l-and-q}, a proper $\sas$ random element does not
have a Gaussian component if $\alpha\neq1,2$ or $\alpha\neq2$ if $\KK$
is a group.

\begin{example}[Positive half-line]
  \label{ex:chrplus}
  Consider $\R_+$ with the arithmetic addition, Euclidean metric and
  identical involution. A separating family of continuous characters
  is given by $\chi_t(x)=e^{-tx}$ for $t\geq0$. By
  Corollary~\ref{cor:max-bound}, the L\'evy measures are homogeneous
  with the order strictly less than 1.  The only $\sas$ laws with
  $\alpha=1$ are deterministic distributions. Note
  that
  \begin{displaymath}
    \phi(\chi)=-\log \E e^{-t\xi}=ct^\alpha
  \end{displaymath}
  is a well-known representation of the Laplace transform of a
  strictly stable law on $\R_+$. The existence of the LePage
  representation follows from Theorem~\ref{thr:alpha01}\textsl{(i)}.
\end{example}

\begin{example}[Banach spaces]
  \label{ex:b-space}
  Let $\KK$ be a Banach space. A separating semigroup $\KDU$ of
  continuous characters consists of $\chi_u(\cdot)=e^{iu(\cdot)}$ for
  all linear continuous functionals $u$. Since the norm is clearly
  invariant, the LePage series converges absolutely for
  $\alpha\in(0,1)$.  Theorem~\ref{th:2} turns into
  \cite[Th.~7.11]{ara:gin}.
  
  The separation property holds by the Hahn-Banach theorem, so that
  \textbf{(C)} is satisfied. The $\KDU$-weak convergence is the
  conventional weak convergence in Banach spaces. The strong
  convergence implies the weak convergence, so that \textbf{(E)}
  holds, hence $\alpha\in(0,2]$.  If the space is reflexive, then the
  unit sphere is weak sequentially compact, \seg
  \cite[Th.~II.3.28]{dunford}.  Furthermore, the weak convergence
  implies the strong boundedness, so that \textbf{(W)} holds and by
  Corollary~\ref{cor:weak-comp}, the L\'evy measure is supported by
  $\KK$. Every $\sas$ random element with $\alpha\in(0,1)$ admits the
  LePage representation by Theorem~\ref{thr:alpha01}\textsl{(ii)}.
  For $1\leq\alpha\leq 2$, one should use the methods that rely on
  symmetrisation arguments, see \cite{linde}, \eg to check the
  tightness of $x_r+\xir$ needed in
  Theorem~\ref{thr:sumnorm}\textsl{(ii)}.
  
  Theorem~\ref{thr:cone-embed} can be applied to any cone which is
  embeddable in a Banach space with separable dual, \eg to the cone of
  non-negative continuous functions with addition operation.
\end{example}

\begin{example}[Compact convex sets with Minkowski addition]
  \label{ex:mconv}
  Let $\KK$ be the family $\co\sK$ of nonempty compact convex sets $K$
  in $\R^d$ with the semigroup operation being the Minkowski
  (elementwise) addition.  Note that $\KK$ is not a group, since the
  Minkowski addition is not invertible. The involution corresponds to
  the symmetry with respect to the origin.  The Hausdorff metric turns
  $\KK$ into a normed cone with $\neutral=\zero=\{0\}$ and the norm
  defined as $\|K\|=\sup\{\|x\|:\; x\in K\}$. Since the Hausdorff
  metric is invariant, \ie the distance between $K_1+L$ and $K_2+L$
  coincides with the distance between $K_1$ and $K_2$, the
  corresponding norm is clearly sub-invariant.  The unit sphere in
  $\KK$ is $\SS=\{K\in\co\sK:\; \|K\|=1\}$; it is compact in the
  Hausdorff metric.
  
  Before describing the characters on $\co\sK$, consider the
  appropriate sub-family of ``centred'' sets. With every convex compact
  set $K$ it is possible to associate its Steiner point
  \begin{displaymath}
    s(K)=\frac{1}{\kappa_d}\int_{S^{d-1}} h(K,u)u\,du\,,
  \end{displaymath}
  where $\kappa_d$ is the volume of the unit ball in $\R^d$, the
  integral is taken with respect to the
  $(d-1)$-dimensional Hausdorff measure on the unit sphere $S^{d-1}$
  in $\R^d$ and $h(K,u)$ is the support
  function of $K\in\co\sK$, \ie the supremum of the scalar product
  $\langle u,x\rangle$ over $x\in K$ for $u\in S^{d-1}$. It is known
  that the Steiner point is linear with respect to Minkowski addition
  and that $s(K)\in K$, see \cite[p.~42]{schn}. Let $\co\sK_0$ be the
  family of convex compact sets with their Steiner points located at the
  origin. This family can be considered as a convex cone with the
  identical involution and the Hausdorff metric. The condition
  \textbf{(C)} holds with a separating family of continuous characters
  given by
  \begin{equation}
    \label{eq:charact-mink-null}
    \chi_\nu(K)=\exp\Bigl\{-\int_{S^{d-1}} h(K,u)\,\nu(du)\Bigr\}\,,
  \end{equation}
  where $\nu$ is a finite measure on $S^{d-1}$. In fact, it suffices
  to consider only counting measures~$\nu$. Note that it does not
  suffice to consider only the characters of the type $e^{-h(K,u)}$,
  since they do not form a closed family with respect to
  multiplication. 
  
  A random element $X$ in $\KK=\co\sK$ is a \emph{random compact set},
  \seg \cite[Sec.~1.5]{mo1}. All non-trivially equal to $\{0\}$ random
  compact sets are proper.  By Theorem~\ref{thr:bound}, a proper
  $\sas$ random element in $\co\sK_0$ necessarily has
  $\alpha\in(0,1]$. Otherwise $X$ is not proper, \ie $X=\{0\}$. If
  $\alpha=1$, then by Corollary~\ref{cor:order2} the corresponding
  Laplace exponent is a linear continuous function of $\chi$ and the
  random set $X$ is deterministic. A conventional argument for this
  relies on the fact that the support function in every direction is a
  non-negative $\sas$ random variable, which is necessarily
  degenerated for $\alpha\geq1$, see \cite{gin:hah85b}.  Similar
  arguments are applicable to the cone of convex compact sets that
  contain the origin (not necessarily having the Steiner point at the
  origin) and are equipped with the identical involution.
  Theorem~\ref{thr:alpha01}\textsl{(i)} applies, implying the LePage
  representation of $\sas$ random elements.
  
  Now return to the cone of $\KK=\co\sK$ of all convex compact sets.
  A convex compact set $K$ can be decomposed as $K=s(K)+K_0$.
  Therefore, the cone of compact convex sets can be decomposed into
  the sum of two cones: $\R^d$ with the conventional addition and the
  cone $\co\sK_0$ of convex sets with the Steiner point at the origin.
  The first has the complex-valued characters, while the second has
  the identical involution and $[0,1]$-valued characters given by
  (\ref{eq:charact-mink-null}).  By combining the two families of
  characters we obtain the separating family of continuous characters
  given by
  \begin{equation}
    \label{eq:charact-mink}
    \chi_{v,\nu}(K)=e^{i\langle v,s(K)\rangle}
    \exp\Bigr\{-\int_{S^{d-1}} h(K-s(K),u)\,\nu(du)\Bigl\}\,,
  \end{equation}
  where $v\in\R^d$. If $X$ is a $\sas$ random convex compact set, then
  $s(K)$ is a $\sas$ random vector, so that $\alpha\in(0,2]$.
  Furthermore, $X-s(X)$ is $\sas$ in $\co\sK_0$, whence
  $\alpha\in(0,1]$. Therefore, a stable random set with
  $\alpha\in[1,2]$ is a sum of an $\sas$ random vector $\xi$ in $\R^d$
  and a deterministic convex compact set, which is necessarily $\{0\}$
  in case $\alpha>1$, \ie $X=\{\xi\}$. If the stability definition
  (\ref{eq:a1alph-a+b1-}) is weakened by allowing an additive
  normalisation, then this deterministic part may be any convex
  compact set, see \cite{gin:hah85b}. 
  
  Theorem~\ref{th:rep} implies that the Laplace exponent
  $\phi(v,\nu)=\E\chi_{v,\nu}(X)$ of a $\sas$ random set $X$ satisfies
  \begin{displaymath}
    \phi(tv,t\nu)=t^\alpha\phi(v,\nu)
  \end{displaymath}
  for all $t>0$, $v$ and $\nu$.  Since the unit sphere $\SS$ is
  compact and \textbf{(S)} holds, every $\sas$ random compact set with
  $\alpha\in(0,1)$ admits the absolutely convergent LePage series
  representation.
  
  Not necessarily strictly stable random compact convex sets for any
  $\alpha\in(0,2]$ in a separable Banach space have been studied in
  \cite{gin:hah85b}.  It is shown in \cite[Th.~1.14]{gin:hah85b} that
  $\alpha$-stable random compact convex sets with $\alpha\in(0,1)$ can
  be represented as a stochastic integral over $\SS$ with respect to a
  certain independently scattered $\alpha$-stable random measure. This
  representation coincides with the LePage representation in $\KK$.

  It should be noted that general semigroups representable as systems
  of compact convex sets have been characterised in \cite{rat:sch77}.
\end{example}

\begin{example}[Upper semicontinuous functions]
  \label{ex:usc}
  Let $\KK$ be the family of upper semicontinuous functions
  $u:\R^d\mapsto[0,1]$ such that $\cup_{t>0} F_t(u)$ is relatively
  compact where $F_t(u)=\{x:\; u(x)\geq t\}$. The metric between $u$
  and $v$ is defined as the supremum of the Hausdorff distances
  between $F_t(u)$ and $F_t(v)$ over $t\in(0,1]$. The sum of $u$ and
  $v$ is defined to be the upper semicontinuous function $w$ such that
  $F_t(w)$ equals the Minkowski sum of $F_t(u)$ and $F_t(v)$ for all
  $t\in(0,1]$. This setting is similar to Example~\ref{ex:mconv}
  reformulated for increasing set-valued function indexed by $[0,1]$,
  see also \cite{mo99o}. Accordingly, the conclusions of
  Example~\ref{ex:mconv} can be transferred to this case practically
  without changes.
\end{example}

\begin{example}[Finite random measures]
  \label{ex:mu1}
  Let $\KK$ be the family of all finite measures $m$ on a locally
  compact topological space $E$ with a countable base. The operations
  are the conventional addition of measures and the multiplication of
  their values by numbers. The neutral element and the origin are both
  equal to the null-measure. The Prohorov metric (\seg \cite{rac91})
  on $\KK$ is sub-invariant and homogeneous at the origin. The
  corresponding norm $\|m\|$ is the total mass of $m$.  Assume that
  the involution is identical. A separating family of characters is
  given by
  \begin{equation}
    \label{eq:muca}
    \chi_{u}(m)=\exp\Bigl\{-\int u\,dm\Bigr\}
  \end{equation}
  for any continuous bounded function $u:E\mapsto\R_+$. Since the
  Prohorov metric metricises the weak convergence, these characters
  are continuous and the $\KDU$-weak convergence is equivalent to the
  metric convergence. The Laplace transform $\E\chi_u(\mu)$ is called
  the Laplace functional of the random measure $\mu$, \seg
  \cite[(6.4.16)]{da:vj}.  The corresponding Laplace exponent is given
  by
  \begin{equation}
    \label{eq:mupsi}
    \phi(u)=-\log\Bigl(\E\exp\Bigl\{-\int u\, d\mu\Bigr\}\Bigr)\,.
  \end{equation}
  All non-trivial random measures have proper distributions.  The
  random measure $\mu$ is $\sas$ if
  \begin{equation}
    \label{eq:mu1}
    a^{1/\alpha}\mu_1(K)+b^{1/\alpha}\mu_2(K)\deq
    (a+b)^{1/\alpha}\mu(K)
  \end{equation}
  for all measurable $K$.  By Theorem~\ref{thr:bound},
  $\alpha\in(0,1]$, which is also due to the fact that $\mu(K)$ is a
  non-negative $\sas$ random variable, \cf \cite{sam:taq94}.  Note
  that $c\circ\chi_{u}$ corresponds to the character $\chi_{cu}$\,.
  By Theorem~\ref{th:rep}, $\mu$ is $\alpha$-stable with some
  $\alpha\in(0,1]$ if and only if its Laplace exponent $\phi(u)$
  satisfies $\phi(su)= s^\alpha\phi(u)$ for all $u$ and $s>0$. The
  integral representation (\ref{eq:ldc}) corresponds to the
  representation for the Laplace functional of infinitely divisible
  random measures, see \cite[Prop.~9.2.VII]{da:vj}. The LePage series
  involves random measures distributed on the unit sphere $\SS$ in
  $\KK$, which is the family of all probability measures on $E$. Thus,
  a $\sas$ measure is the weighted sum of \iid random probability
  measures with weights $\Gamma_k^{-1/\alpha}$, $k\geq1$. This
  representation is a particular case of \cite[Th.~3.9.1]{sam:taq94}
  for strictly $\alpha$-stable random measures.

  If $E$ is a compact space, then the unit sphere $\SS$ is compact, so
  that each $\alpha$-stable random measure with $\alpha\in(0,1)$
  admits the LePage representation by
  Theorem~\ref{thr:alpha01}\textsl{(i)}.
\end{example}

\begin{example}[Locally finite random measures]
  \label{ex:lfrm}
  Let $\KK$ be the family of locally finite measures on $\R^d$ with
  the topology of vague convergence and the same operations as in
  Example~\ref{ex:mu1}. This space is important in the studies of
  point processes without accumulation points.
  
  A separating family of continuous characters is given by
  (\ref{eq:muca}), but with $u$ being a measurable bounded function
  with bounded support. The conditions \textbf{(C)} and \textbf{(E)}
  hold, and all non-trivial random measures are proper. Any $\sas$
  random measure $\mu$ has the Laplace exponent that is homogeneous,
  \ie $\phi(su)=s^\alpha\phi(u)$, where $\alpha\in(0,1]$ by
  Theorem~\ref{th:rep}.
  
  The extension of the Prohorov metric typically used to metricise the
  vague convergence of locally finite measures (\seg
  \cite[Sec.~A2.6]{da:vj} or \cite[Sec.~10.2]{rac91}) is not
  homogeneous, since it is constructed using the sums of the type
  $2^{-i}d_i/(1+d_i)$, where $d_i$ is a certain distance between the
  restrictions of the measures onto the balls $B_i$, $i\geq1$.
  
  The unit sphere and the LePage series can be constructed for a
  sub-family of $\KK$ that consists of locally finite measures $m$
  with a certain growth restriction at the infinity, \eg those which
  satisfy
  \begin{displaymath}
    \int_0^\infty e^{-r} m(B_r)dr <\infty\,.
  \end{displaymath}
  The value of this integral may serve as a norm of $m$ with the
  corresponding sub-invariant metric given by the integral of the
  Prohorov distance between the measures restricted on $B_r$. In
  particular, scale the norm if necessary so that atomic measures
  $\delta_x$ belongs to the unit sphere $\SS$ in $\KK$ for all $x$
  from the unit sphere $S^{d-1}$ in $\R^d$. Then the LePage series is
  $\mu_\alpha=\sum_k \Gamma_k^{-1/\alpha}\delta_{\eps_k}$, where
  $\{\eps_k,\, k\geq1\}$ are \iid random variables on $S^{d-1}$.  Note
  that this representation does not coincide with the representation
  of a stable Poisson process $\Pi_\alpha$ from Theorem~\ref{prop:2}.
  The random measure $\mu_\alpha$ is supported by a subset of the unit
  sphere $S^{d-1}$. It does not take integer values and has the total
  mass being the $\sas$ random variable $\xi_\alpha=\sum_k
  \Gamma_k^{-1/\alpha}$ in $(\R_+,+)$.
\end{example}

\subsection{Cones without the second distributivity law}
\label{sec:exampl-with-second}

These examples typically appear if $\KK$ is not a group, but only a
semigroup. In these cases one might have positive or negative
stability parameter $\alpha$ unless an element of $\KK$ possesses a
non-trivial inverse, which implies, under conditions \textbf{(C)} and
\textbf{(E)}, that $\alpha$ is necessarily positive, see
Theorem~\ref{th:alpha}.

\begin{example}[Positive cone in Euclidean space with coordinatewise
  maximum] 
  \label{ex:max}
  For some $d\geq1$, let $\KK$ be $[0,\infty)^d$ with the
  coordinatewise maximum operation, \ie $x\vee y=(x_1\vee
  y_1,\dots,x_d\vee y_d)$, and the conventional multiplication by
  numbers. In this case $\zero=\neutral=0$ and $\KK(\alpha)=\{0\}$ for
  all $\alpha\neq0$. The corresponding $\sas$ laws are called
  \emph{max-stable}, see \cite{gal78,kot:nad00}.
  
  The semigroup $\KK$ has the identical involution and moreover is
  idempotent, so that the characters take only values $0$ or $1$, and
  the Gaussian component is degenerated by Lemma~\ref{lem:q-idempo}.
  The Euclidean metric is sub-invariant.  A separating family of
  series continuous (and also upper semicontinuous) characters is
  given by $\chi_z(x)=\one_{[0,z]}(x)$ for
  $z=(z_1,\dots,z_d)\in(0,\infty)^d$, where
  $[0,z]=[0,z_1]\times\cdots\times[0,z_d]$.  If $z$ has rational
  coordinates, these characters form a countable strictly separating
  family such that \textbf{(C)} and \textbf{(E)} hold.  The
  multiplication operation acts on characters as
  $a\circ\one_{[0,z]}=\one_{[0,z/a]}$.  The Laplace transform of $\xi$
  is then
  \begin{displaymath}
    \E\one_{[0,z]}(\xi)=\P\{\xi\leq z\}\,,
  \end{displaymath}
  hence any non-trivial $\xi$ with the support not separated from zero
  is proper.  By Theorem~\ref{th:alpha}, the corresponding parameter
  $\alpha$ is always positive.
  
  The unit sphere is compact and \textbf{(S)} holds.  The L\'evy
  function and $\ell$ vanish, so that $\gamma_r=1$, see
  (\ref{eq:gammar-1}).  By Theorem~\ref{thr:alpha01}\textsl{(i)}, each
  $\sas$ random element with $\alpha\in(0,1)$ admits the LePage
  representation. The LePage series absolutely converges \as for all
  $\alpha>0$ and yields a max-stable random element.  The
  corresponding representation of a max-stable law as the maximum of
  Poisson points is well known in the theory of multivariate extremes,
  see \cite[Ch.~5]{res87}.
\end{example}

\begin{example}[Half-line with minimum operation]
  \label{ex:rplus-min}
  Consider the extended half-line $[0,\infty]$ with the minimum
  operation and the conventional multiplication by numbers. Then
  $\zero$ is conventional zero $0$, while $\neutral=\infty$ and so has
  infinite norm.  The upper semicontinuous characters are given by
  $\one_{[z,\infty]}(x)$ for $z\in[0,\infty]$.  They are also series
  continuous and form a strictly separating family, so that
  Theorem~\ref{th:rep} applies.  The $\sas$ random elements
  necessarily have $\alpha<0$ and the corresponding stable Poisson
  process $\Pi_\alpha$ has realisations from $\sM$.  The LePage series
  $\min_{k\geq1}(\Gamma_k^{-1/\alpha}\eps_k)$ with $\eps_k=1$
  evidently converges to $\Gamma_1^{-1/\alpha}$, so that each $\sas$
  random element has the LePage representation by
  Theorem~\ref{thr:id-invol}.
  
  In the same way it is possible to handle the coordinatewise minimum
  operation that leads to multivariate min-stable laws,
  \cf Example~\ref{ex:max}.
\end{example}

\begin{example}[Union-stable random compact sets]
  \label{ex:usrs}
  Let $\KK$ be the family $\sK$ of compact subsets $K\subset\R^d$ with
  the Hausdorff metric and the union operation that turns $\KK$ into
  an idempotent semigroup. The multiplication by numbers is defined as
  the corresponding homothetical transformation of sets. The neutral
  element of $\KK$ is the empty set, while $\zero$ is the origin
  $\{0\}$. Note that the norm of the sum is always larger than the
  norm of each summand, so that every non-trivial $\sas$ law should
  have $\alpha>0$ by Lemma~\ref{lem:alpha-norm}.
  
  A separating family $\KDU$ of series continuous characters consists
  of $\chi_G(K)=\one_{G\cap K=\emptyset}$ for all open sets $G$.  It
  is possible to extract its countable separating sub-family by
  considering $G$ that are unions of finite number of open balls with
  rational centres and radii, \ie from the so-called separating class
  \cite[Sec.~1.4]{mo1}.  The product of two characters $\chi_G$ and
  $\chi_{G'}$ is the character $\chi_{G\cup G'}$, so that these
  characters indeed build a semigroup $\KDU$. 
  
  The Laplace transform of a random compact set $X$ is given by
  \begin{displaymath}
    \E\chi_G(X)=\P\{X\cap G=\emptyset\}\,,
  \end{displaymath}
  which is usually called \emph{the avoidance functional} of $X$ and
  denoted by $Q_X(G)$. If $X$ is not proper, then $Q_X(G)=0$ for some
  $G$, \ie $X$ almost surely hits $G$. This implies that $X$ has a
  fixed point, \ie there exists $x\in\R^d$ such that $x\in X$ a.\ssp
  s., see \cite[Lemma~4.1.8]{mo1}. Thus proper random elements are
  exactly those that do not possess fixed points. Conditions
  \textbf{(C)} and \textbf{(E)} hold.  By Theorem~\ref{th:rep}, the
  union-stable random compact sets can be characterised as those
  having homogeneous Laplace exponents $\phi(G)=-\log Q_X(G)$, \ie
  $\phi(sG)=s^\alpha\phi(G)$.  This characterisation has been obtained
  in \cite[Th~4.1.12]{mo1} using direct proofs, but in a more general
  case for not necessarily compact random closed sets that possibly
  possess fixed points, \ie with not necessarily proper distributions.

  Since $\neutral\neq\zero$, Theorem~\ref{th:alpha} implies that
  $\alpha<0$, contrary to the above mentioned conclusion of
  Lemma~\ref{lem:alpha-norm}. Thus, this cone does not possess any
  non-trivial \emph{proper} $\sas$ law.  Indeed, the corresponding
  LePage series would involve the union of random compact sets scaled
  by arbitrarily large factors $\Gamma_k^{-1/\alpha}$, so that this
  union is no longer compact if $\alpha<0$ and this union contains the
  origin (and thereupon is not proper) if $\alpha>0$.
\end{example}

\begin{example}[Union-stable random compact sets containing the origin]
  \label{ex:usrs-fp}
  Let $\KK$ be the cone of compact sets in $\R^d$ that contain the
  origin. We keep the same cone operations as in
  Example~\ref{ex:usrs}. However in this case both the neutral element
  and the origin are $\{0\}$, so that $\alpha>0$ by
  Theorem~\ref{th:alpha} and also by Lemma~\ref{lem:alpha-norm}. Since
  all elements of $\KK$ contain the origin, the characters are given
  by $\one_{G\cap X=\emptyset}$ indexed by open sets $G$ that do not
  contain the origin.
  
  The Hausdorff metric on $\sK$ is sub-invariant with the unit sphere
  $\SS=\{K\in\sK:\; \|K\|=1\}$. The LePage series is
  \begin{equation}
    \label{eq:lp-union}
    \bigcup_{k\geq1} \Gamma_k^{-1/\alpha} X_k\,,
  \end{equation}
  where $\{X_k,\, k\geq1\}$ is any sequence of \iid random compact
  sets from $\SS$. This series converges for all $\alpha>0$, although
  its convergence for $\alpha\geq1$ does not follow from
  Theorem~\ref{th:1}, but is easily seen by the direct proof. If $X_k$
  are independent singletons $\eps_k$ independently distributed on the
  unit sphere in $\R^d$, then $X$ is the support of the stable Poisson
  process $\Pi_\alpha$ with $\alpha>0$.
  
  Theorem~\ref{thr:lc-k} is applicable, so that the L\'evy measures
  are supported by $\KK$ and every proper $\sas$ random element admits
  the LePage representation for each $\alpha>0$, see
  Theorem~\ref{thr:id-invol}. Note that the Gaussian component is
  always degenerated by Lemma~\ref{lem:q-idempo}.
\end{example}

\begin{example}[Continuous functions with addition and argument
  rescaling]
  \label{ex:rescaled-functins}
  Let $\KK$ be the family of continuous functions $f:\R^d\mapsto\R$
  such that $f(0)=0$ and $\|x\|^{-1}|f(x)|$ is bounded in
  $x\in\R^d\setminus\{0\}$.  The cone operations are the arithmetic
  addition and the rescaling of the argument, \ie $(D_af)(x)=f(ax)$
  for all $x$. Define invariant and homogeneous metric by
  \begin{equation}
    \label{eq:norm-resc}
    d(f_1,f_2)=\sup_{x\in\R^d\setminus\{0\}}\|x\|^{-1} |f_1(x)-f_2(x)|\,.
  \end{equation}
  Then $\|f\|$ is the supremum of $\|x\|^{-1}|f(x)|$, $x\neq 0$. The
  neutral element is the zero function, while the imposed condition 
  $\|f\|<\infty$ implies that $\zero=\neutral$. The involution is
  given by $f^\star=-f$. It is interesting to note that $\KK$ is a
  group, where nonetheless the second distributivity law does not
  hold.  Since the metric is invariant, the LePage series converges
  absolutely for $\alpha\in(0,1)$.
  
  Define continuous characters on $\KK$ as
  \begin{equation}
    \label{eq:char-f}
    \chi_\nu(f)=\exp\Bigl\{i \int f d\nu\Bigr\}\,,
  \end{equation}
  where $\nu$ is a signed measure on $\R^d$ with bounded support and
  finite total variation.  Their continuity follows from the fact that
  the convergence in norm implies the pointwise convergence, while the
  pointwise convergence is equivalent to the $\KDU$-weak convergence.
  Then $c\circ\chi_\nu=\chi_{c\circ\nu}$, where
  $(c\circ\nu)(A)=\nu(c^{-1}A)$ for measurable $A$. The Laplace
  exponent $\phi$ of a $\KK$-valued random element $\xi$ is given by
  $\E\chi_\nu(\xi)=e^{-\phi(\nu)}$. Conditions~\textbf{(C)} and
  \textbf{(E)} hold, so that Theorem~\ref{th:rep} yields that $\sas$
  elements in $\EE$ are characterised by having homogeneous Laplace
  exponents, \ie $\phi(c\circ\nu)=c^\alpha\phi(\nu)$, where $\alpha>0$
  by Theorem~\ref{th:alpha}.
  
  The unit ball in $\KK$ is $\KDU$-weak compact, since
  $\|x\|^{-1}f_n(x)\leq 1$, $n\geq1$, implies that $\|x\|^{-1}f_n(x)$
  has a pointwisely convergent subsequence, whence $f_n$ has a
  pointwise convergent subsequence too. However Condition \textbf{(W)}
  does not hold. Neither is the unit sphere compact in the metric
  topology on $\KK$. Thus the available results of
  Section~\ref{sec:support-levy-measure} do not lead to the conclusion
  that the L\'evy measure is supported by $\KK$.
\end{example}

\begin{example}[Integrable random probability measures with 
  convolution operation]
  \label{ex:rprobm}
  Let $\KK$ be the family of probability measures $m$ on $\R$ with
  finite first moment, the addition operation being the convolution of
  measures (denoted by $*$) and the multiplication by scalars given by
  $(D_a m)(K)=m(a^{-1}K)$.  Assume that the corresponding random
  variables are defined on a non-atomic probability space.  If $\xi_1$
  and $\xi_2$ are independent random variables with distributions
  $m_1$ and $m_2$ respectively, then $(D_{c_1}m_1)*(D_{c_2}m_2)$ is
  the distribution of $c_1\xi_1+c_2\xi_2$. The unit measure $\delta_0$
  concentrated at $0$ is both the neutral element of $\KK$ and the
  origin.  The second distributivity law does not hold in this case,
  since the convolution of two measures (\ie the distribution of the
  sum of two \iid random variables) does not generally equal the
  rescaled measure (\ie the distribution of the rescaled random
  variable). The involution corresponds to the central symmetry, \ie
  $m^\star$ is the distribution of $-\xi$ if $m$ is the distribution
  of $\xi$.
  
  Note that $\KK(\alpha)$ are all non-trivial for $1< \alpha\leq 2$
  consisting of $\sas$ probability distributions in $\R$. The
  integrability condition (that has not been used so far) implies the
  absence of non-trivial elements of $\KK(\alpha)$ for
  $\alpha\in(0,1]$. Indeed, the corresponding $\sas$ laws are
  non-integrable, \cf Example~\ref{ex:all-prob-m}.
  
  The problem of defining a sub-invariant norm on the family of
  probability measures can be easily solved using the tools from the
  theory of probability metrics, see \cite{rac91,zol}. Recall that a
  probability metric $d$ defines a distance between probability
  measures $m_1$ and $m_2$ or between the corresponding random
  variables $\xi$ and $\eta$. The metric is called \emph{simple}, if
  it does not take into account the joint distribution of these random
  variables. A metric $d$ is called $r$-homogeneous if
  $d(c\xi,c\eta)=|c|^rd(\xi,\eta)$ for all $c\neq0$; this metric is
  said to be regular if $d(\xi+\zeta,\eta+\zeta)\leq d(\xi,\eta)$ for
  any random variable $\zeta$ that is independent of both $\xi$ and
  $\eta$. The regularity property written using the addition in $\KK$
  (\ie the convolution) turns into (\ref{triangleq}) meaning that the
  simple regular metric is sub-invariant. In the theory of probability
  metrics, an $r$-homogeneous and regular simple metric is called
  \emph{ideal} of order~$r$.
  
  An ideal probability metric of order $1$ between distributions of
  random variables $\xi$ and $\eta$ can be defined as
  \begin{displaymath}
    d(\xi,\eta)=\sup\{|\E(f(\xi)-f(\eta))|:\; f\in\mathrm{Lip}_1\}\,,
  \end{displaymath}
  where $\mathrm{Lip}_1$ is the family of functions with Lipschitz
  constant 1, \ie $|f(x)-f(y)|\leq |x-y|$ for all $x$ and $y$.  This
  is the so-called \emph{Kantorovich metric} or \emph{Wasserstein
    metric}, \seg \cite[Th.~1.3.1]{zol}.  Note also that
  \begin{displaymath}
    d(\xi,\eta)=\int_{-\infty}^\infty |\P\{\xi\leq x\}-\P\{\eta\leq
    x\}| dx\,,
  \end{displaymath}
  so that the corresponding norm of $\xi$ (or of the measure $m$) is
  given by
  \begin{displaymath}
    \|\xi\|=\|m\|=\int_{-\infty}^0 \P\{\xi\leq x\}dx +
    \int_0^\infty \P\{\xi>x\} dx=\E|\xi|\,.
  \end{displaymath}  
  The norm of any integrable measure is clearly
  finite. The Kantorovich metric is complete, see
  \cite[Th.~6.3.3]{rac91}. Furthermore $d(\xi_k,\xi)\to0$ if and only
  if $\xi_k\wto\xi$ weakly and $\E|\xi_k|\to\E|\xi|$, see
  \cite[Th.~3.9.4]{kal:rac90}. 

  Characters on $\KK$ are the characteristic functions, \ie
  \begin{displaymath}
    \chi_u(m)=\int e^{iux}m(dx)\,, \quad u\in\R\,.
  \end{displaymath}
  They form a separating family of continuous characters with values
  in the unit complex disk $\DD$, but they do not build a
  semigroup, since their products are not necessarily of the same
  type. However, if one defines 
  \begin{equation}
    \label{eq:char-uuu}
    \chi_{u_1,\dots,u_k}(m)=\prod_{i=1}^k \chi_{u_i}(m)
  \end{equation}
  for any finite set $\{u_1,\dots,u_k\}\subset \R$, $k\geq1$, then
  this family does constitute a semigroup.  Note that
  $a\circ\chi_{u_1,\dots,u_k}=\chi_{au_1,\dots,au_k}$.  Conditions
  \textbf{(C)} and \textbf{(E)} hold, since the characters are
  continuous. The Laplace transform of a random measure $\mu$ from
  $\KK$ is given by
  \begin{displaymath}
    \E \chi_{u_1,\dots,u_k}(\mu)=
    \idotsint e^{i(u_1x_1+\cdots+u_kx_k)}
    \overline{\mu}_k(dx_1\times\cdots\times dx_k)\,,
  \end{displaymath}
  where $\overline{\mu}_k$ is the $k$th order moment measure of the
  random measure $\mu$, \seg \cite{da:vj,m:k:m}. In other words, the
  Laplace exponent 
  \begin{math}
    \phi(u_1,\dots,u_k)=-\log(\E \chi_{u_1,\dots,u_k}(\mu))
  \end{math}
  is the log-characteristic function of the $k$th moment measure
  $\overline{\mu}_k$.   A random probability measure $\mu$ is $\sas$ if
  \begin{equation}
    \label{eq:mu-st}
    \bigl((a^{1/\alpha}\mu_1)* (b^{1/\alpha}\mu_2)\bigr)(B)\deq
    (a+b)^{1/\alpha}\mu(B)
  \end{equation}
  for any Borel $B\in\sB(\R)$, where $\mu_1,\mu_2$ are two \iid copies
  of $\mu$.  By Theorem~\ref{th:rep}, a non-trivial random
  measure is $\sas$ if and only if
  \begin{math}
    \phi(au_1,\dots,au_k)=a^\alpha\phi(u_1,\dots,u_k)
  \end{math}
  for all $u_1,\dots,u_k\in\R$, $k\geq1$, where $\alpha>0$ by
  Theorem~\ref{th:alpha}. Thus all moment measures have homogeneous
  log-characteristic functions. In particular, the first moment
  measure (also called the intensity measure)
  $\overline{\mu}=\overline{\mu}_1$ has the characteristic function of
  the form $e^{-c|u|^\alpha}$. This means that $\overline{\mu}$ is an
  $\sas$ law in $\R$. Similarly, $\overline{\mu}_k$ is a stable law in
  $\R^{k}$. Thus, $0<\alpha\leq 2$, \ie the stability properties of
  the cone of probability measures are similar to those of a linear
  space.  Theorem~\ref{th:max-alpha} applied to $\rho=\rho_x$ yields
  that the order of homogeneity of the L\'evy measure is necessarily
  smaller than $2$. 
  
  However the introduced family of characters is still not rich
  enough, since $m_n$ (or the corresponding random variables $\xi_n$)
  may $\KDU$-weakly converge to $\neutral$, whereas their norms are
  separated from zero. Thus, Condition~\textbf{(S)} does not hold.
  Moreover a weakly convergent sequence may have unbounded norms, so
  that \textbf{(W)} neither holds. The unit sphere is not compact,
  since a weak convergent sequence of probability measures does not
  necessary converges in the Kantorovich metric. Consequently, the
  results of Section~\ref{sec:support-levy-measure} does not lead to
  the conclusion that the L\'evy measure is supported by $\KK$. Thus,
  in general the LePage series representation of a $\sas$ random
  probability measure might involve a point process on the second dual
  semigroup, see Section~\ref{sec:lepage-series-second}.
  
  The sub-invariance property of the Kantorovich metric implies that
  the LePage series in $\KK$ converges absolutely for $\alpha\in(0,1)$
  and thereupon defines a $\sas$ random probability measure. A simple
  example of a $\sas$ random probability measure is obtained if we
  take $\eps_k$ in (\ref{eq:9}) to be a deterministic measure
  representing the probability distribution of a random variable $\xi$
  with $\E|\xi|=1$. The corresponding LePage series defines a random
  probability measure $\mu$ that, for a given sequence $\Gamma_k$, is
  the distribution of $\sum \Gamma_k^{-1/\alpha}\xi_k$ with \iid
  random variables $\xi_1,\xi_2,\ldots$ distributed as $\xi$. If
  $\psi(u)=\log\E e^{iu\xi}$, then the characteristic function of
  every realisation of $\mu$ is given by
  \begin{displaymath}
    \chi_u(\mu(\omega))=\exp\{\sum_{k=1}^\infty
    \psi(u\Gamma_k^{-1/\alpha}(\omega))\}\,.
  \end{displaymath}
  In particular, if $\xi$ is concentrated at a single point 1, \ie if
  $\epsilon_k=\delta_{1}$, then 
  \begin{displaymath}
    \chi_u(\mu)=\exp\{iu\sum_{k=1}^\infty
    \Gamma_k^{-1/\alpha}\}=\exp \{iu\xi_\alpha\}\,,
  \end{displaymath}
  meaning that $\mu$ is the unit mass measure concentrated at a
  realisation of $\sas$ random variable $\xi_\alpha$ in $(\R_+,+)$.
  Similarly, if $\xi=\sqrt{\pi/2}\,\nu$, where $\nu$ is standard
  normal, then $\mu$ is the unit mass measure concentrated at a
  realisation of a normal random variable with mean 0 and variance
  $\xi_{\alpha/2}\pi/2$. Note that this random variable exists if
  $\alpha<2$, so in this particular case the LePage series converges
  for $\alpha\in(0,2)$. 
  
  In contrast, a random probability measure $\mu=\delta_{\nu}$
  concentrated at the standard normal distribution $\nu$ provides an
  example of a Gaussian element in $\KK$.  Indeed, its Laplace
  transform
  \begin{displaymath}
    \E \chi_{u_1,\dots,u_k}(\delta_{\nu})
    =\E [e^{iu_1\nu}\cdots e^{iu_k\nu}]= \E e^{i(u_1+\cdots+u_k)\nu}
    =\exp\bigl\{-\frac{1}{2}(u_1+\cdots+u_k)^2\bigr\}
  \end{displaymath}
  is the exponent of a quadratic form $q(\chi_{u_1,\dots,u_k})$.
  
  In the same spirit, it is possible to consider the semigroup of
  probability measures where the addition operation corresponds to the
  maximum of independent random variables rather than the sum. As a
  distance, one can take any probability metric that is an ideal with
  respect to the maximum operation, see \cite{zol}. 
\end{example}

\begin{example}[Uniformly integrable probability measures on $\R_+$]
  \label{ex:unif-integr}
  Consider the cone $\KK$ which is a sub-family of integrable
  probability measures $m$ on $\R_+$ such that $\KK$ is closed with
  respect to the cone operations from Example~\ref{ex:rprobm}, $\KK$
  is closed in the Kantorovich metric, and such that the unit sphere
  $\SS$ constitutes a uniformly integrable family of probability
  measures with expectation (and Kantorovich norm) being $1$.  By the
  well-known sufficient condition \cite[Lemma~II.6.3]{shir84}, this is
  the case if
  \begin{displaymath}
    \sup_{m\in\SS} \int G(x)\,m(dx) <\infty\,,
  \end{displaymath}
  for some non-negative increasing function $G$ such that $G(x)/x\ti$
  as $x\ti$. For instance, it suffices to assume that for all random
  variables $\xi$ corresponding to $m\in\KK$ one has $\E
  \xi^{1+\epsilon}\leq C(\E\xi)^{1+\epsilon}$ for some fixed constants
  $C,\epsilon>0$ and all $\xi$.  The uniform integrability condition
  ensures that the unit sphere is compact.
  
  Recall that $\neutral=\zero=\delta_0$. Assume that the involution is
  identical. The cone $\KDU$ of continuous separating characters
  consists of the $[0,1]$-valued characters given by
  \begin{equation}
    \label{eq:char-var}
    \chi_{c,u_1,\dots,u_k}(m)=e^{-c\|m\|}
    \prod_{i=1}^k \int e^{-u_ix} m(dx)\,, 
  \end{equation}
  where $c,u_1,\dots,u_k\geq 0$ and $k\geq1$.  Then the $\KDU$-weak
  convergence is identical to the metric convergence. Furthermore,
  \textbf{(S)} holds, so that every $\sas$ random probability measure
  from $\KK$ admits the LePage representation by
  Theorem~\ref{thr:alpha01}\textsl{(i)}.
\end{example}

\subsection{Cones violating basic assumptions}
\label{sec:exampl-viol-basic}

The examples below clarify the influence of the crucial assumptions on
$\KK$ such as the existence of the origin, the norm, its
sub-invariance property or the existence of a family of separating
characters satisfying \textbf{(C)}. In some of these cases the LePage
series converges, but its convergence usually has to be confirmed by
means of methods specific to the particular situation.

\begin{example}[Real line with maximum operation]
  \label{ex:R-max}
  Consider the extended line $[-\infty,\infty]$ with the maximum
  operation, conventional multiplication by non-negative numbers and
  the Euclidean metric. The neutral element is $-\infty$ and the
  origin is $0$. The corresponding norm is not sub-invariant, \eg
  $d(\max\{-1,1\},-1)=2$ is not smaller than $d(1,0)=1$. However a
  direct argument shows that the LePage series with $\alpha<0$ and
  $\eps_k=-1$ defines a max-stable random element with the Weibull
  distribution. The same holds for $\alpha>0$ and $\eps_k=1$.
  Semicontinuous characters are indicators $\one_{[-\infty,a)}$ and
  $\one_{[-\infty,a]}$. It should be noted that the Laplace transform
  of a non-negative random variable $\xi$ vanishes for characters with
  $a<0$, so that such $\xi$ does not have a proper distribution.
\end{example}

\begin{example}[Cylinder]
  Let $\KK=\R_+\times[0,2\pi)$ with the Euclidean topology. The
  addition is defined coordinatewise with the second coordinates added
  modulo $2\pi$.  The multiplication defined as
  $D_a(x_1,x_2)=(ax_1,x_2)$ acts only on the first coordinate. Such a
  convex cone does not possess the origin, since $D_a(x_1,x_2)\to
  (0,x_2)$, \ie the limit as $a\downarrow0$ is not unique. Because of
  non-uniqueness of the origin, it is unclear how to define the unit
  sphere in this example.  A natural replacement for the unit sphere
  is the set given by
  \begin{displaymath}
    \{(x_1,x_2):\; x_1=1\}=\{x:\; d(x,\lim_{c\downarrow 0}cx)=1\}
  \end{displaymath}
  for a (non sub-invariant) metric
  $d((x_1,x_2),(y_1,y_2))=|x_1-x_2|+|y_1-y_2|\mod 2\pi$. However, the
  corresponding LePage series does not converge unless the spectral
  measure is concentrated on a singleton, since adding of summands of
  type $\Gamma_k^{-1/\alpha}\eps_k$ changes the sum considerably even
  for a large $k$ because the second coordinate is not scaled.
  
  In contrast, if we look at the pair $(x_1,x_2)$ as the polar
  coordinates of a point in $\R^2$ with the Euclidean distance metric,
  then all the points of the type $(0,x_2)$ are equivalent and this
  equivalence class can be used as the origin. All classical results
  are now valid for this cone showing the principal role uniqueness of
  the origin plays in our considerations.
\end{example}

\begin{example}[Non-compact unit sphere]
  \label{ex:non-complete}
  Take $\KK$ to be the positive quadrant $\R_+^2$ with both axes
  $\{0\}\times (0,+\infty)$ and $(0,+\infty)\times \{0\}$ excluded.
  For some $\alpha\in(0,1)$, consider two independent $\alpha$-stable
  positive random variables $\xi_1$ and $\xi_2$. Then the vector
  $(\xi_1,\xi_2)$ is $\sas$ in $\KK$, whereas its spectral measure is
  concentrated on two points $(1,0)$ and $(0,1)$ which are not part of
  $\KK$. Consequently, no LePage representation is possible for
  $(\xi_1,\xi_2)$. This example shows importance of compactness of the
  unit sphere even when all other sufficient conditions for existence
  of the LePage series are met, \seg Theorem~\ref{thr:alpha01}.
  \end{example}

\begin{example}[Half-line with harmonic mean as addition]
  \label{ex:harmon}
  Let $\KK$ be the extended positive half-line $[0,\infty]$ with
  addition $x\oplus y=(x^{-1}+y^{-1})^{-1}$, usual multiplication and
  norm.  Then $\zero$ is the conventional $0$, while the neutral
  element is $\infty$.  Although the norm is sub-linear, the
  corresponding metric $d$ cannot be sub-invariant.  Indeed, by the
  required continuity of operations, $x\oplus ch\to x\oplus0$ as
  $c\downarrow 0$, so that the sub-linearity would imply $d(x\oplus
  ch,x)\leq c\|h\|\to 0$. However $d(x\oplus ch,x)\to d(x\oplus
  0,x)=d(0,x)=\|x\|$, since $x\oplus0=0$.
  
  Take the Euclidean distance as the metric. The LePage series is
  $\bigl(\sum_{k=1}^\infty \Gamma_k^{1/\alpha}\bigr)^{-1}$, thus for
  positive $\alpha$ it gives identical 0 while for $-1<\alpha<0$ it
  converges to $\xi_{-\alpha}^{-1}$, where $\xi_{-\alpha}$ is the
  $(-\alpha)$-stable random variable in $(\R_+,+)$.  Continuous
  characters are given by $\chi_a(x)=e^{-a/x}$ for $a>0$.  The bijection
  $x\mapsto x^{-1}$ provides a homomorphism of semigroups $(\R_+,
  \oplus)$ and $(\R_+,+)$ that enables one to show directly that
  $\alpha$ takes values from $(-1,0)$ and that any $\sas$ random
  variable with $\alpha\in(-1,0)$ admits the LePage representation.
\end{example}

\begin{example}[Half-plane]
  \label{ex:R-max-plus}
  Let $\KK=\R\times\R_+$ be the upper half-plane. The addition is
  defined as the arithmetic addition of the first coordinates and as
  the maximum of the second coordinates of the points. Condition
  \textbf{(C)} holds, since the characters
  $\chi_{a,t}(x)=\one_{(-\infty,a]}(x_2)\,e^{itx_1}$, $x=(x_1,x_2)\in\KK$,
  have semicontinuous modulus and continuous argument. The $\sas$ laws
  in this case have the first coordinate which is stable in the
  conventional sense in $\R$, and the max-stable second coordinate.
  Thus, the parameter of such a stable law belongs to $(0,2]$.
  
  If the cone operation is altered, so that instead of the maximum we
  take the minimum of the second coordinates, then the second
  coordinate should have a negative parameter $\alpha$ unless the
  second coordinate is identical 0. The only non-trivial stable laws
  in this case are degenerate in the second coordinate and are stable
  with $\alpha\in(0,2]$ in the first coordinate.  The Euclidean metric
  is not sub-invariant for this addition operation. Indeed, adding a
  small $h=(h_1,h_2)$ to $x=(x_1,x_2)$ results in
  $(x_1+h_1,\min(x_2,h_2))$, which may be quite far away from $x$ if
  $x_2$ is large and $h_2$ is small.
\end{example}

\begin{example}[Compact sets with Minkowski addition]
  \label{ex:comp-non-conv}
  Let us drop the convexity requirement in Example~\ref{ex:mconv}, so
  that $\KK$ is the family $\sK$ of nonempty compact sets $K$ in
  $\R^d$ with the Hausdorff metric.  Then (\ref{eq:cone4}) does not
  hold and $\KK$ does not possess a family of separating characters.
  Indeed, $K+K=L+L$ and $K+K+K=L+L+L$ is possible, for instance if
  $K=[0,1]$ and $L=[0,0.4]\cup[0.6,1]$ on the real line.  As a result,
  this example cannot be investigated using the harmonic analysis
  tools. However the Hausdorff metric is sub-invariant and so the
  LePage series still defines a $\sas$ random compact set in $\R^d$
  with $\alpha\in(0,1)$, thereby positively answering a question about
  a definition of a stochastic integral that defines non-convex stable
  random sets, see \cite[p.~457]{gin:hah85b}.  Theorem~\ref{th:2}
  complements a similar result obtained in \cite{gin:hah85b} for the
  convex case.
  
  As an example, let $\eps_k=\{0,1\}$, $k\geq1$, be deterministic.
  Then the LePage series (\ref{eq:9}) defines a $\sas$ random compact
  set $X$ in $[0,\infty)$ given by all sums of the type
  $\sum_{k=1}^\infty \Gamma_k^{-1/\alpha}i_k$, where $i_k$ is either
  $0$ or $1$, i.e. all sub-series sums of $\sum_{k=1}^\infty
  \Gamma_k^{-1/\alpha}$.
\end{example}

\begin{example}[Union-stable random closed sets]
  \label{ex:un-stable-closed}
  Let $\KK$ be the family $\sF$ of closed subsets $F\subset\R^d$ with
  the union operation as addition and homothety as multiplication by
  scalars.  Then $\neutral=\emptyset$, whereas no unique origin
  exists, since $aF$ may have a limit as $a\downarrow0$ that varies
  with the choice of $F$, \eg it can be $F$ itself if $F$ is a cone in
  $\R^d$. Any cone is an element of $\KK(\alpha)$ for each
  $\alpha\neq0$, so that a sub-invariant norm cannot be finite
  everywhere by Lemma~\ref{lemma:sub-invar}. 
  
  An analogue of the LePage series in this case is constructed by
  means of a locally finite homogeneous measure $\Lambda$ on $\sF$
  that defines a Poisson process on $\sF$. The union of the sets from
  this Poisson process defines a $\sas$ random element in $\KK$, \ie a
  union-stable random closed set. The union-stable random closed sets
  with possibly non-proper distributions have been extensively studied
  in \cite[Ch.~4]{mo1} and \cite{mo93l}.  Note that $\KK$ possesses
  the same family of separating characters as described in
  Example~\ref{ex:usrs}, so that the characterisation results from
  Section~\ref{sec:distr-stable-rand} (that do not require the
  existence of the norm) hold. Here $\sas$ random elements (or
  union-stable random closed sets) exist for every $\alpha\neq0$.
\end{example}

\begin{example}[Group with Gaussian distribution being $\sas$ of any
  $\alpha$] 
  \label{ex:new-group}
  Let $\KK$ be the family of continuous functions on $\R$ with the
  pointwise addition and the multiplication by scalars corresponding
  to the rescaling of the argument as in
  Example~\ref{ex:rescaled-functins}. Consider the characters given by
  (\ref{eq:char-f}). If $\eta$ is a standard Gaussian random variable,
  then the Laplace transform of the random function $\xi(x)=\eta
  x^\beta$ for some $\beta>0$ is given by
  $\E\chi_\nu(\xi)=e^{-q(\nu)}$, where
  \begin{displaymath}
    q(\nu)=\frac{1}{2}\Big(\int x^\beta d\nu\Big)^2
  \end{displaymath}
  is a quadratic form, since
  $q(\nu_1+\nu_2)+q(\nu_1-\nu_2)=2q(\nu_1)+2q(\nu_2)$.  Thus $\xi$ is
  a Gaussian element in $\KK$. 
  
  Let $\alpha=2\beta$. Since $(D_{a^{1/\alpha}}\xi)(x)=\eta
  a^{1/2} x$, it is easily seen that $\xi$ satisfies
  (\ref{eq:a1alph-a+b1-}), \ie is $\sas$. By varying $\beta$ it is
  possible to obtain $\sas$ Gaussian elements in such a cone with an
  arbitrary characteristic exponent $\alpha$. It should be noted
  however that the considered cone does not possess the origin and
  norm.
\end{example}

\begin{example}[Intrinsically stable random measures]
  \label{ex:isrm}
  Consider the family $\KK$ of all locally finite measures on
  $\R^d\setminus\{0\}$ with the vague topology and the conventional
  addition operation as in Example~\ref{ex:mu1}, but with the
  multiplication defined as $(D_am)(K)=m(a^{-1}K)$. The second
  distributivity law does not hold in this case, since $(m+m)$ is
  generally not $D_2m$. The neutral element is zero measure, while a
  unique origin does not exist, since the rescaled measures $D_am$ may
  have various limits as $a\to0$. A random measure $\mu$ is called
  intrinsically stable if
  \begin{equation}
    \label{eq:mu2}
    \mu_1(a^{-1/\alpha}K)+\mu_2(b^{-1/\alpha}K)\deq
    \mu((a+b)^{-1/\alpha}K)
  \end{equation}
  for all measurable $K$. Note that this definition combines the
  addition operation as in Example~\ref{ex:lfrm} with the scaling used
  in Example~\ref{ex:rprobm}. 
  
  The continuous $(0,1]$-valued characters $\chi_u$ are given by
  (\ref{eq:muca}) for $u$ being a continuous function with bounded
  support. The multiplication by scalars is uplifted to the characters
  as $(c\circ\chi_{u(\cdot)})=\chi_{u(c^{-1}\cdot)}$. Since
  \textbf{(C)} holds, an intrinsically stable non-trivial random
  measure has a homogeneous Laplace exponent $\phi(u)=\phi(\chi_u)$,
  \ie $\phi(u(s^{-1}\cdot))= s^\alpha\phi(u(\cdot))$ for all
  continuous $u$ with bounded support and $s>0$, \cf
  Theorem~\ref{th:rep}.  Theorem~\ref{th:alpha} is not applicable to
  show that $\alpha$ has a particular sign, since the origin is not
  defined.
  
  In this case, any $\alpha\neq0$ is a possible parameter for a $\sas$
  random measure. For instance, the counting measures $\Pi_\alpha$ for
  all $\alpha\neq0$ are stable in this cone $\KK$.
\end{example}

\begin{example}[Random probability measures]
  \label{ex:all-prob-m}
  Let $\KK$ be the family of all probability measures on $\R$ with the
  same cone operations as in Example~\ref{ex:rprobm}, convergence in
  distribution and the characters given by~(\ref{eq:char-uuu}).
  Theorem~\ref{th:alpha} implies that $\alpha>0$.  In this case
  $\KK(\alpha)$ are all non-trivial for $0<\alpha\leq 2$ consisting of
  $\sas$ probability distributions in $\R$. As a consequence of
  Lemma~\ref{lemma:sub-invar}, it is not possible to define a
  sub-invariant metric in this cone unless assigning infinite norm for
  all elements of $\KK(\alpha)$ with $\alpha<1$. Recalling the
  discussion in Example~\ref{ex:rprobm}, this shows that it is not
  possible to construct an ideal probability metric of order $1$ that
  is finite on all random variables.
  
  Convergence of the LePage series should be checked in each
  particular case. For this, it is simpler to work with random
  characteristic functions rather than with random probability
  measures, since the convolutions of measures become the product of
  characteristic functions.  Let $\KK$ be the family of characteristic
  functions $\theta(t)$, $t\in\R$, with addition being the product of
  functions and multiplication by numbers corresponding to the
  rescaling of the argument, \ie $(D_a\theta)(t)=\theta(at)$. The
  LePage series then yields the random characteristic function
  \begin{equation}
    \label{eq:prod-theta}
    \prod_{k=1}^\infty \theta_k(\Gamma_k^{-1/\alpha}t)\,,
  \end{equation}
  where $\{\theta_k,\, k\geq1\}$ are \iid random characteristic
  functions from a unit sphere in $\KK$.  The unit sphere $\SS$ in
  $\KK$ can defined by using any polar decomposition of $\KK$.  For
  instance, a homogeneous (but not sub-linear) norm of a probability
  measure $m$ with distribution function $F_m(x)=m\{(-\infty,x]\}$ and
  the corresponding percentiles $q_m(t) =\inf\{x:\ F_m(x)\geq t\}$ can
  be defined as 
  \begin{displaymath}
    \|m\|^2=q_m^2\bigl(F_m(0-)/2\bigr) +
    q_m^2\bigl((1+F_m(0))/2\bigr)\,.
  \end{displaymath}
  Notice that $\|m\|=0$ implies~$m=\delta_0$.  In case of integrable
  centred probability measures, one can use the sub-invariant
  homogeneous (\ie ideal) metric given by
  \begin{displaymath}
    d(\xi,\eta)=d(\theta_1,\theta_2)
    =\sup_{u\in\R} |u|^{-1}|\theta_1(u)-\theta_2(u)|\,,
  \end{displaymath}
  where $\theta_1$ and $\theta_2$ are the characteristic functions of
  $\xi$ and $\eta$ respectively.
  
  Assume that $\{\theta_k,\,k\geq1\}$ are non-random characteristic
  functions of the strictly stable distribution with characteristic
  exponent $\beta\in(0,1]$, \ie $\theta_k(t)=e^{-|t|^\beta}$.  In case
  $\beta=1$ one obtains the symmetric Cauchy distribution.
  Convergence of the infinite product (\ref{eq:prod-theta}) is then
  equivalent to convergence of the series
  \begin{displaymath}
    \sum_{k=1}^\infty \Gamma_k^{-\beta/\alpha}\,.
  \end{displaymath}
  This series converges \as if and only if $\beta>\alpha$, and its
  limit is a strictly stable random variable $\zeta$ with
  characteristic exponent $\alpha/\beta$. Thus, the LePage series
  produces a random probability measure $D_\zeta m$ that corresponds
  to the strictly stable probability measure $m$ with characteristic
  exponent $\beta\in(0,1)$ rescaled by a random $\alpha/\beta$-stable
  coefficient $\zeta$.
  
  As another example, assume that $\{\theta_k,\,k\geq1\}$ are
  characteristic functions of the Gamma-distributions with random
  shape parameter $p_k>0$ and scale parameter $1$. The corresponding
  measures are integrable, so that this case is covered by
  Example~\ref{ex:rprobm}. Since the unit sphere consists of all
  probability measures on $\R_+$ with unit mean, the considered random
  measures (or random characteristic functions) need to be rescaled by
  $p_k$ to lie on the unit sphere.  Then~(\ref{eq:prod-theta}) becomes
  \begin{displaymath}
    \prod_{k=1}^\infty (1-i\Gamma_k^{-1/\alpha}t/p_k)^{-p_k}\,,
  \end{displaymath}
  so that the LePage series converges for every $\alpha\in(0,1)$ and
  every distribution of the $p_k$'s. This is also confirmed by the
  existence of a sub-invariant norm for integrable probability
  measures, see Example~\ref{ex:rprobm}.
\end{example}

\section*{Acknowledgements}
\setlength{\fboxsep}{0.5ex}

YuD is grateful to V.~A.~Egorov, \fbox{A.~V.~Nagaev} and V.~Tarieladze
for useful discussions and interest to this work.  SZ also
acknowledges hospitality of Institute Mittag-Leffler, the Royal
Swedish Academy of Sciences, where a part of this work has been
carried out.

\bibliography{/home/ilya/tex/bibdata/abbrev,/home/ilya/tex/bibdata/x/gesamt}

\begin{center}
  \parbox{7.5cm}{
    Youri Davydov\\
    Universit\'e  Lille~1\\
    Laboratoire de Paule Painlev\'e\\
    UFR de Math\'ematiques\\
    Villeneuve d'Ascq Cedex\\
    France\\
    E-mail: youri.davydov@univ-lille1.fr}
  \parbox{7.5cm}{
    Ilya Molchanov\\
    University of Berne\\
    Department of Mathematical
    Statistics\\ and Actuarial Science\\
    Sidlerstr. 5\\
    CH-3012 Berne\\
    Switzerland\\
    E-mail: ilya@stat.unibe.ch}

  \bigskip
  \parbox{8cm}{
    Sergei Zuyev\\
    Department of Statistics \\ and Modelling Science\\
    University of Strathclyde\\
    Glasgow G1~1XH\\
    Scotland, U.~K.\\
    E-mail: sergei@stams.strath.ac.uk}  
\end{center}

\end{document}